\documentclass[a4paper,11pt]{amsart}

\usepackage{amssymb,xspace,amscd,yhmath,mathdots,txfonts,
mathrsfs}
\usepackage[matrix,arrow,curve]{xy}\CompileMatrices

\usepackage{hyperref}

\DeclareMathOperator{\rank}{rank}
\DeclareMathOperator{\kt}{{\kappaup}}

\DeclareMathOperator{\Li}{\mathfrak{L}}
\DeclareMathOperator{\Lio}{\mathfrak{L}_0}
\DeclareMathOperator{\nil}{\mathrm{nil}}
\DeclareMathOperator{\rad}{\mathrm{rad}}
\DeclareMathOperator{\radn}{\mathrm{rad}_n}
\DeclareMathOperator{\rn}{\mathfrak{n}_{\mathfrak{g}}}
\DeclareMathOperator{\nr}{\mathfrak{n}_{{\kt}}}

\DeclareMathOperator{\ad}{\mathrm{ad}}
\DeclareMathOperator{\Ad}{\mathrm{Ad}}
\DeclareMathOperator{\Lie}{\mathrm{Lie}}
\DeclareMathOperator{\lie}{\mathscr{L}}
\DeclareMathOperator{\Par}{{\mathbb{P}_{\kt}}}

\DeclareMathOperator{\ses}{\mathfrak{s}}

\DeclareMathOperator{\zt}{\mathfrak{z}}
\DeclareMathOperator{\st}{\mathfrak{s}}
\DeclareMathOperator{\zto}{\mathfrak{z}_0}
\DeclareMathOperator{\N}{\mathrm{N}}
\DeclareMathOperator{\Z}{\mathrm{Z}}

\DeclareMathOperator{\pt}{\mathfrak{p}}
\DeclareMathOperator{\ft}{\mathfrak{f}}
\DeclareMathOperator{\Ind}{\mathrm{Int}}
\DeclareMathOperator{\Zt}{\mathscr{Z}}
\DeclareMathOperator{\It}{\mathscr{I}}
\DeclareMathOperator{\Vt}{\mathscr{V}}
\DeclareMathOperator{\Ot}{\mathscr{O}}
\DeclareMathOperator{\Dt}{\mathscr{D}}

\title[Reductive compact homogeneous CR~manifolds]
{Reductive compact homogeneous CR~manifolds}
\author[A.~Altomani]{Andrea Altomani}
\address{A.\ Altomani:
University of Luxembourg \\ Campus Kirchberg
Mathematics Research Unit, BLG \\
6, rue Richard Coudenhove-Kalergi
L-1359 Luxembourg
Grand-Duchy of Luxembourg}
\email{andrea.altomani@uni.lu}

\author[C.~Medori]{Costantino Medori}
\address{C.\ Medori:
Dipartimento di Matematica\\ Universit\`a di Parma\\ Viale G.P.
Usberti, 53/A
\\ 43100 Parma (Italy)} \email{costantino.medori@unipr.it}

\author[M.~Nacinovich]{Mauro Nacinovich}
\address{M.\ Nacinovich:
Dipartimento di Matematica\\ II Universit\`a di Roma
``Tor Ver\-ga\-ta''\\ Via della Ricerca Scientifica\\ 00133 Roma
(Italy)}
\email{nacinovi@mat.uniroma2.it}
\date{\today}

\subjclass[2000]{Primary: 32V05
Secondary: 32L05, 53C30}
\keywords{Compact homogeneous $CR$ manifold,
$CR$ algebra, equivariant fibration}

\numberwithin{equation}{section}

\theoremstyle{plain}
\newtheorem{thm}{Theorem}[section]
\newtheorem{lem}[thm]{Lemma}
\newtheorem{cor}[thm]{Corollary}
\newtheorem{prop}[thm]{Proposition}

\theoremstyle{definition}
\newtheorem{exam}[thm]{Example}

\newtheorem{dfn}[thm]{Definition}
\newtheorem{rmk}[thm]{Remark}
\newtheorem{defn}[thm]{Definition}
\newtheorem{ntz}[thm]{Notation}

\begin{document}
\begin{abstract}
We consider a class of compact homogeneous $CR$ manifolds, 
that we call $\mathfrak n$-reductive, which includes the orbits 
of minimal dimension of a compact Lie group $\mathbf{K}_0$ 
in an algebraic homogeneous variety of its complexification 
$\mathbf{K}$. For these manifolds we define 
canonical equivariant fibrations onto complex flag manifolds. 
The simplest example is the Hopf fibration $S^3\to\mathbb{CP}^1$.
\par
\noindent
In general these fibrations are not $CR$ submersions, 
however they satisfy a weaker condition that we introduce here, 
namely they are $CR$-deployments.
\end{abstract}
\maketitle
\tableofcontents
\section{Introduction}
In this article we focus our attention 
on a class of compact homogeneous $CR$ manifolds. Our objects 
will turn out to be
minimal orbits of a compact Lie group
in their equivariant canonical 
$CR$-embedding (see \S\,\ref{sec46}). We call them 
here $\mathfrak{n}$-\textit{reductive}
(Definition\,\ref{df33}), because the complex Lie algebra 
describing 
their homogeneous $CR$ structure 
is the direct sum of a nilpotent ideal and of the complexification of
their isotropy.
Compact homogeneous
$CR$ manifolds and associated canonical bundles were 
also considered in \cite{AS02, AS03, GH09}. 
We shall also consider various canonical
fibrations. Let us consider a simple example.\par
Let $Q$ be the non degenerate 
Hermitian quadric in $\mathbb{CP}^n$, with index $q\leq[\tfrac{n}{2}]$.
It is the minimal orbit of the natural action of $\mathbf{SU}_{q+1,n-q}$ in
$\mathbb{CP}^{n}$.  The choice of 
two $Q$-polar subspaces $S_1,S_2$ of $\mathbb{CP}^{n}$,
of dimension $q$ and $n\! -\! q\!-\! 1$, respectively, that do not intersect
$Q$, is equivalent to that of a maximal compact subgroup
$\mathbf{K}_0\simeq\mathbf{S}(\mathbf{U}_{q+1}\times\mathbf{U}_{n-q})$ of
$\mathbf{SU}_{q+1,n-q}$. This choice
defines a map $\piup:Q\to S_1\times{S}_2$ by
\begin{gather*}
  Q\ni p \longrightarrow \piup(p)=(\piup_1(p),\piup_2(p))\in S_1\times{S}_2,\\
\text{with $\{\piup_1(p)\}=S_1\,\cap\, pS_2$, 
$\{\piup_2(x)\}=S_2\,\cap\, pS_1$,}
\end{gather*}
which presents ${Q}$  as a circle bundle over
$\mathbb{CP}^{q}\times\mathbb{CP}^{n-q-1}$, and reduces, for
$q=0$, to the classical Hopf fibration of the odd dimensional sphere.
 The projection $\piup$ is a $CR$-submersion, being the restriction
of a holomorphic submersion from 
$V=\mathbb{CP}^{n}\setminus (S_1\cup{S}_2)$,
with  fibers biholomorphic to $\mathbb{C}^*$. We note that 
$V$ is a homogeneous space of the complexification 
$\mathbf{K}\simeq \mathbf{S}(\mathbf{GL}_{q+1}(\mathbb{C})\times\mathbf{GL}_{n-q})$
of $\mathbf{K}_0$ and that the immersion $Q\hookrightarrow{V}$ is induced by
the inclusion $\mathbf{K}_0\subset\mathbf{K}$. 
\par\smallskip
We shall generalize this construction to 
general $\mathfrak{n}$-reductive
homogeneous compact $CR$ manifolds. 
This class contains the minimal orbits 
of real forms in complex flag manifolds (see \cite{AMN06}), 
and, more in general, the
compact homogeneous $CR$ manifold arising as
intersections of general orbits with their Matsuki dual (see \cite{Mats88,
AMN06b}). \par
Our generalization of the Hopf fibration 
will yield $CR$ maps which are smooth submersions, but which may fail,
in general, to be $CR$-submersions. They satisfy, however, 
a weaker requirement, which lead us to the notions of 
$CR$-spread and
$CR$-deployment. 
A $CR$ map $\piup:M\to{N}$ is a $CR$-spread
if $N$ has the minimal $CR$ structure for which $\piup$ is a $CR$ map,
and a $CR$-deployment if, moreover, $d\piup$ is injective on the 
analytic tangent bundle. We prove that an $\mathfrak{n}$-reductive $M$
has a $CR$-deployment over a complex flag manifold by a map that
restricts a holomorphic fibrations with Stein fibers (see Theorems\,\ref{thm68},
\,\ref{thm69}). This notion has also a counterpart for foliated
complex manifolds, connected to the Segre varieties associated to
a $CR$ manifold (see Theorem\,\ref{th713}).
\par\smallskip
Let $\piup:M\to{N}$ be a smooth submersion of an $m$-dimensional smooth manifold
$M$ onto an $n$-dimensional
complex manifold $N$.  An Ehresmann connection
$\mathcal{H}$ on $M\xrightarrow{\;\piup\;}N$ 
uniquely defines an almost $CR$ structure on $M$ for which the analytic
tangent space is the horizontal distribution and $\piup$ a $CR$-submersion.
The $CR$ structure on the nondegenerate Hermitian 
quadric $Q$  is obtained
from its Hopf fibration 
$\piup$ if we define $\mathcal{H}$ to be the analytic tangent space to
${Q}$. 
\begin{defn}
Let $M$ be a $CR$ manifold of $CR$ dimension $n$ and $CR$ codimension $k$,
$N$ an $n$-dimensional complex manifold. If 
$\piup:M\to{N}$ 
a $CR$-submersion, for which the analytic tangent space $HM$ is
the horizontal distribution of an
Ehresmann connection, then we say that the $CR$ structure of $M$ 
\emph{lifts} the complex structure of~$N$. 
\end{defn}
In general, an Ehresmann connection $\mathcal{H}$ on 
$M\xrightarrow{\;\piup\;}N$
needs to satisfy strong integrability
conditions for the lifted structure to be $CR$. 
Our construction in \S\ref{lif} shows
that any $\mathfrak{n}$-reductive
$CR$ structure on $M$ can be strengthened to the lift of
a complex flag manifold. This also solves the problem of determining
all homogeneous $\mathfrak{n}$-reductive $CR$ structures on a given
smooth compact homogeneous manifold. Similar problems were considered
for instance in \cite{Charb04, AMN10a}. \par
Let us briefly describe the contents of this paper. 
We collected in \S\ref{sec:hom} and \S\ref{sechm} the preliminary notions on
$CR$ manifolds and homogeneous $CR$ manifolds. In \S\ref{sc4} we consider
specifically those which are homogeneous for the action of a compact
Lie group, and introduce the notion of \textit{$\mathfrak{n}$-reductive}.
In \S\ref{sec:mf} we show that the intersection of the orbits of the real
forms in complex flag manifolds with their Matsuki duals are
$\mathfrak{n}$-reductive. In \S\ref{sec:rg} we recall a classical construction
of Chevalley (see e.g. \cite{BT71}) which is preliminary to the construction
of the deployments and the lifts of \S\ref{sec:dp} and  \S\ref{lif}.
In \S\ref{gen} we make some remarks for the case of compact
homogeneous $CR$ manifolds which are not $\mathfrak{n}$-reductive. In the final
\S\ref{exa} we discuss some examples.
\section{$CR$ manifolds and $CR$ maps}
\label{sec:hom}
Let $M$ be a smooth manifold of real dimension $m$, and 
$n_{\muup},k_{\muup}\geq{0}$ two
integers with $2n_{\muup}+k_{\muup}=m$.
A \emph{$CR$ structure} $\muup$ of type 
$(n_{\muup},k_{\muup})$ on $M$ is the datum of
a  rank $n_{\muup}$ smooth complex vector 
subbundle $T^{0,1}M$ of its complexified tangent bundle
${\mathbb{C}}TM$, with 
\begin{equation}
  \label{eq:h0} T^{0,1}M\cap\overline{T^{0,1}M}=\underline{0},\quad 
  [\Gamma(M,T^{0,1}M),\Gamma(M,T^{0,1}M)]\subset\Gamma(M,T^{0,1}M).
\end{equation}
The integer $n_{\muup}$ is the $CR$-dimension, and 
$k_{\muup}=\mathrm{dim}_{\mathbb{R}}{M}-2n_{\muup}$ the
$CR$-codimension of $\muup=(M,T^{0,1}M)$.
If $n_{\muup}=0$, we say that $\muup$ is totally real; if 
$k_{\muup}=0$, $\muup$ is
a complex structure on $M$,
in view of the Newlander-Nirenberg theorem.\par
Let $\muup=({M},T^{1,0}{M})$ and $\nuup=({N},T^{1,0}{N})$ be two
$CR$ manifolds.\par
A smooth map $\phiup:{M}\to{N}$ is $CR$
if $d\phiup(T^{0,1}{M})\subset{T}^{0,1}{N}$.
A \emph{$CR$-immersion} (\emph{$CR$-embedding}) is a smooth immersion 
(embedding) \mbox{$\phiup:{M}\to{N}$}
such that
\begin{equation}\label{eq:h3}
  d\phiup(T^{1,0}{M})=T^{1,0}{N}\cap
d\phiup(T^{\mathbb{C}}{M}).
\end{equation}
A $CR$-immersion is called \emph{generic} if, moreover, 
$d\phiup({\mathbb{C}}T_p{M})+T^{0,1}_{\phiup(p)}N
=T^{\mathbb{C}}_{\phiup(p)}N$ for all $p\in{M}$. 
This is the case when $n_{\muup}+k_{\muup}=n_{\nuup}+k_{\nuup}$.
\par
A smooth submanifold $M$ of a $CR$ manifold $N$ is a $CR$-submanifold if 
$T^{1,0}M=T^{\mathbb{C}}M\cap{T}^{1,0}N$ 
defines a $CR$ structure on $M$, i.e. if $T^{1,0}M$
is a complex distribution
of constant rank on~$M$. \par
A \emph{$CR$-submersion}, or \emph{$CR$-bundle} is a smooth submersion 
$\piup:{M}\to{{N}}$ such that $d\piup(T_{p}^{0,1}M)=T^{0,1}_{\piup(p)}{N}$
for all $p\in{M}$. In this case the fibers are $CR$-submanifolds
of ${M}$. They are totally real if and only if ${M}$ and $N$ have the
same $CR$-dimension,
and totally complex if and only if ${M}$ and $N$ have the
same $CR$-codimension. \par
A $CR$-diffeomorphism is a smooth diffeomorphism $\phiup:M\to{N}$ which
is at the same time a $CR$-immersion and a $CR$-submersion.\par
We recall the concept of strengthening a $CR$ structure
(see \cite[Def.\,5.8]{AMN06b}).
\begin{defn}
Let $M$ be a smooth manifold and $\muup_1,\muup_2$ two $CR$ structures on
$M$. If the identity map from $\muup_1$ to $\muup_2$ is $CR$,
we say that $\muup_2$ is a \emph{strengthening} of $\muup_1$
or, equivalently, that $\muup_1$ is a \emph{weakening} of $\muup_2$.
\end{defn}

\par\smallskip
Let $M$ be a smooth manifold and denote by 
$\varOmega_M$ the sheaf of germs of complex valued smooth exterior forms
on $M$ and by $\mathfrak{X}^{\mathbb{C}}_M$ the sheaf of germs of complex valued 
smooth vector fields on $M$. 
\begin{defn}[ideal sheaf and characteristic distribution]
The \emph{ideal sheaf $\It_{\muup}$} of 
$\muup=(M,T^{0,1}M)$ is
the subsheaf of $\varOmega_M$ consisting of the germs of forms of
positive degree that vanish on $T^{0,1}M$. 
It is graded and $d$-closed.
The \textit{characteristic distribution} ${\Zt}_{\muup}$ 
of ${\muup}$ is the sheaf
of germs of smooth sections of $T^{0,1}M$.
\end{defn}
If $\muup_1,\muup_2$ are $CR$ structures on the same smooth manifold $M$, 
necessary and sufficient conditions for $\muup_2$ to be a 
strengthening of $\muup_1$ is that 
${\Zt}_{\muup_1}\subset{\Zt}_{\muup_2}$, 
or, equivalently, that
${\It}_{\muup_2}\subset{\It}_{\muup_1}$.
\par\smallskip
Let $M$ be a smooth manifold, $\nuup=(N,T^{0,1}N)$ a smooth $CR$ manifold
and $\phiup:M\to{N}$ a smooth map. The pullback 
$\phiup^*{\Zt}_{\nuup}$ of ${\Zt}_{\nuup}$
is the sheaf associated to the presheaf
 \begin{equation}\label{eq:h3b}
 U\longrightarrow
 \phiup^*{\Zt}_{\nuup}(U)=\{X\in\mathfrak{X}^{\mathbb{C}}_M(U)\mid 
d\phiup(X_p)\in{T}^{0,1}N,\;
 \forall p\in{U}\}.
\end{equation}
Denote by $\phiup^*{\It}_{\nuup}$, the 
$\mathcal{C}^{\infty}_M$-module generated by the pull-backs $\phiup^*\omega$
of germs $\omega\in{\It}_{\nuup}$.
The elements of $\phiup^*{\It}_{\nuup}$ vanish on
$\phiup^*{\Zt}_{\nuup}$ and
the system  $\phiup^*{\Zt}_{\nuup}$ is \textsl{formally integrable}. 
Indeed, we have 
\begin{equation}
 [\phiup^*{\Zt}_{\nuup},\phiup^*{\Zt}_{\nuup}]\subset
 \phiup^*{\Zt}_{\nuup},\qquad d(\phiup^*{\It}_{\nuup})
\subset\phiup^*{\It}_{\nuup}.
\end{equation}
One easily verifies the following
\begin{lem} Let 
$\muup=(M,T^{0,1}M)$ and $\nuup=(N,T^{0,1}N)$ be $CR$ manifolds, 
$\phiup:M\to{N}$ a smooth map and $p_0\in{M}$. The following are equivalent 
\begin{enumerate}
 \item $\phiup$ is $CR$ on a neighborhood of $p_0$;
 \item $\phiup^*{\It}_{\nuup_{\phiup(p_0)}}\subset{\It}_{\muup_{(p_0)}}$;
 \item ${\Zt}_{\muup_{(p_0)}}\subset\phiup^*{\Zt}_{\nuup_{\phiup(p_0)}}$.
\end{enumerate}
\end{lem}
If $\phiup:M\to{N}$ is a smooth map, we denote by 
${\Vt}_{\phiup}$
the sheaf of germs
of $\phiup$-vertical smooth complex vector fields on $M$. This is 
the sheaf associated to the
presheaf 
\begin{equation*}
 U\longrightarrow {\Vt}_{\phiup}(U)=\{X\in\mathfrak{X}^{\mathbb{C}}_M(U)\mid
 d\phiup(X_p)=0,\;\;\forall p\in{U}\}.
\end{equation*}
If $\nuup=(N,T^{0,1}N)$ is a $CR$ structure on $N$, then 
${\Vt}_{\phiup}\subset\phiup^*{\Zt}_{\nuup}$, 
and if $\phiup$ is $CR$ for 
$\muup=(M,T^{0,1}M)$, then $\phiup^*{\Zt}_{\nuup}$
contains the Lie subalgebra of $\mathfrak{X}^{\mathbb{C}}_M$ generated by
${\Zt}_{\muup}+{\Vt}_{\phiup}$.\par
It is therefore natural to introduce the following generalization of the notion of 
$CR$-submersion.
\begin{defn}[$CR$-spread, $CR$-deployment] \label{df24} 
Let $\muup=(M,T^{0,1}M)$, $\nuup=(N,T^{0,1}N)$ be $CR$ manifolds and 
$\phiup:M\to{N}$ a smooth map. \par\smallskip
We say that $\phiup$ is a \emph{$CR$-spread}
at $p_0\in{M}$ if it is a smooth submersion at $p_0$ and, moreover
\begin{equation}\label{eq:h4}
 \phiup^*{\Zt}_{\nuup_{(p_0)}} \;\;\text{is the Lie algebra generated by\;
$ {\Zt}_{\muup_{(p_0)}}\!\!+{\Vt}_{\phiup_{(p_0)}}$}.
\end{equation}
\par
We say simply that $\phiup$ is a \emph{$CR$-spread} if 
$\phiup:M\to{N}$ is a smooth submersion and 
\begin{equation}\label{eq:h4b}
 \phiup^*{\Zt}_{\nuup} \;\;\text{is the Lie algebra generated by
$ {\Zt}_{\muup}\!+{\Vt}_{\phiup}$}.
\end{equation}
\par
We say that $\phiup$ is a \emph{$CR$-deployment} at
$p_0\in{M}$ if it is a 
$CR$ spread at $p_0$ and moreover $d\phiup:T^{0,1}_{p_0}M\to T^{0,1}_{q_0}N$
is injective.
\par
We say simply that $\phiup$ is a \emph{$CR$-deployment} if it is
a $CR$ deployment at all points. \par
Equivalently, a $CR$-deployment is a $CR$-spread with totally real fibers.
\end{defn}
Roughly speaking, a $CR$ map $\phiup:M\to{N}$ is a $CR$-spread when $\phiup$ is
a submersion and $N$ carries the \textsl{weakest} $CR$ structure for which 
$\phiup$ is $CR$.\par
\smallskip
The $CR$ spread has a nice interpretation in terms of the sheaves 
\begin{equation*}
{\Ot}_{\muup}=\{f\in\mathcal{C}^{\infty}_M\mid df\in{\It}_{\muup}\},\quad
{\Ot}_{\nuup}=\{f\in\mathcal{C}^{\infty}_N\mid df\in{\It}_{\nuup}\}
\end{equation*}
of germs of
smooth $CR$ functions on $M,N$. 
\begin{prop}
Let $M,N$ be $CR$ manifolds, with $CR$ structures $\muup$, $\nuup$, respectively,
 and $\phiup:M\to{N}$ a smooth submersion. \par
\begin{enumerate}
\item $\phiup$ is a $CR$-spread at $p_0\in{M}$ if and only if
  \begin{equation}
    \label{eq:h77}
    \omega\in{\It}_{\nu_{\phiup(p_0)}} 
\Longleftrightarrow \phiup^*\omega\in{\It}_{\mu_{p_0}}.
  \end{equation}
\item In particular, 
 if $\phiup$ is a 
 $CR$-spread at
$p_0\in{M}$, then
\begin{equation}\label{eq:h71}
  {\Ot}_{N_{\phiup(p_0)}}=\{u\in\mathcal{C}^{\infty}_{N_{\phiup(p_0)}}\mid
[\phi^*u]_{(p_0)}\in{\Ot}_{M_{(p_0)}}\}.
\end{equation}\\
Condition
\eqref{eq:h71} is also sufficient when $N$ is locally embeddable at 
\mbox{$\phiup(p_0)$}.
\end{enumerate}
\end{prop}
We recall that $N$ is locally embeddable at $q_0\in{N}$ 
if the real parts of germs of smooth
$CR$ functions provide coordinate patches at $q_0$. This is always the case when
the $CR$ structure $\nuup$ of $N$ is real-analytic.
\section{Homogeneous $CR$ manifolds and their $CR$ algebras}
\label{sechm}
\subsection{Homogeneous $CR$ manifolds}
Let $\mathbf{G}_0$ be a Lie group, $\mathfrak{g}_0$ its Lie algebra and
$\mathfrak{g}=\mathbb{C}\otimes_{\mathbb{R}}\mathfrak{g}_0$ its complexification. 
\begin{dfn} \label{dfn:aa}
A ${{\mathbf{G}_0}}$-homogeneous $CR$ manifold is a
${{\mathbf{G}_0}}$-homogeneous space, endowed with a 
${{\mathbf{G}_0}}$-equivariant
$CR$ structure.
\end{dfn}
Let $M$ be a $\mathbf{G}_0$-homogeneous $CR$ manifold. Fix $p_0\in{M}$
and denote by 
${\mathbf{M}}_0=\{g\in\mathbf{G}_0\mid g\!\cdot\! p_0=p_0\}$ 
the isotropy subgroup of $M$ at
$p_0$. Let $\piup:\mathbf{G}_0\to{M}\simeq\mathbf{G}_0/{\mathbf{M}}_0$
be the canonical projection and
$d\piup:\mathfrak{g}\to{T}_{p_0}^{\mathbb{C}}M$
the complexification of the differential of $\piup$ at the identity. 
We showed (see e.g. \cite[\S{1}]{AMN06b})
that
\begin{equation}\label{eq:h1}
  \mathfrak{q}=d\piup^{-1}({T^{0,1}_{p_0}M})
\end{equation}
is a complex Lie subalgebra of $\mathfrak{g}$, with
$\mathfrak{q}\cap\mathfrak{g}_0=\mathfrak{m}_0=\Lie(\mathbf{M}_0)$.
\begin{defn}
The pair 
$(\mathfrak{g}_0,\mathfrak{q})$ is the $CR$-algebra at $p_0$
of the $\mathbf{G}_0$-homogeneous $CR$ manifold $M$. 
\end{defn}
Clearly the datum of the $CR$ algebra at any single point of $M$ completely
defines its $CR$ structure and we have:
\begin{prop}\label{prop28}
The $\mathbf{G}_0$-homogeneous
$CR$ structures on $M$ are in one-to-one correspondence
with the set of complex Lie subalgebras
$\mathfrak{q}$ of $\mathfrak{g}$ such that
\begin{enumerate}
\item
$\mathfrak{q}\cap\mathfrak{g}_0$ equals the Lie algebra $\mathfrak{m}_0$
of $\mathbf{M}_{0}$,
\item $\mathfrak{q}$ is $\Ad_{\mathfrak{g}}(\mathbf{M}_0)$-invariant.\qed
\end{enumerate}
\end{prop}
This was our motivation to introduce and discuss 
\emph{$CR$ algebras} in \cite{MN05}. \par\smallskip
Let $M$ be a $\mathbf{G}_0$-homogeneous $CR$ manifold, with $CR$ algebra
$(\mathfrak{g}_0,\mathfrak{q})$ at $p_0\in{M}$. The complex valued
$\mathbf{G}_0$-equivariant one-forms on $M$ are in one-to-one correspondence
with the elements of the dual $\mathfrak{g}^*$ of $\mathfrak{g}$ which
annihilate $\mathfrak{q}\!\cap\!\bar{\mathfrak{q}}$. For
$\etaup\in(\mathfrak{q}\!\cap\!\bar{\mathfrak{q}})^0$, 
we define $\etaup_{p_0}^*\in{^{\mathbb{C}}T^{*}_{p_0}M}$ by the commutative 
diagram
\begin{equation*}
  \xymatrix{
    \mathfrak{g}\ar[rr]^{d\piup}\ar[dr]_{\etaup}&& T^{\mathbb{C}}_{p_0}M
\ar[dl]^{\etaup_{p_0}^*}\\
&\,\mathbb{C}},
\end{equation*}
where $\piup:\mathbf{G}_0\ni{g}\to{g}\cdot{p_0}\in{M}$ is
the natural projection.
Then $\etaup_{p_0}^*$ can be extended to a  $\mathbf{G}_0$-equivariant one-form
by setting $\etaup^*_{g\cdot{p_0}}=g_*\etaup^*_{p_0}$, for $g\in\mathbf{G}_0$.
\par
The ideal sheaf ${\It}_M$ is generated, as a $\mathcal{C}^\infty_M$ 
module, by the forms
$\etaup^*$, for $\etaup$ varying in 
the annihilator $\mathfrak{q}^0=\{\etaup\in\mathfrak{g}^*
\mid \etaup^*(Z)=0,\;\forall Z\in\mathfrak{q}\}$ of
$\mathfrak{q}$ in $\mathfrak{g}^*$.
\begin{dfn} \label{dfn:ad}
A \emph{$CR$ algebra} is a pair $(\mathfrak{g}_0,\mathfrak{q})$, consisting
of a real Lie algebra $\mathfrak{g}_0$ and of a complex Lie subalgebra
$\mathfrak{q}$ of its complexification $\mathfrak{g}$, such that
the quotient $\mathfrak{g}_0/(\mathfrak{q}\cap\mathfrak{g}_0)$ is a finite
dimensional real vector space.\par
If $M$ is a ${{\mathbf{G}_0}}$-homogeneous $CR$ manifold and $\mathfrak{q}$
is defined by
\eqref{eq:h1}, 
we say that the $CR$ algebra
$(\mathfrak{g}_0,\mathfrak{q})$ is \emph{associated} to $M$ at $p_0$.\par
\end{dfn}
\begin{rmk}\label{rmk:ab}
If $(\mathfrak{g}_0,\mathfrak{q})$ is 
associated
to 
$M$ at  $p_0$,
and $g\in\mathbf{G}_0$, then $(\mathfrak{g}_0,\Ad_{\mathfrak{g}}(g)(\mathfrak{q}))$
is associated to $M$ at $g\!\cdot\!{p}_0$. 
\end{rmk}
\begin{rmk}\label{rmk:ac}
The $CR$-dimension and $CR$-codimension of $M$ can be computed in terms
of its associated 
$CR$ algebra $(\mathfrak{g}_0,\mathfrak{q})$. We have indeed
\begin{equation}
  CR\text{-}\mathrm{dim}\,{M}\!=\!\mathrm{dim}_{\mathbb{C}}\mathfrak{q}\! -\! 
\mathrm{dim}_{\mathbb{C}}(\mathfrak{q}\cap\bar{\mathfrak{q}}), \quad
CR\text{-}\mathrm{codim}\,{M}\! =\! \mathrm{dim}_{\mathbb{C}}\mathfrak{g}\! -\!
\mathrm{dim}_{\mathbb{C}}(\mathfrak{q}+\bar{\mathfrak{q}}).  
\end{equation}
\end{rmk}
\begin{defn}\label{dfn:26}
 Let $\mathbf{G}_0$ be a real form of a complex Lie group $\mathbf{G}$.
We say that a $\mathbf{G}_0$-homogeneous $CR$ manifold $M$ is 
\emph{$\mathbf{G}$-realizable} if $M$ admits a $\mathbf{G}_0$-equivariant
generic $CR$-embedding into a complex homogeneous space
$N$ of $\mathbf{G}$.
\end{defn}
\begin{rmk} Let $M$ be a $\mathbf{G}_0$-homogeneous $CR$ manifold,
with isotropy $\mathbf{M}_0$ and $CR$ algebra
$(\mathfrak{g}_0,\mathfrak{q})$ at $p_0\in{M}$. Then
$M$ is $\mathbf{G}$-realizable if and only if there is a closed complex Lie subgroup
$\mathbf{Q}$ of $\mathbf{G}$ with $\Lie(\mathbf{Q})=\mathfrak{q}$ and
$\mathbf{Q}\cap\mathbf{G}_0=\mathbf{M}_0$.
\end{rmk}
\begin{ntz} For a linear subspace $\mathfrak{l}$ of a Lie algebra $\mathfrak{g}$,
we will denote by $\lie(\mathfrak{l})$ the Lie subalgebra of $\mathfrak{g}$
generated by $\mathfrak{l}$: 
\begin{equation}
  \label{eq:xbx}
  \lie(\mathfrak{l})=\mathfrak{l}+[\mathfrak{l},\mathfrak{l}]+
[\mathfrak{l},[\mathfrak{l},\mathfrak{l}]]+\cdots.
\end{equation}
\end{ntz}
Let $M$, $N$ be $\mathbf{G}_0$-homogeneous $CR$ manifolds and 
$\phiup:M\to{N}$ a   $\mathbf{G}_0$-equivariant smooth map.
Fix $p_0\in{M}$ and let  $(\mathfrak{g}_0,\mathfrak{q})$ be
the $CR$ algebra of $M$ at $p_0$, 
$(\mathfrak{g}_0,\mathfrak{e})$ the $CR$ algebra of $N$ at
$\phiup(p_0)$, respectively. Then
\begin{equation}
  \label{eq:a5}
  \mathfrak{q}\cap\bar{\mathfrak{q}}\subset\mathfrak{e}\cap\bar{\mathfrak{e}}
\end{equation} 
and we have
\begin{prop}\label{prop:214}
The $\mathbf{G}_0$-equivariant 
map $M\xrightarrow{\;\phiup\;}N$ is $CR$ if and only if 
$\mathfrak{q}\subset\mathfrak{e}$. \par 
Assume $M\xrightarrow{\;\phiup\;}N$
is $CR$. Then: \begin{align}
\label{eq:a1}
&\text{$\phiup$ is a $CR$-submersion if and only if 
$\mathfrak{e}=\mathfrak{q}+\mathfrak{e}\cap\bar{\mathfrak{e}}$.}\\
 \label{eq:n0}&\text{$\phiup$ is a $CR$-spread if and only if 
$  \mathfrak{e}=\lie(\mathfrak{q}+
\mathfrak{e}\cap\bar{\mathfrak{e}})$.}\\
 \label{eq:n1}&\text{
The fibers of $M\xrightarrow{\;\phiup\;}N$ are totally real if and only if
$  \mathfrak{q}\cap\bar{\mathfrak{e}}=\mathfrak{e}\cap\bar{\mathfrak{q}}
=\mathfrak{q}\cap\bar{\mathfrak{q}}$.}\\
\notag &\text{The fibers of 
$M\xrightarrow{\;\phiup\;}N$ are totally complex if and only if}\qquad\\
&\label{eq:n2}\qquad\qquad\qquad\qquad
\text{$  \mathfrak{q}\cap\bar{\mathfrak{e}}+\mathfrak{e}\cap\bar{\mathfrak{q}}=
 \mathfrak{e}\cap\bar{\mathfrak{e}}.$}
\end{align}
\end{prop}
\begin{proof}
For \eqref{eq:a1}, \eqref{eq:n1}, \eqref{eq:n2} we refer to \cite{MN05}.
\par
Let us prove \eqref{eq:n0}. Denote by $\muup,\nuup$ the $CR$ structures of
$M$, $N$, respectively.
The map $\phiup$ fits into a commutative diagram 
\begin{equation*}
 \xymatrix@R=5pt{&\mathbf{G}_0\ar[dl]_{\piup_M}\ar[dr]^{\piup_N}\\
 M\ar[rr]_{\phiup}&&N.}
\end{equation*}
It is convenient here to identify the elements of $\mathfrak{g}$ 
to the corresponding left-invariant
complex vector fields on $\mathbf{G}_0$.
Then $\piup_M^*{\Zt}_{\muup}$ and 
$\piup_N^*{\Zt}_{\nuup}$ are the distributions of complex vector fields on
$\mathbf{G}_0$ which are generated by the elements
of $\mathfrak{q},\mathfrak{e}$, respectively, and, to verify that
\eqref{eq:n0} is necessary and sufficient
for $\phiup$ being a $CR$ spread, it suffices to note that 
 $\piup^*_M(\lie({\Zt}_\muup+{\Vt}_{\phiup}))=
\mathcal{C}^{\infty}_{\mathbf{G}_0}\otimes
 \lie(\mathfrak{q}+\mathfrak{e}\cap\bar{\mathfrak{e}})$.
\end{proof}
\begin{prop}\label{prop:h215}
Let $M$ be a $\mathbf{G}_0$-homogeneous smooth manifold, 
with isotropy $\mathbf{M}_0$ at $p_0\in{M}$, and  
$\muup_1,\muup_2$ two $\mathbf{G}_0$-equivariant $CR$ structures on
$M$.
Let $(\mathfrak{g}_0,\mathfrak{q}_1)$,
 $(\mathfrak{g}_0,\mathfrak{q}_2)$ be
the $CR$ algebras at $p_0$ corresponding
to $\muup_1,\muup_2$, respectively. Then
\begin{equation}\label{eq:ysc0}
 \mathfrak{q}_1\cap\bar{\mathfrak{q}}_1= \mathfrak{q}_2\cap\bar{\mathfrak{q}}_2=
\mathfrak{m},\;\;\text{with $\mathfrak{m}=\mathbb{C}\otimes\Lie(\mathbf{M}_0)$}, 
\end{equation}
and a necessary and sufficient condition
for $\muup_2$ to strengthen $\muup_1$  is that
\begin{equation}
  \label{eq:xsc0}
 \mathfrak{q}_1\subset\mathfrak{q}_2.
\end{equation}
\end{prop}
The notion of strengthening a $CR$ structure translates into
a corresponding notion for $CR$ algebras:
\begin{defn} We say that $(\mathfrak{g}_0,\mathfrak{q}_2)$
\emph{strengthens}
$(\mathfrak{g}_0,\mathfrak{q}_1)$ if 
\eqref{eq:ysc0} and 
\eqref{eq:xsc0} are satisfied.
\end{defn}
\section{Compact homogeneous $CR$ 
manifolds and $\mathfrak{n}$-reductiveness} \label{sc4}
\subsection{Splittable subalgebras of a complex reductive
Lie algebra}
Let ${\kt}$ be a reductive complex Lie algebra, and 
\begin{equation*}
  {\kt}=\mathfrak{z}\oplus\mathfrak{s},\quad
\mathfrak{z}=\{X\in\kt\mid [X,{\kt}]=\{0\}\},\;\;
\mathfrak{s}=[{\kt},{\kt}]
\end{equation*}
its decomposition into the direct sum of its center and its semisimple
ideal. We say that
$X\in{\kt}$ is \emph{semisimple} if
$\mathrm{ad}(X)$ is a semisimple derivation of ${\kt}$, 
and \emph{nilpotent} if $X\in\mathfrak{s}$ and $\mathrm{ad}(X)$ is
nilpotent.\par 
An equivalent formulation is obtained by
considering a faithful matrix representation of $\kt$ in
which the elements of $\mathfrak{z}$ are diagonal. Semisimple and nilpotent
elements correspond to semisimple and nilpotent matrices, respectively.
\par
Each $X\in{\kt}$ admits a unique Jordan-Chevalley decomposition
\begin{equation}
  \label{eq:10a}
  X=X_s+X_n,\quad \text{with $X_s$ semisimple, $X_n$ nilpotent, and
$[X_s,X_n]=0$.}
\end{equation}\par
A Lie subalgebra $\mathfrak{v}$ of $\kt$ is \emph{splittable} if, 
for each $X\in\mathfrak{v}$,
both $X_s$ and $X_n$ belong to~$\mathfrak{v}$. \par
If $\mathfrak{v}$ is a Lie subalgebra of ${\kt}$, the nilpotent
elements of its radical $\rad(\mathfrak{v})$ form a nilpotent ideal 
$\nr(\mathfrak{v})$
of $\mathfrak{v}$, with
\begin{equation*}
  \radn(\mathfrak{v})=\rad(\mathfrak{v})\cap [\mathfrak{v},\mathfrak{v}]
\subset \nr(\mathfrak{v})\subset\nil(\mathfrak{v}),
\end{equation*}
where $\nil(\mathfrak{v})$ is the nilradical, i.e. the maximal nilpotent ideal
of $\mathfrak{v}$, and $\radn(\mathfrak{v})$ its nilpotent radical, i.e. the 
intersection of the kernels of all irreducible finite dimensional
linear representations of
$\mathfrak{v}$. Note that the nilpotent ideal $\nr(\mathfrak{v})$, unlike
$\nil(\mathfrak{v})$ and $\radn(\mathfrak{v})$, depends on the inclusion
$\mathfrak{v}\subset{\kt}$ (cf. \cite[\S 5.3]{Bou75}).\par
We recall (see \cite[\S 5.4]{Bou75}) 
\begin{prop}
 Every splittable Lie subalgebra $\mathfrak{v}$ 
admits a Levi-Chevalley decomposition 
\begin{equation}\label{eq:42}
 \mathfrak{v}=\nr(\mathfrak{v})\oplus\mathfrak{m},
\end{equation}
where $\mathfrak{m}$ is reductive in $\kt$, and is uniquely determined modulo
conjugation by elementary automorphisms of $\mathfrak{v}$, i.e. finite products of
automorphisms of the form $\exp(\ad(X))$, with $X\in\mathfrak{v}$ and nilpotent.
\end{prop}
\subsection{Definition of $\mathfrak{n}$-reductive}
Let ${\kt}$ be the complexification of 
a compact Lie algebra ${\kt}_0$. 
Conjugation in ${\kt}$ is understood
with respect to the real form ${\kt}_0$. Note that all subalgebras of
a compact Lie algebra are compact and hence reductive.
\begin{prop} \label{prop:31}
For any complex Lie subalgebra
$\mathfrak{v}$ of ${\kt}$, the intersection  
$\mathfrak{v}\cap\bar{\mathfrak{v}}$ is reductive and splittable. 
In particular,
$  \mathfrak{v}\cap\bar{\mathfrak{v}}\cap\nr(\mathfrak{v})=\{0\}$. A
splittable $\mathfrak{v}$ 
admits a Levi-Chevalley decomposition with a reductive Levi factor containing
  $\mathfrak{v}\cap\bar{\mathfrak{v}}$.
\end{prop}
\begin{proof} We recall that $\mathfrak{v}$ is splittable if and only if 
$\rad({\mathfrak{v}})$ is splittable. When this is the case, 
$\mathfrak{v}$ admits a Levi-Chevalley decomposition, and all maximal reductive
Lie subalgebras of $\mathfrak{v}$ can be taken as 
reductive Levi factors.
The intersection $\mathfrak{v}\cap\bar{\mathfrak{v}}$ is reductive, because it
is the complexification of the compact Lie algebra
$\mathfrak{v}\cap{{\kt}}_0$. Then, the reductive Levi factor in the
Levi-Chevalley decomposition of $\mathfrak{v}$ can be taken to contain
$\mathfrak{v}\cap\bar{\mathfrak{v}}$ (see e.g. \cite{OV90}).
\end{proof}

Let $\mathbf{K}_0$ be a compact Lie group with Lie algebra
$\kt_0$ and $M$ a $\mathbf{K}_0$-homogeneous $CR$ manifold,
with isotropy $\mathbf{M}_0$ and $CR$ algebra
$({\kt}_0,\mathfrak{v})$ at a point
$p_0\in{M}$.
\begin{defn}\label{df33}
We say that $M$, and the $CR$ algebra $({\kt}_0,\mathfrak{v})$, 
are $\mathfrak{n}$-reductive if 
$\mathfrak{v}\cap\bar{\mathfrak{v}}$
is a reductive complement of
$\nr(\mathfrak{v})$ in $\mathfrak{v}$.  
\end{defn}
\begin{rmk}
If $(\kt_0,\mathfrak{v})$ is $\mathfrak{n}$-reductive, then
$\mathfrak{v}$ is splittable. Indeed all elements of 
$\nr(\mathfrak{v})$ are nilpotent and all elements of 
$\mathfrak{v}\cap\bar{\mathfrak{v}}$ are splittable, because
this subalgebra is the complexification of the subalgebra 
$\mathfrak{v}\cap\kt_0$, 
which is
splittable because consists of semisimple elements.  
\end{rmk}
\begin{thm}[realization]\label{thm45}
Every $\mathfrak{n}$-reductive $\mathbf{K}_0$-homogeneous $CR$ manifold 
is $\mathbf{K}$-realizable
(see Definition\,\ref{dfn:26}). \par
In fact, if $M$ is a $\mathbf{K}_0$-homogeneous $CR$ manifold
with isotropy $\mathbf{M}_0$ and $CR$ algebra $(\kt_0,\mathfrak{v})$
at $p_0\in{M}$, then
there is a closed complex
subgroup $\mathbf{V}$ of $\mathbf{K}$ with $\Lie(\mathbf{V})=\mathfrak{v}$
\;
and \; $\mathbf{V}\cap\mathbf{K}_0=\mathbf{M}_0$. 
Set $\mathfrak{m}_0=\Lie(\mathbf{M}_0)$, 
$\mathfrak{m}=\mathbb{C}\otimes_{\mathbb{R}}\mathfrak{m}_0$. \par
The Lie subgroup $\mathbf{V}$ admits 
a Levi-Chevalley decomposition
\begin{equation}
  \label{eq:6x16}
  \mathbf{V}=\mathbf{V}_n\mathbf{M},\quad\text{with \;
$\Lie(\mathbf{V}_n)=\nr(\mathfrak{v})$, \; $\Lie(\mathbf{M})=
\mathfrak{m}$}.
\end{equation}
\end{thm}
\begin{proof} The linear subspace 
$i\mathfrak{m}_0$ of $\kt$ 
is $\Ad_{\kt}(\mathbf{M}_0)$-invariant, and 
the unique 
complexification $\mathbf{M}$ of $\mathbf{M}_0$ in $\mathbf{K}$
is defined by
\begin{equation*}
 \mathbf{M}=\{g\exp(iX)\mid g\in\mathbf{M}_0,\; X\in\mathfrak{m}_0\}.
\end{equation*}
Clearly $\Lie(\mathbf{M})=\mathfrak{m}$. 
The complex analytic Lie subgroup $\mathbf{V}_n$ 
of $\mathbf{K}$ with Lie algebra
$\nr(\mathfrak{v})$ is closed, is normalized by $\mathbf{M}$
and $\mathbf{V}_n\cap\mathbf{M}$ is finite. Hence
\begin{equation*}
  \mathbf{V}=\mathbf{V}_n\mathbf{M}=
\{g_1g_2\mid g_1\in\mathbf{V}_n,\; g_2\in\mathbf{M}\}
\end{equation*}
is a closed Lie subgroup of $\mathbf{K}$, with $\mathbf{V}\cap\mathbf{K}_0=
\mathbf{M}_0$ and $\Lie(\mathbf{V})=\mathfrak{v}$. With 
$V=\mathbf{K}/\mathbf{V}$, the inclusion $\mathbf{K}_0\hookrightarrow
\mathbf{K}$ passes to the quotients, defining a $\mathbf{K}_0$-equivariant
$CR$-embedding $M\hookrightarrow{V}$.
\end{proof}
Later, in \S\ref{sec:mf}, we will encounter a large class of 
$\mathfrak{n}$-reductive compact homogeneous $CR$
manifolds. 
We exhibit here
an example of a compact homogeneous $CR$ manifold  $M$ which is not 
$\mathfrak{n}$-reductive.
\begin{exam}
Let $\mathbf{K}_0=\mathbf{SU}(n)$, $n\geq 3$. Fix a complex symmetric 
non degenerate $n\times{n}$ matrix $S$ 
and consider the subgroup $\mathbf{V}=\{a\in\mathbf{SL}(n,\mathbb{C})\mid
a^tSa=S\}$ of $\mathbf{SL}(n,\mathbb{C})$, with Lie algebra
$\mathfrak{v}=\{X\in\mathfrak{sl}(n,\mathbb{C})\mid
X^tS+SX=0\}$. 
Set $\mathbf{M}_0=\mathbf{V}\cap\mathbf{K}_0$ and
$M=\mathbf{K}_0/\mathbf{M}_0$. This is a $\mathbf{K}_0$-homogeneous 
$CR$ manifold with $CR$ algebra 
$(\kt_0,\mathfrak{v})$, where
$\kt_0\simeq\mathfrak{su}(n)$, $\mathfrak{v}\simeq\mathfrak{so}(n,\mathbb{C})$.
If $S,S^*$ are linearly independent, then $\mathfrak{v}$
is a semisimple Lie subalgebra of $\kt$ distinct from
$\mathfrak{v}\cap\bar{\mathfrak{v}}$.
\end{exam}
\subsection{Regular and parabolic subalgebras}
Let $\mathfrak{t}_0$ be a maximal torus of $\kt_0$.
\begin{defn}  We call
$\mathfrak{t}_0$-\emph{regular} a subalgebra 
   $\mathfrak{v}$ of ${\kt}$ which is 
normalized by  $\mathfrak{t}_0$.
\end{defn}
The complexification $\mathfrak{t}$ of 
$\mathfrak{t}_0$ is a Cartan subalgebra of $\kt$.
Denote by $\mathcal{K}$ the root system of $(\kt,\mathfrak{t})$
and by 
\begin{displaymath}
\kt^{\alpha}=\{X\in\kt\mid [H,X]=\alpha(H)X,\;\forall H\in\mathfrak{t}\}
\end{displaymath}
the root space of $\alpha\in\mathcal{K}$. Note that
$\overline{\kt^{\alpha}}=\kt^{\,-\alpha}$.
\begin{lem}\label{lem:a2}  
 A $\mathfrak{t}_0$-regular subalgebra $\mathfrak{v}$ of ${\kt}$ is splittable and 
 admits a Levi-Chevalley decomposition 
\begin{equation}
  \label{eq:10}
  \mathfrak{v}=\nr(\mathfrak{v})\oplus\mathfrak{m},\;\;
\text{with reductive Levi factor
  \;\;$\mathfrak{m}\, =\, \mathfrak{t}\!\cap\!\mathfrak{v}
\,+\,\mathfrak{v}\!\cap\!\bar{\mathfrak{v}}$}.
\end{equation}
In particular, $(\kt_0,\mathfrak{v})$ is $\mathfrak{n}$-reductive if 
and only if \; $\mathfrak{t}\cap\mathfrak{v}\subset\bar{\mathfrak{v}}$.
\end{lem}
\begin{proof} A $\mathfrak{t}_0$-regular $\mathfrak{v}$ 
admits a root-space decomposition
\begin{equation*}
  \mathfrak{v}=\mathfrak{a}\oplus{\sum}_{\alpha\in\mathcal{V}}{\kt}^{\alpha},
\;\;\text{with $\mathfrak{a}=\Z_{\mathfrak{v}}(\mathfrak{t})
=\mathfrak{t}\cap\mathfrak{v}$, 
$\mathcal{V}=\{\alpha\in\mathcal{K}\mid {\kt}^{\alpha}\subset
\mathfrak{v}\}$}.
\end{equation*}
Thus the Lie algebra $\mathfrak{v}$ admits a set of generators 
that are either semisimple or
nilpotent and is splittable by  \cite[Ch.VII,\S{5}]{Bou75}.
\par
The conjugation with respect to the real form ${\kt}_0$ yields by duality 
the conjugation $\bar{\alpha}=-\alpha$ in $\mathcal{K}$. 
Set
\begin{align*}
\nr(\mathfrak{v})&=  {\sum}_{\alpha\in\mathcal{V}_n}{\kt}^{\alpha},
&&\quad \text{with $\mathcal{V}_n=\{\alpha\in\mathcal{V}\mid -
\alpha\notin\mathcal{V}\}$},\\ \mathfrak{m}=
\mathfrak{a}+\mathfrak{v}\cap\bar{\mathfrak{v}}&=
\mathfrak{a}\oplus{\sum}_{\alpha\in\mathcal{V}_r}{\kt}^{\alpha},
&&\quad
\text{with $\mathcal{V}_r=\{\alpha\in\mathcal{V}\mid -\alpha\in\mathcal{V}\}=
\mathcal{V}\cap\bar{\mathcal{V}}$},
\end{align*}
and \eqref{eq:10} follows, as one easily checks that $\nr(\mathfrak{v})$
is the ideal of nilpotent elements of $\rad(\mathfrak{v})=
\mathfrak{a}\oplus\nr(\mathfrak{v})$ and 
$\mathfrak{m}$ is reductive.
\end{proof}
If $A\in\kt_0$, then $\ad_{\kt}(A)$ is semisimple with 
purely imaginary eigenvalues. The sum of the corresponding eigenspaces 
with positive imaginary part
\begin{equation}
  \label{eq:p0}
  \mathfrak{q}_A={\sum}_{\lambda\geq{0}}\{X\in\kt\mid [A,X]=i\lambda{X}\}
\end{equation}
is a \emph{parabolic} Lie subalgebra of $\kt$, and all parabolic subalgebras
of $\kt$ have this form. 
The intersection of 
$\mathfrak{q}_A$ with $\kt_0$ is the centralizer of $A$ in $\kt_0$:
\begin{equation}
\Li_0(\mathfrak{q}_A)=\mathfrak{q}_A\cap\kt_0=\Z_{\kt_0}(A)=\{X\in\kt_0\mid
[A,X]=0\}
\end{equation}
and is the union of the maximal tori of $\kt_0$ contained in $\mathfrak{q}_A$.
Its complexification is the centralizer of $A$ in $\kt$ and equals the intersection
with its $\kt_0$-conjugate:
\begin{equation}
 \Li(\mathfrak{q}_A)=\Z_{\kt}(A)=
 \mathfrak{q}_A\cap\bar{\mathfrak{q}}_A.
\end{equation}
This subalgebra $\Li(\mathfrak{q}_A)$ is the unique conjugation invariant
reductive Levi factor in the Levi-Chevalley decomposition 
\begin{equation}
 \mathfrak{q}_A=\nr(\mathfrak{q}_A)\oplus\Li(\mathfrak{q}_A)=\nr(\mathfrak{q}_A)
 \,\oplus\,\mathfrak{q}_A\!\cap\!\bar{\mathfrak{q}}_A.
\end{equation}
\begin{ntz} In the following, for $\mathfrak{v}$ complex subalgebra of $\kt$,
we shall use the notation
\begin{equation*}
  \Li_0(\mathfrak{v})=\mathfrak{v}\cap\kt_0,\quad \Li(\mathfrak{v})=
\mathfrak{v}\cap\bar{\mathfrak{v}}.
\end{equation*}
\par
 We denote by $\Par$ the set of complex parabolic Lie subalgebras of $\kt$.
\end{ntz}
\begin{rmk}
 If $\mathfrak{q}\in\Par$, then the set 
$\{A\in\kt_0\mid \mathfrak{q}_A=\mathfrak{q}\}$
 is an open cone in the real vector space
 $\Z_{\kt_0}(\Li_0(\mathfrak{q}))=\{A\in\kt_0
 \mid [A,\mathfrak{q}\cap\bar{\mathfrak{q}}]=\{0\}\}$.
\end{rmk}
\begin{lem}\label{lem:11} 
If $\mathfrak{q}\in\Par$, then $\nr(\mathfrak{q})=\radn(\mathfrak{q})$
and $(\kt_0,\mathfrak{q})$ is $\mathfrak{n}$-reductive
and totally complex.
\end{lem}
\begin{proof} An element $\mathfrak{q}\in\Par$ contains a 
a maximal torus $\mathfrak{t}_0$ of $\kt_0$ (see e.g. \cite[Ch.VIII]{Bou75}). 
Let
$\mathfrak{t}=\mathbb{C}\!\otimes_{\mathbb{R}}\!\mathfrak{t}_0$ be its complexification,
$\mathcal{K}$ the root system of $(\kt,\mathfrak{t})$ and 
$\mathcal{Q}=\{\alpha\in\mathcal{K}\mid \kt^{\alpha}\subset\mathfrak{q}\}$
the parabolic set of $\mathfrak{q}$. 
With $\mathcal{Q}_r=\{\alpha\in\mathcal{Q}\mid
-\alpha\in\mathcal{Q}\}=\mathcal{Q}\cap\bar{\mathcal{Q}}$, 
$\mathcal{Q}_n=\{\alpha\mid
-\alpha\notin\mathcal{Q}\}=\mathcal{Q}\setminus\bar{\mathcal{Q}}$,
we have
\begin{equation*}
\mathfrak{q}=\nr(\mathfrak{q})\oplus\Li(\mathfrak{q}),\;\;\text{with}\;\;
\nr(\mathfrak{q})={\sum}_{\alpha\in\mathcal{Q}_n}{\kt}^{\alpha},
\;\;
\Li(\mathfrak{q})=\mathfrak{t}\oplus{\sum}_{\alpha\in\mathcal{Q}_r}{\kt}^{\alpha}.
\end{equation*}
Since
\begin{equation*}
  \overline{\nr(\mathfrak{q})}={\sum}_{\alpha\in\mathcal{Q}_n}\overline{{\kt}^{\alpha}}
  ={\sum}_{\alpha\in\mathcal{Q}_n}{\kt}^{-\alpha},
\end{equation*}
and 
$\mathcal{Q}\cup(-\mathcal{Q})=\mathcal{K}$ because $\mathcal{Q}$ is parabolic,
we obtain
\begin{equation*}
  \kt=\nr(\mathfrak{q})\oplus\Li(\mathfrak{q})\oplus\overline{\nr(\mathfrak{q})}
=\mathfrak{q}+\bar{\mathfrak{q}},
\end{equation*}
showing that $({\kt}_0,\mathfrak{q})$ is totally complex
(see e.g. \cite[p.237]{MN05}).
\end{proof}
\subsection{Morphisms of $\mathfrak{n}$-reductive $CR$ algebras}
Let $M$ be a $\mathbf{K}_0$-homogeneous $CR$ manifold with isotropy subgroup
$\mathbf{M}_0$ and $CR$ algebra $(\kt_0,\mathfrak{v})$ at  $p_0\in{M}$.
If $M$ is $\mathfrak{n}$-reductive, we can identify $T^{0,1}_{p_0}M$ with
$\nr(\mathfrak{v})$. Therefore, a $\kt_0$-equivariant morphism of
$\mathfrak{n}$-reductive
$CR$ algebras $(\kt_0,\mathfrak{v})\to(\kt_0,\mathfrak{q})$ induces a 
$\mathbb{C}$-linear map $\nr(\mathfrak{v})\to\nr(\mathfrak{q})$, which is
defined by the commutative diagram 
\begin{equation*}
 \begin{CD}
 \mathfrak{v} @>>> \mathfrak{q}\\
 @VVV @VVV\\
 \nr(\mathfrak{v})@>>>\nr(\mathfrak{q})
\end{CD}
\end{equation*}
in which the top horizontal arrow is inclusion and the two vertical arrows are projections
along $\mathfrak{v}\cap\bar{\mathfrak{v}}$ and $\mathfrak{q}\cap\bar{\mathfrak{q}}$,
respectively.
\begin{prop}
 Let $\mathfrak{v},\mathfrak{q}$ be complex Lie subalgebras of $\kt$, with
 $\mathfrak{v}\subset\mathfrak{q}$, and assume that
 $(\kt_0,\mathfrak{v})$ is
$\mathfrak{n}$-reductive. Set $\Li(\mathfrak{v})=\mathfrak{v}\cap\bar{\mathfrak{v}}$,
$\Li(\mathfrak{q})=\mathfrak{q}\cap\bar{\mathfrak{q}}$.
 The $\kt_0$-equivariant
 fibration 
\begin{equation}\label{eq:bb}
 (\kt_0,\mathfrak{v})\longrightarrow(\kt_0,\mathfrak{q})
\end{equation}
\begin{list}{-}{}
\item is a $CR$-submersion if and only if $\mathfrak{q}=\nr(\mathfrak{v})+
\Li(\mathfrak{q})$;
 \item is a $CR$-submersion with totally real fibers 
 if and only if $\mathfrak{q}=\nr(\mathfrak{v})\oplus
\Li(\mathfrak{q})$.
\end{list}
In particular, \begin{list}{-}{}
\item
  is a $CR$-submersion if 
 $\nr(\mathfrak{q})\subset\nr(\mathfrak{v})$;
 \item
 is a $CR$-submersion with totally real fibers 
 if 
 $\nr(\mathfrak{q})=\nr(\mathfrak{v})$.
\end{list}
\end{prop} 
\begin{proof}
If $(\kt_0,\mathfrak{v})$ is $\mathfrak{n}$-reductive, we have 
\begin{equation*}
 \mathfrak{v}+\Li(\mathfrak{q})=\nr(\mathfrak{v})+\Li(\mathfrak{q})
\end{equation*}
because $\Li(\mathfrak{v})\subset\Li(\mathfrak{q})$. 
In view of \eqref{eq:a1},
this equality yields the above
characterization of $CR$-submersion from a $\mathfrak{n}$-reductive
$CR$ algebra.\par
If $\mathfrak{q}=\nr(\mathfrak{v})\oplus
\Li(\mathfrak{q})$, then $\bar{\mathfrak{q}}=\overline{\nr(\mathfrak{v})}\oplus
\Li(\mathfrak{q})$ and we have  a direct sum decomposition 
\begin{equation*}
 \mathfrak{q}+\bar{\mathfrak{q}}=\nr(\mathfrak{v})\oplus
\Li(\mathfrak{q})\oplus\overline{\nr(\mathfrak{v})}.
\end{equation*}
An element $X\in\mathfrak{v}\cap\bar{\mathfrak{q}}$ uniquely decomposes as
sums $X=X_1+X_2=Y_1+Y_2$ with $X_1\in\nr(\mathfrak{v})$, $X_2\in\Li(\mathfrak{v})
\subset\Li(\mathfrak{q})$, $Y_1\in\Li(\mathfrak{q})$, 
$Y_2\in\overline{\nr(\mathfrak{v})}$.
By comparing the two decompositions, we obtain that $X_1=Y_2=0$ and 
$X=Y_1=X_2\in\Li(\mathfrak{v})$.\par
Vice versa, if $\mathfrak{v}\cap\bar{\mathfrak{q}}=\Li(\mathfrak{v})$, then 
$\nr(\mathfrak{v})\cap\Li(\mathfrak{q})=\{0\}$ and the sum $\nr(\mathfrak{v})+
\Li(\mathfrak{q})$ is direct.
\end{proof}
\begin{exam}
Let us consider the $CR$ algebras $(\mathfrak{u}_3,\mathfrak{v})$ and
 $(\mathfrak{u}_3,\mathfrak{q})$ with
\begin{align*}
 \mathfrak{v}&=\left.\left\{\left( 
\begin{smallmatrix}
 \lambda&\zeta_1&\zeta_2\\
 0&\lambda&\zeta_1\\
 0&0&\lambda
\end{smallmatrix}\right)\right| \lambda,\zeta_1,\zeta_2\in\mathbb{C}\right\},\\
\mathfrak{q}&=\left.\left\{\left( 
\begin{smallmatrix}
 \lambda_1&\zeta_1&\zeta_2\\
 0&\lambda_2&\zeta_3\\
 0&\zeta_4&\lambda_3
\end{smallmatrix}\right)\right| \lambda_1,\lambda_2,\lambda_3
,\zeta_1,\zeta_2,\zeta_4,\zeta_2\in\mathbb{C}\right\}.
\end{align*}
We have 
\begin{equation*}
 \nr(\mathfrak{v})=\left\{\left( 
\begin{smallmatrix}
 0&\zeta_1&\zeta_2\\
 0&0&\zeta_1\\
 0&0&0
\end{smallmatrix}\right)\right\},\quad \nr(\mathfrak{q})=\left\{\left( 
\begin{smallmatrix}
 0&\zeta_1&\zeta_2\\
 0&0&0\\
 0&0&0
\end{smallmatrix}\right)\right\},
\end{equation*}
so that $\mathfrak{q}=\nr(\mathfrak{v})\oplus\Li(\mathfrak{q})$, showing that
$(\mathfrak{u}_3,\mathfrak{v})\to(\mathfrak{u}_3,\mathfrak{q})$ is a 
$CR$-submersion
with totally real fibers, for which
$\nr(\mathfrak{q})\not\subset\nr(\mathfrak{v})$.
\end{exam}
\subsection{Regularity type for 
$\mathfrak{n}$-reductive $CR$ algebras}
\begin{defn}
 We say that an $\mathfrak{n}$-reductive $CR$ algebra $(\kt_0,\mathfrak{v})$ is
of
\begin{list}{-}{}
\item type $\mathrm{I}$ if $\mathfrak{v}$ is $\mathfrak{t}_0$-regular, for
some maximal torus $\mathfrak{t}_0$ of $\kt_0$,\\
\item type $\mathrm{I\! I}$ if $\Li(\mathfrak{v})$ 
is $\mathfrak{t}_0$-regular, for
some maximal torus $\mathfrak{t}_0$ of $\kt_0$,\\
\item type $\mathrm{I\! I\! I}$, or general, if we make no assumption on
normalizing tori. 
\end{list}
\end{defn}
\subsection{Minimal orbits of compact forms}\label{sec46}
Let $\mathbf{K}$ be a connected
reductive complex Lie group, with $\Lie(\mathbf{K})=\kt$. If
$\mathbf{V}$ is an algebraic subgroup of $\mathbf{K}$, then $\mathbf{V}$
admits a Levi-Chevalley decomposition
\begin{equation*}
 \mathbf{V}= \mathbf{V}_n\mathbf{V}_r,\;\;\text{with
$\mathbf{V}_n$ unipotent, $\mathbf{V}_r$ reductive, $|
\mathbf{V}_n\cap\mathbf{V}_r|<\infty$. }
\end{equation*}
The maximal compact subgroup $\mathbf{M}_0$ of 
$\mathbf{V}_r$ 
is contained in a compact form $\mathbf{K}_0$ of
$\mathbf{K}$. The orbit $M$ of $\mathbf{K}_0$ through the base point
of the complex 
homogeneous space 
$V=\mathbf{K}/\mathbf{V}$
inherits from $V$ 
the structure of a $\mathbf{K}_0$-homogeneous
$\mathfrak{n}$-reductive $CR$ manifold,
with isotropy $\mathbf{M}_0$ and $CR$ algebra $(\kt_0,\mathfrak{v})$
at the base point, where $\kt_0=\Lie(\mathbf{K}_0)$ and 
$\mathfrak{v}=\Lie(\mathbf{V})$. Since all compact forms of $\mathbf{K}$
are conjugated by an inner automorphism, we obtain the
\begin{prop}
Let $V$ be the homogeneous complex manifold obtained as the quotient 
of a connected reductive complex manifold $\mathbf{K}$ by an algebraic
subgroup $\mathbf{V}$, and 
$\mathbf{K}_0$ a compact form of $\mathbf{K}$.
Then $\mathbf{K}_0\backslash V$ contains $\mathfrak{n}$-reductive orbits.
All $\mathfrak{n}$-reductive orbits in $\mathbf{K}_0\backslash V$ are
$CR$-diffeomorphic. \par
The $\mathfrak{n}$-reductive orbits are the orbits of minimal dimension
in $\mathbf{K}_0\backslash V$.
\end{prop}
\begin{proof}
We already proved existence.
Let us
prove equivalence. We can assume that the orbit $M$ through the base
point $p_0$ is $\mathfrak{n}$-reductive. Let $p_1=g_1V$, $g_1\in\mathbf{K}$,
and assume that also $M'=K_0p_1$ is $\mathfrak{n}$-reductive. This 
means that $\mathbf{K}_0\cap{\ad(g_1)(\mathbf{V})}$ is 
a maximal compact subgroup of $\ad(g_1)(\mathbf{V})$,
i.e. that $\ad(g_1^{-1})(\mathbf{K}_0)\cap{\mathbf{V}}$ is a  maximal compact subgroup of
$\mathbf{V}$. By the conjugacy of  maximal compact subgroups,
there is an element $g_2\in\mathbf{V}$ such that
$\ad(g_1^{-1})(\mathbf{K}_0)\cap\mathbf{V}=\ad(g_2)(\mathbf{K}_0)\cap\mathbf{V}$. 
Then we obtain,
for $k\in\mathbf{K}_0$, 
\begin{equation*}
  g_1\!\cdot\!{g}_2\cdot k\!\cdot\! p_0=
\ad(g_1\!\cdot\!{g}_2)(k)\cdot g_1\!\cdot\!{g}_2\cdot p_0 
= k'\cdot p_1,\;\;\text{with $k'=\ad(g_1\!\cdot\!{g}_2)(k)\in\mathbf{K}_0$}.
\end{equation*}
This shows that the translation in $V$ by $g_1\!\cdot\!{g}_2$ transforms
$M$ onto $M'$. This map is a biholomorphism of $V$ and then its restriction
to $M$ is a $CR$-diffeomorphism  of $M$ onto $M'$. \par
The fact that the $\mathfrak{n}$-reductive orbits have minimal dimension
follows because the isotropy $\mathbf{M}_p$ of a $\mathbf{K}_0$-orbit
$M$ at a point $p=g\cdot{p}_0$ is a reductive subgroup of 
$\ad(g)(\mathbf{V})$, and therefore $\dim_{\mathbb{R}}\mathbf{M}_0\leq
\dim_{\mathbb{C}}\mathbf{V}_r$, and we have equality if and only if 
$\mathbf{M}_0$ is a maximal compact subgroup of $\ad(g)(\mathbf{V})$.
\end{proof}
\subsection{$\mathfrak{n}$-reductive structures on compact
homogeneous spaces}\label{sub47}
Let $\mathbf{K}_0$ be a connected compact Lie group and $M$ a
$\mathbf{K}_0$-homogeneous space, with isotropy 
$\mathbf{M}_0$ at $p_0\in{M}$. With $\Lie(\mathbf{M}_0)=\mathfrak{m}_0$,
$\mathfrak{m}=\mathbb{C}\otimes_{\mathbb{R}}\mathfrak{m}_0$, the 
$CR$ algebra $(\kt_0,\mathfrak{m})$ defines on $M$ the trivial totally real
$CR$ structure. In accordance with
Theorem\,\ref{thm45}, there is a reductive complex Lie subgroup 
$\mathbf{M}$ of the complexification $\mathbf{K}$ of $\mathbf{K}_0$
such that, with $M^{\mathbb{C}}=\mathbf{K}/\mathbf{M}$, we have a
$\mathbf{K}_0$-equivariant immersion $M\hookrightarrow{M}^{\mathbb{C}}$
of $M$ into its complexification $M^{\mathbb{C}}$, 
which is Stein because $\mathbf{K}$ is Stein and 
$\mathbf{M}$  is the complexification of the compact subgroup
$\mathbf{M}_0$ (see \cite{Mts60}).
\par
We have
\begin{prop}
The $\mathfrak{n}$-reductive $\mathbf{K}_0$-equivariant $CR$ structures
on $M$ are in a one-to-one correspondence with the 
$\Ad_{\kt}(\mathbf{M}_0)$-invariant complex subalgebras
$\mathfrak{n}$ of $\kt$ consisting of nilpotent elements. \par
If $(\kt_0,\mathfrak{v})$ is the $CR$ algebra at $p_0$ of an 
$\mathfrak{n}$-reductive $\mathbf{K}_0$-equivariant $CR$ structure
on $M$, then, with the notation
of Theorem\,\ref{thm45}, we have a commutative diagram
\begin{equation*}
  \begin{CD}
    M_{(\kt_0,\mathfrak{m})} @>>> M^{\mathbb{C}}\\
@VVV  @VV{\piup_V}V \\
M_{(\kt_0,\mathfrak{v})} @>>> V
  \end{CD}
\end{equation*}
where $ M_{(\kt_0,\mathfrak{m})}$ and $M_{(\kt_0,\mathfrak{v})}$ are the
$CR$ manifolds with underlying smooth manifold $M$ and 
$\mathbf{K}_0$-equivariant $CR$ structures 
defined by the $CR$ algebras $(\kt_0,\mathfrak{m})$,
$(\kt_0,\mathfrak{v})$, respectively; the horizontal arrows are
$\mathbf{K}$-realizations, the left vertical arrow is a strengthening
of the $CR$ structure, 
and the right vertical arrow is a holomorphic 
submersion with Euclidean fibers.
\end{prop}
\begin{proof}
Indeed, the fibers of the $\mathbf{K}$-equivariant holomorphic submersion
$M^{\mathbb{C}}\to{V}$ are complex nilmanifolds conjugate to $\mathbf{V}_n$.  
\end{proof}
\begin{defn} \label{defn417}
We call the fibration $M^{\mathbb{C}}\xrightarrow{\;\piup_V\;}V$ 
a \emph{Stein lift} of $V$.  
\end{defn}
\subsection{Complexification. Foliated complex manifolds}
\label{sub48}
We keep the notation of~\S{\ref{sub47}}. 

\begin{rmk}
The Stein lift of Definition\,\ref{defn417}
is characterized by the universal property that
every $\mathbf{K}$-equivariant fibration $E\xrightarrow{\;\varpi\;}V$ of a
$\mathbf{K}$-homogeneous Stein manifold $E$ onto $V$ factors through  
$M^{\mathbb{C}}\xrightarrow{\;\piup_V\;}V$.      
\end{rmk}
The fibers of $M^{\mathbb{C}}\xrightarrow{\;\piup_V\;}V$ are the integral
manifolds of a holomorphic distribution ${\Dt}_{\mathfrak{v}}$ on 
$M^{\mathbb{C}}$.
\begin{defn} We call the \textit{foliated complex manifold}
$(M^{\mathbb{C}},{\Dt}_{\mathfrak{v}})$ a \emph{complexification}
of the $CR$ manifold $M$.
\end{defn}
This construction can be generalized to the case of a real analytic
$CR$ manifold $M$. Indeed, $M$ has a complexification 
$M\xrightarrow{\;\iotaup_M\;}M^{\mathbb{C}}$
and a generic embedding $M\xrightarrow{\;\iotaup_V\;}{V}$ (see \cite{AF79}). 
Since $\imath_V$ is real analytic, it extends uniquely to a holomorphic
submersion $M^{\mathbb{C}}\xrightarrow{\;\piup_V}V$ (after shrinking $V$ and
$M^{\mathbb{C}}$, if needed). Since $M$ has in $M^{\mathbb{C}}$ a fundamental 
system of Stein neighborhoods, we can say that $M$ admits 
a complexification 
$(M^{\mathbb{C}},{\Dt})$ with a Stein $M^{\mathbb{C}}$ and ${\Dt}$
yielding a fibration on a generic complex neighborhood $V$ of $M$.\par 
In this more general context, the leaves of ${\Dt}$ are also called 
\textsl{Segre varieties} of~$M$.\par
For further reference, we collect below some general notions on 
foliated complex manifolds.
\begin{defn}
A \emph{foliated complex manifold} is a pair $(E,{\Dt})$ consisting
of a complex manifold $E$ and an integrable holomorphic distribution
(of constant rank) ${\Dt}$ on $E$. We set 
$D_pE=\{X_p\mid X\in{\Dt}\}\subset{T}^{\mathbb{C}}E$, and
$DE={\bigsqcup}_{p\in{E}}D_pE$. 
\par
If $(E',{\Dt}')$ is another foliated complex manifold, a holomorphic
map $\phiup:E\to{E}'$ is a ${\Dt}$-\emph{morphism} 
if $\phiup_*(DE)\subset D'E'$. We say that $\phiup$ is
\begin{list}{-}{}
\item a ${\Dt}$-submersion if $\phiup_*(D_pE)=D'_{\phiup(p)}E'$
for all $p\in{E}$;
\item a ${\Dt}$-spread if it is a ${\Dt}$-submersion and
$\phiup^*{\Dt}'$ is generated by ${\Dt}+\mathcal{V}_{\phiup}$,
where $\mathcal{V}_{\phiup}$ is the vertical holomorphic distribution of
$\phiup$;
\item a ${\Dt}$-deployment if it is a ${\Dt}$-spread and
the maps $\phiup_*(p):D_pE\to D'_{\phiup(p)}E'$ are injective for all
$p\in{E}$. 
\end{list}
\end{defn}
\section{A class of $\mathfrak{n}$-reductive
compact homogeneous $CR$ manifolds} 
\label{sec:mf}
In this section we describe a large class of $\mathfrak{n}$-reductive
compact homogeneous $CR$ manifold, 
which are $CR$ submanifolds of 
complex flag manifolds. 
Although it could be proved that these manifolds are $\mathfrak{n}$-reductive 
by using the results in \S\ref{sec46}, we give a direct 
and more explicit proof of this fact.
\subsection{Complex flag manifolds} \label{su41}
 Let $\mathbf{G}$ 
be a connected complex algebraic 
Lie group and $\mathbf{F}$ a parabolic subgroup of $\mathbf{G}$. 
If $\mathfrak{g}=\Lie(\mathbf{G})$, ${\ft}=\Lie(\mathbf{F})$, 
we have \begin{equation}
  \label{eq:4a}
  \mathbf{F}=\N_{\mathbf{G}}({\ft})=
\{g\in\mathbf{G}\mid\Ad_{\mathfrak{g}}(g)({\ft})={\ft}\}.
\end{equation}
The group $\mathbf{F}$ is connected and $F=\mathbf{G}/\mathbf{F}$ is
a complex projective variety. \par
Let $\mathfrak{h}$ be a Cartan subalgebra of $\mathfrak{g}$,
contained in ${\ft}$. Its elements are semisimple.
Let $\mathfrak{h}_{\mathbb{R}}$ be the set of $H\in\mathfrak{h}$
for which $\ad({H})$ has real eigenvalues. Then $i\mathfrak{h}_{\mathbb{R}}$
is a maximal torus in $\mathfrak{g}$, and there is a compact form
$\mathfrak{u}_0$ of $\mathfrak{g}$ containing 
$i\mathfrak{h}_{\mathbb{R}}$. The analytic subgroup $\mathbf{U}_0$ with
$\Lie(\mathbf{U}_0)=\mathfrak{u}_0$ is a maximal compact subgroup of
$\mathbf{G}$ and acts transitively on $F$. Let $\tauup$ be the conjugation
on $\mathfrak{g}$ defined by the real form $\mathfrak{u}_0$. There is
a unique Levi-Chevalley decomposition 
\begin{equation}
  \label{eq:43}
  {\ft}=\rn({\ft})\oplus\Li_{\tauup}({\ft}),\quad
\Li_{\tauup}({\ft})=\zt_{\tauup}({\ft})\oplus\st_{\tauup}({\ft}),
\end{equation}
where $\rn({\ft})$ is the ideal 
of the nilpotent elements of $\rad({\ft})$, 
and coincides with its nilradical, and $\Li_{\tauup}({\ft})$ is its
$\tauup$-invariant reductive Levi factor, which further decomposes
into the direct sum of its center
$\zt_{\tauup}({\ft})$ and 
its semisimple ideal $\st_{\tauup}({\ft})$. Let
$A$ be any element
of maximal rank in $\zt_{\tauup}({\ft})\cap\mathfrak{h}_{\mathbb{R}}$. 
Then the isotropy at the base point $p_0$ of the transitive action of
$\mathbf{U}_0$ on $F$ is the group
\begin{equation}
  \label{eq:43i}
  \mathbf{F}_0=\mathbf{F}\cap\mathbf{U}_0=\{g\in\mathbf{U}_0\mid
\Ad_{\mathfrak{g}}(g)(A)=A\}
=\Z_{\mathbf{U}_0}(\zt_{\tauup}({\ft})).
\end{equation}
We observe that $F$ can be viewed as an $\mathbf{U}_0$-homogeneous
totally complex $CR$ manifold, with isotropy $\mathbf{F}_0$ and
$CR$ algebra $(\mathfrak{u}_0,{\ft})$ at the base point.
\subsection{Orbits of a real form}
Let $\mathbf{G}_0$ be a real form of $\mathbf{G}$, 
$\mathfrak{g}_0=\Lie(\mathbf{G}_0)$, and denote by  $\sigmaup(X)=\bar{X}$
the conjugation in $\mathfrak{g}$ defined by $\mathfrak{g}_0$. \par
We consider the $\mathbf{G}_0$-orbit
$E$ through the base point $p_0$ of $F=\mathbf{G}/\mathbf{F}$. It is  a 
$\mathbf{G}_0$-homogeneous $CR$ manifold, having at the base point
$p_0$ isotropy 
\begin{equation}
  \label{eq:42x}
  \mathbf{E}_0=\mathbf{F}\cap\mathbf{G}_0=\N_{\mathbf{G}_0}({\ft})=
\{g\in\mathbf{G}_0\mid \Ad_{\mathfrak{g}}(g)({\ft})={\ft}\}.
\end{equation}
and $CR$ algebra $(\mathfrak{g}_0,{\ft})$. With this structure, 
$E$ is a 
generic $CR$-submanifold of $F$.
\subsection{Mostow fibration} 
A pair $(\thetaup,\mathfrak{h}_0)$, consisting of a Cartan involution 
$\thetaup$ of $\mathfrak{g}_0$ and a $\thetaup$-invariant Cartan
subalgebra $\mathfrak{h}_0$ of $\mathfrak{g}_0$
contained in ${\ft}$ was called in \cite{AMN06b} 
\emph{adapted} to $(\mathfrak{g}_0,{\ft})$. If $(\thetaup,
\mathfrak{h}_0)$ is adapted, the involutions $\thetaup$ and
$\sigmaup$ commute and their composition 
$\tauup=\sigmaup\circ\thetaup$ is the conjugation with respect to
a compact form $\mathfrak{u}_0$ of $\mathfrak{g}$, for which
we have \eqref{eq:43}.
\par
The complexification $\mathfrak{h}$ of $\mathfrak{h}_0$ is a Cartan 
subalgebra of $\mathfrak{g}$ contained in ${\ft}$. Let 
$\mathcal{R}\subset\mathfrak{h}^*_{\mathbb{R}}$ be the root system
of $(\mathfrak{g},\mathfrak{h})$,
$\mathfrak{g}^{\alpha}=\{X\in\mathfrak{g}\mid [H,X]=\alpha(H)X,\;\forall
H\in\mathfrak{h}_{\mathbb{R}}\}$ the corresponding root spaces.
There is an element 
$A\in\zt_{\tauup}({\ft})\cap\mathfrak{h}_{\mathbb{R}}$ such that
\begin{equation}
  \label{eq41}
  \rn({\ft})={\sum}_{\alpha(A)>0}\mathfrak{g}^{\alpha},\;\;
\Li_{\tauup}({\ft})=\mathfrak{h}\oplus{\sum}_{\alpha(A)=0}\mathfrak{g}^{\alpha}.
\end{equation}
The set $\mathcal{F}=\{\alpha\in\mathcal{R}\mid\alpha(A)\geq{0}\}=
\{\alpha\in\mathcal{R}\mid\mathfrak{g}^{\alpha}\subset{\ft}\}$ is
the parabolic set of ${\ft}$. It  decomposes into the disjoint
union $\mathcal{F}=\mathcal{F}_n\cup\mathcal{F}_r$, with
$\mathcal{F}_n=\{\alpha\mid\alpha(A)>0\}$ and
$\mathcal{F}_r=\{\alpha\mid\alpha(A)=0\}=\{\alpha\in\mathcal{F}\mid
-\alpha\in\mathcal{F}\}$.\par
Consider the Cartan decomposition associated to $\thetaup$:
\begin{equation}
  \label{eq:51}
  \mathfrak{g}_0={\kt}_0\oplus\pt_0,\quad \kt_0=\{X\in\mathfrak{g}_0\mid
\thetaup(X)=X\},\;\; \pt_0=\{X\in\mathfrak{g}_0\mid
\thetaup(X)=-X\}.
\end{equation}
Then $\mathbf{K}_0=\mathbf{U}_0\cap\mathbf{G}_0$ has Lie algebra
$\kt_0$ and is a maximal compact subgroup of $\mathbf{G}_0$.
\par
In \cite[\S{10}]{AMN06b} we noticed that $E$ admits a Mostow fibration
(see \cite{Most55, Most62}),
whose basis $M$ is the
$\mathbf{K}_0$ orbit through the base point $p_0$; it is a
$\mathbf{K}_0$-homogeneous $CR$ manifold, with isotropy
\begin{equation}\label{eq:40}
  \mathbf{M}_0=\{g\in\mathbf{K}_0\mid \Ad_{\kt}(g)({\ft})={\ft}\}
=\{g\in\mathbf{K}_0\mid \Ad_{\kt}(g)(A)=A\},
\end{equation}
and $CR$ algebra $({\kt}_0,\mathfrak{v})$, with 
$\mathfrak{v}={\ft}\cap{\kt}$, at $p_0$. \par
We observe
that $\mathbf{M}_0$ is maximally compact in $\mathbf{E}_0$, and is
its deformation retract. In particular, 
$\mathbf{M}_0$ and $\mathbf{E}_0$ have the same number
of connected components (see e.g. \cite[\S{10}]{AMN06b}). We recall that
$\mathbf{M}_0$ is connected when $M$ is a minimal
orbit (see \cite[\S{8}]{AMN06}).\par
The $\mathbf{K}_0$-orbit $M$ in  $F$  
is the intersection of a $\mathbf{K}$-orbit $E^*$ 
with the $\mathbf{G}_0$-orbit $E$ in $F$ 
corresponding to $E^*$ in 
the Matsuki duality (cf. \cite{Mats88}).
The $\mathbf{K}$-orbit $E^*$ provides
a \textsl{generic} $CR$-embedding and a $\mathbf{K}$-realization of $M$.
In general,
the $CR$-embedding $M\hookrightarrow{F}$
is not generic, unless $E$ is the minimal orbit: indeed in this case
$M=E$ and $E^*$ is open in $F$ (see \cite{wolf69, AMN06}).
\subsection{The base of the Mostow fibration}
In this subsection we will show that $M$ is $\mathfrak{n}$-reductive,
i.e. that $M$ is a minimal $\mathbf{K}_0$-orbit in $E^*$, 
and investigate some features of 
the algebra $\mathfrak{v}={\ft}\cap{\kt}$ which describes
its $CR$ structure. 
\smallskip\par
Having fixed  an adapted $(\thetaup$,
$\mathfrak{h}_0)$, we use \eqref{eq:51} to split
$\mathfrak{h}_0$ into its toroidal and vectorial parts:
\begin{equation}
  \label{eq:52}
  \mathfrak{h}_0=\mathfrak{h}_0^+\oplus\mathfrak{h}_0^-,\quad\text{with}\quad
\mathfrak{h}_0^+=\mathfrak{h}_0\cap{\kt}_0,\quad
\mathfrak{h}_0^-=\mathfrak{h}_0\cap\pt_0.
\end{equation} \par
Then $\mathfrak{h}_{\mathbb{R}}=\mathfrak{h}_0^-\oplus{i}\mathfrak{h}^+_0$. 
The involutions $\sigmaup,\thetaup,\tauup$ act on 
$\mathfrak{h}_{\mathbb{R}}$, yielding by duality involutions
$\sigmaup^*,\thetaup^*,\tauup^*$ on $\mathfrak{h}_{\mathbb{R}}^*$,
which restrict to permutations of $\mathcal{R}$.
Set $\sigmaup^*(\alpha)=\bar{\alpha}$. Then
\begin{gather*}
  \sigmaup^*(\alpha)=\bar{\alpha},\;\;\thetaup^*(\alpha)=-\bar\alpha,\;\;
\tauup^*(\alpha)=-\alpha,\qquad \\
  \sigmaup({\mathfrak{g}^{\alpha}})=\mathfrak{g}^{\bar\alpha},\quad 
\thetaup(\mathfrak{g}^{\alpha})=\mathfrak{g}^{-\bar\alpha},\quad
\tauup(\mathfrak{g}^{\alpha})=\mathfrak{g}^{-\alpha},\quad \forall
\alpha\in\mathcal{R}.
\end{gather*}
We use the notation
\begin{equation*}
  \mathcal{R}_{\mathrm{re}}=\{\alpha\in\mathcal{R}\mid \bar\alpha=\alpha\},\quad
 \mathcal{R}_{\mathrm{im}}=\{\alpha\in\mathcal{R}\mid \bar\alpha=-\alpha\},\quad
 \mathcal{R}_{\mathrm{cx}}=\{\alpha\in\mathcal{R}\mid \bar\alpha\neq\pm\alpha\},
\end{equation*}
for the sets of \textsl{real}, \textsl{imaginary} and \textsl{complex} 
(meaning neither real nor imaginary) roots.
Imaginary roots are further subdivided 
into  \textsl{compact} and \textsl{noncompact}:
\begin{equation*} 
\mathcal{R}_{\mathrm{im}}= \mathcal{R}_{\mathrm{cp}}\cup \mathcal{R}_{\mathrm{nc}},
\quad\text{with}\quad \mathcal{R}_{\mathrm{cp}}=
\{\alpha\in\mathcal{R}_{\mathrm{im}}\mid\mathfrak{g}^{\alpha}\subset{\kt}\}, 
\quad\mathcal{R}_{\mathrm{nc}}
=\{\alpha\in\mathcal{R}_{\mathrm{im}}\mid\mathfrak{g}^{\alpha}\subset\pt\}.
\end{equation*}

Let 
\begin{equation}
\piup:\mathfrak{g}\ni X \longrightarrow \tfrac{1}{2}(X+\thetaup(X))
\in{\kt}  .
\end{equation}
be the projection on ${\kt}$ along $\pt$.
\begin{lem} \label{lem:73} There is a set $\mathcal{R}^+$ 
of positive roots with the property that
\begin{equation}\label{eq:75}
  \bar\alpha\in\mathcal{R}^+,\;
\forall \alpha\in\mathcal{R}^+\cap\mathcal{R}_{\mathrm{cx}},
\end{equation}
yielding  a direct sum decomposition of ${\kt}$:
\begin{equation}
  \label{eq:76}
  {\kt}=\mathfrak{h}^+
\oplus{\sum}_{\alpha\in\mathcal{R}^+\cap\mathcal{R}_{\mathrm{cp}}}
\mathfrak{g}^{-\alpha}
\oplus{\sum}_{\alpha\in\mathcal{R}^+\setminus\mathcal{R}_{\mathrm{nc}}}
\piup(\mathfrak{g}^{\alpha}).
\end{equation}
\end{lem}\begin{proof}
The existence of a set of positive roots satisfying \eqref{eq:75}
is, e.g., a consequence of \cite[Lemma\,4.3]{AMN06b}. Then \eqref{eq:76}
follows because
$\piup(\mathfrak{g}^{\alpha})=
\piup(\mathfrak{g}^{-\bar{\alpha}})$ for all \mbox{$\alpha\in\mathcal{R}$}
and $\piup(\mathfrak{g}^{\alpha})=\{0\}$ for $\alpha\in\mathcal{R}_{\mathrm{nc}}$.
\end{proof}
\begin{lem}\label{lem:82} Let $\mathcal{F}$ be the parabolic set of \;
${\ft}$ \; and $\mathcal{F}^*=\mathcal{F}\setminus\mathcal{R}_{\mathrm{nc}}$. Then
\begin{equation}
  \label{eq:qk1}
 \mathfrak{v}\! =\! {\ft}\cap{\kt}
\! =\!\mathfrak{h}^+\oplus{\sum}_{\alpha\in\mathcal{F}^{\thetaup}}
\piup(\mathfrak{g}^{\alpha}),
\;\; \text{with $ \mathcal{F}^{\thetaup}\!=\!
\{\alpha\in\mathcal{F}^*
\mid -\bar\alpha\in\mathcal{F}\}=\mathcal{F}^*\cap\thetaup^*(\mathcal{F}^*)$.}
\end{equation}
\end{lem}
\begin{proof}
An $X\in{\ft}$ uniquely decomposes as
\begin{equation*}
  X=H+{\sum}_{\alpha\in\mathcal{F}}X_{\alpha},\quad\text{with
$H\in\mathfrak{h}$ and $X_{\alpha}\in\mathfrak{g}^{\alpha}$}.
\end{equation*}
If $\thetaup(X)=X$, from
\begin{equation*}
  X=\thetaup(X)=\thetaup({H})+{\sum}_{\alpha\in\mathcal{F}}\thetaup({X}_{\alpha})
\end{equation*}
we obtain 
\begin{equation*}
  X_{\alpha}\neq{0}\;\Longrightarrow \; -\bar\alpha\in\mathcal{F} 
\quad\text{and}\quad
\thetaup(X_{\alpha})=X_{-\bar{\alpha}}\in{\ft}.
\end{equation*}
This yields \eqref{eq:qk1}.
\end{proof}
Split $\mathcal{F}^{\thetaup}$ into the two subsets
\begin{equation}
  \label{eq:712}
  \mathcal{F}^{\thetaup}_n=\mathcal{F}^{\thetaup}\cap\mathcal{F}_n,\quad
\mathcal{F}^{\thetaup}_r=\mathcal{F}_r\cap\thetaup^*(\mathcal{F}_r)
=\{\alpha\in\mathcal{F}_r\mid
-\bar{\alpha}\in\mathcal{F}_r\}.
\end{equation}
\begin{lem}
  $\mathcal{F}^{\thetaup},\mathcal{F}^{\thetaup}_r,\mathcal{F}^{\thetaup}_n$ are
closed set of roots.\qed
\end{lem}
\begin{lem} Let
$\rn({\ft})$ be the ideal of the nilpotent elements of the radical
of ${\ft}$. Then
\begin{equation}
  \label{eq:qk2}
  \nr(\mathfrak{v})
=\{X+\thetaup(X)\mid X\in\rn({\ft}),\;\thetaup(X)\in
{\ft}\}
={\sum}_{\alpha\in\mathcal{F}^{\thetaup}_n}\piup(\mathfrak{g}^{\alpha}).
\end{equation}\end{lem}
\begin{proof} Set $\mathfrak{y}=\{X+\thetaup(X)\mid X\in\rn({\ft}),\;
\thetaup(X)\in
{\ft}\}$.
 First we
prove that $\mathfrak{y}$ is an ideal in $\mathfrak{v}$. In fact, we have
\begin{align*}
  [X+\thetaup(X),Y+\thetaup(Y)]&=[X,Y]+\thetaup([X,Y])+[X,\thetaup(Y)]+
\thetaup([X,\thetaup(Y)]),
\end{align*}
with $[X,Y]$ and $[X,\thetaup(Y)]$ both belonging to $\rn({\ft})$.
\par
If $X+\thetaup(X)$ is an element of $\mathfrak{y}$, then
$[X,\thetaup(X)]\in\rn({\ft})$. We need to prove that also
$X+\thetaup(X)$ is nilpotent. To this aim, we observe that
$(\ad(X+\thetaup(X)))^m(Y)$ is a sum of terms of the form
$[Z_1,\hdots,Z_m,Y]$ with $Z_i\in\{X,\thetaup(X)\}$ for $i=1,\hdots,m$.
\par
Since $\rn({\ft})$ is an algebra of nilpotent elements of 
$\mathfrak{g}$, there is some integer $\nu$ such that 
$[X_1,\hdots,X_h,Y]=0$ for all 
$Y\in\mathfrak{g}$ if $X_1,\hdots,X_h\in\rn({\ft})$ and
$h\geq\nu$. We claim that
\begin{equation}\tag{$*$} \label{eq:99}
  [X_1,\hdots,X_m,Y]=0,\;\forall Y\in\mathfrak{g}\;\;\text{if}\;\;
X_1,\hdots,X_m\in{\ft}\;\;\text{and}\;\;
\#\{i\mid X_i\in\rn({\ft})\}\geq{\nu}.
\end{equation}
To prove our claim,
we argue by
recurrence on $m$. In fact, when $m=\nu$, \eqref{eq:99} is clearly
satisfied. Assume now that $m>\nu$ and that \eqref{eq:99} is true for
commutators of a lesser number of elements.
Clearly $[X_1,\hdots,X_m,Y]=0$
when all $X_i$'s belong to $\rn({\ft})$. Take now
$X_1,\hdots,X_m\in{\ft}$ with at least $\nu$ elements in
$\rn({\ft})$ and at least one, say $X_{i_0}$, in
${\ft}\setminus\rn({\ft})$. The commutator
$[X_1,\hdots,X_m,Y]$ is a linear combination of the commutators
$[X_1,\hdots,X_{i_0-1},X_{i_0+1},\hdots,[X_{i_0},Y]]$ and
$[X_1,\hdots,[X_{i_0},X_j],\hdots,Y]$. All these terms are zero
by the recursive assumption. \par
This shows that the terms of the form $[Z_1,\hdots,Z_m,Y]$
with at least $\nu$ of the $Z_i's$ equal to $X$ are zero for any
choice of $Y$. \par
Using the fact that $\thetaup$ is an automorphism of $\mathfrak{g}$,
it follows that also the terms $[Z_1,\hdots,Z_m,Y]$ for which
at least $\nu$ elements equal $\thetaup(X)$ are zero,
showing that $(\ad(X+\theta(X)))^{2\nu-1}=0$. \par
This proves that $\mathfrak{y}$ is an ideal of nilpotent elements 
in $\mathfrak{v}$. Hence, $\mathfrak{y}\subset\nr(\mathfrak{v})$.
From
\begin{equation}
  \label{eq:714}
  \mathfrak{y}={\sum}_{\alpha\in\mathcal{F}^{\thetaup}_n}\piup(\mathfrak{g}^{\alpha}),
\end{equation}
and the decomposition \eqref{eq:qk1} of $\mathfrak{v}$
in Lemma\,\ref{lem:82} 
we obtain that
\begin{equation}
  \label{eq:713}
  \mathfrak{L}(\mathfrak{v})=\mathfrak{h}^+\oplus
{\sum}_{\alpha\in\mathcal{F}^{\thetaup}_r}\piup(\mathfrak{g}^{\alpha})=\mathfrak{v}\cap
\bar{\mathfrak{v}}
\end{equation}
is a reductive complement of $\mathfrak{y}$ in $\mathfrak{v}$. 
Hence $\mathfrak{y}=\nr(\mathfrak{v})$. \end{proof}
We summarize the results obtained so far in this section by the following
\begin{thm}\label{thm:AAA}
Let $(\mathfrak{g}_0,{\ft})$ be a $CR$ algebra with
$\mathfrak{g}_0$ semisimple and ${\ft}$ parabolic in the
complexification $\mathfrak{g}$ of $\mathfrak{g}_0$, 
$(\thetaup,\mathfrak{h}_0)$ an adapted pair, 
$\kt_0=\{X\in\mathfrak{g}_0\mid \thetaup(X)=X\}$
the corresponding maximal compact subalgebra of
$\mathfrak{g}_0$, 
${\kt}$ its complexification and
$\mathfrak{v}={\ft}\cap{\kt}$. Then 
$\mathfrak{v}\cap\bar{\mathfrak{v}}$ is a reductive complement of
$\nr(\mathfrak{v})$  in $\mathfrak{v}$. \par
Hence, the $\mathbf{K}_0$-orbit $M$ of \; ${\ft}$\;  in $\mathfrak{g}$ is
an $\mathfrak{n}$-reductive $CR$ manifold. \qed
\end{thm}
We also have
\begin{prop}
  The closed Lie group $\mathbf{V}=\mathbf{F}\cap\mathbf{K}$ admits
a Levi-Chevalley decomposition
\begin{equation}
  \label{eq:756}
  \mathbf{V}=\mathbf{V}_n\mathbf{V}_r,\;\;\text{with}\;\;
\Lie(\mathbf{V}_n)=\nr(\mathfrak{v}), \;\;
\Lie(\mathbf{V}_r)=\Li(\mathfrak{v})=\mathfrak{m}.
\end{equation}
\end{prop}
\begin{proof} 
Indeed $\mathbf{G}$ is an algebraic group and $\mathbf{V}$ is an
algebraic subgroup of $\mathbf{G}$, hence admitting a
Levi-Chevalley decomposition \eqref{eq:756} (see e.g. \cite[Ch.6]{OV90}).
\end{proof}

\begin{exam}\label{ex:b}
Consider the minimal orbit $M$ of $\mathbf{SU}(2,3)$ in the complex flag
manifold of lines and $3$-planes of $\mathbb{C}^5$:
\begin{equation*}
 M^{(2,6)}=\{\text{flags $(L_1\subset{L}_3)$ with $L_1\subset{L}_3^\perp\subset{L}_3$}
\}, 
\end{equation*}
where the perpendicularity is taken with respect to a Hermitian symmetric
form $\mathrm{h}$ of signature $(2,3)$ in $\mathbb{C}^5$. \par
This is the minimal orbit described by the cross marked Satake diagram
\bigskip
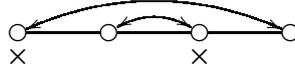
\begin{equation*}
  \xymatrix@R=-.3pc{\!\!\medcirc\!\!\ar@{-}[r]\ar@{<->}@/^1pc/[rrr]
&\!\!\medcirc\!\!
\ar@{-}[r]\ar@{<->}@/^/[r]&\!\!\medcirc\!\! \ar@{-}[r]&\!\!\medcirc\!\!\\
\times&&\times}
\end{equation*}
Let $\mathrm{h}$ be defined by the anti-diagonal matrix 
$\left( \begin{smallmatrix}
    &&&&1\\
&&&1\\
&&1\\
&1\\
1
  \end{smallmatrix}\right)$.
\par\smallskip
We give below a matrix description of the complexification $\kt$ of
the maximal compact subalgebra 
$\kt_0=\mathfrak{s}(\mathfrak{u}_2\times\mathfrak{u}_3)$ 
of $\mathfrak{su}_{2,3}$, of
$\mathfrak{v}={\ft}\cap\kt$, and of its normalizer in~$\kt$: 
\begin{align*}
  {\kt}&:
  \begin{pmatrix}
    \lambda_2&\zeta_1&\eta_1&\eta_2&\mu_2\\
z_1&\lambda_1&z_2&\mu_1&w_2\\
w_1&\zeta_2&\lambda_0&\zeta_2&w_1\\
w_2&\mu_1&z_2&\lambda_1&z_1\\
\mu_2&\eta_2&\eta_1&\zeta_1&\lambda_2
  \end{pmatrix},&\mathfrak{v}=\kt\cap{\ft}&:
\begin{pmatrix}
    \lambda_2&\zeta_1&0&0&0\\
0&\lambda_1&0&0&0\\
0&\zeta_2&\lambda_0&\zeta_2&0\\
0&0&0&\lambda_1&0\\
0&0&0&\zeta_1&\lambda_2
  \end{pmatrix}.
\end{align*}
The normalizer  
\begin{equation*}
  \mathrm{N}_{{\kt}}({\kt}\cap{\ft}):
  \begin{pmatrix}
      \lambda_2&\zeta_1&0&0&\mu\\
0&\lambda_1&0&\mu&0\\
0&\zeta_2&\lambda_0&\zeta_2&0\\
0&\mu&0&\lambda_1&0\\
\mu&0&0&\zeta_1&\lambda_2 
  \end{pmatrix}
\end{equation*}
of ${\kt}\cap{\ft}$ 
does not contain any maximal torus in 
${\kt}$, and hence
$\mathfrak{v}$
is not \textsl{regular}.
This example shows that we cannot expect, in general, that the basis $M$
of the Mostow fibration of a $\mathbf{G}_0$-orbit has a
$\mathfrak{t}_0$-regular $\mathfrak{v}$. \par

\end{exam}
We have indeed
\begin{prop} Let $\mathfrak{v}={\ft}\cap\kt$ and set
$\Li(\mathfrak{v})=\mathfrak{v}\cap\bar{\mathfrak{v}}$. Then the
following are equivalent:
\begin{enumerate}
\item \label{t1} $(\kt_0,\mathfrak{v})$ is of type $\mathrm{I\! I}$, i.e.
$\Li(\mathfrak{v})$ is normalized by a maximal torus $\mathfrak{t}_0$ 
of $\kt_0$;
\item \label{t3}
there exists a maximal set of strongly orthogonal real roots 
$\{\alpha_1,\hdots,\alpha_r\}$ in
$\mathcal{F}$ such that
\begin{equation}
  \label{eq:afx}
  \beta\in\mathcal{F}_r^{\thetaup},\;\;1\leq{i}\leq{r},\;\; \beta+\alpha_i\in
\mathcal{R}\;\Longrightarrow \beta+\alpha_i\in\mathcal{F}_r^{\thetaup}.
\end{equation}
\item \label{t4} if $(\thetaup,\mathfrak{h}_0)$ is an adapted pair
with $\mathfrak{h}_0\subset{\ft}$ 
maximally compact, then there exists a maximal system
of strongly orthogonal real roots $\{\alpha_1,\hdots,\alpha_r\}$ in
$\mathcal{F}_n$ such that
\begin{equation}
  \label{eq:afy}
\beta+\alpha_i\notin\mathcal{R},\;\;\forall \beta\in\mathcal{F}_r^{\thetaup},\;
1\leq{i}\leq{r}.
\end{equation}
\end{enumerate}
The following are equivalent
\begin{enumerate}
\item[(a)]\label{ta} $(\kt_0,\mathfrak{v})$ is of type $\mathrm{I}$, i.e.
$\mathfrak{v}$ is normalized by a maximal torus $\mathfrak{t}_0$ 
of $\kt_0$;
\item[(b)]\label{tb} there is a maximal system of strongly orthogonal 
real roots $\{\alpha_1,\hdots,\alpha_r\}$ in $\mathcal{F}$ such that
\begin{equation}
  \label{eq:afz}
  1\leq{i}\leq{r},\;\;\beta\in\mathcal{F}^\thetaup,\;\;
\alpha_i+\beta\in\mathcal{R}\Longrightarrow
\alpha_i+\beta\in\mathcal{F}^{\thetaup};
\end{equation}
\item[(c)]\label{tcx} if $(\thetaup,\mathfrak{h}_0)$ is an adapted pair
with $\mathfrak{h}_0\subset{\ft}$ 
maximally compact,
then there
exists  a maximal system of strongly orthogonal 
real roots $\{\alpha_1,\hdots,\alpha_r\}$ in $\mathcal{F}_n$ such that
\begin{equation}
  \label{eq:afw}
  1\leq{i}\leq{r},\;\;\beta\in\mathcal{F}^\thetaup,\;\;
\alpha_i+\beta\in\mathcal{R}\Longrightarrow
\alpha_i+\beta\in\mathcal{F}^{\thetaup}.
\end{equation}
\end{enumerate}
\par
In particular, $\mathfrak{v}$ is regular in $\kt$ 
when
${\ft}\cap\mathfrak{g}_0$ contains 
a maximal torus of $\mathfrak{g}_0$.
\end{prop} 
\begin{proof} Fix a maximally compact Cartan subalgebra $\mathfrak{h}_0$
adapted to $(\mathfrak{g}_0,{\ft})$. Then $\mathcal{F}_r$ does not
contain any real root and
every $\thetaup$-invariant maximal torus of $\mathfrak{g}_0$ can be
obtained by applying to $\mathfrak{h}_0$ the Cayley transform with respect to
a maximal set $\{\alpha_1,\hdots,\alpha_r\}$ of strongly orthogonal real roots
contained in $\mathcal{F}_n$ (see \cite{Kn:2002}). If $\beta\in\mathcal{F}$
and $\beta+\alpha_i\in\mathcal{R}$, then $\beta+\alpha_i\in\mathcal{F}_n$.
This gives the equivalence \eqref{t1}$\Leftrightarrow$\eqref{t4}.
We prove the equivalence  \eqref{t3}$\Leftrightarrow$\eqref{t4} using
the Cayley transform with respect to a maximal system of strongly 
orthogonal real roots in $\mathcal{F}_r$.\par\smallskip
Assume that (c) holds true. We want to prove that we
have also \eqref{t4}. Indeed, assume by contradiction that there are
$\beta\in\mathcal{F}^{\thetaup}$ and 
$1\leq{i}\leq{r}$, with  $\alpha_i+\beta\in\mathcal{R}$. 
Since $\pm\beta,\pm\bar{\beta}\in\mathcal{F}$, we have
$\alpha_1+\beta,\alpha_i+\bar\beta\in\mathcal{F}_n$. By the assumption
that $\alpha_1+\beta,\alpha_i+\bar\beta\in\mathcal{F}^\thetaup$ we
get a contradiction, because also 
$\thetaup^*(\alpha_1+\beta)=-\alpha_1-\bar\beta,\thetaup^*(\alpha_i+\bar\beta)
=-\alpha_i-\beta\in\mathcal{F}$, and this gives in particular
$\alpha_i+\beta\in\mathcal{F}_r$.\par
The equivalence (a)$\Leftrightarrow$(b)$\Leftrightarrow$(c) is then proved
with arguments similar to those used to show that 
\eqref{t1}$\Leftrightarrow$\eqref{t3}$\Leftrightarrow$\eqref{t4}.
\end{proof}
In Example\,\ref{ex:b}, denoting by $\alpha_1,\hdots,\alpha_4$ the simple
roots in the Satake diagram, we have
\begin{gather*}
  \mathcal{F}_n=\{\alpha_1,\alpha_1\! + \! \alpha_2,\alpha_1\! + \! 
\alpha_2\! + \! \alpha_3,
\alpha_1\! + \! \alpha_2\! + \! \alpha_3\! + \! \alpha_4,\alpha_2\! + 
\! \alpha_3, \alpha_2\! + 
\! \alpha_3\! + \! \alpha_4,
\alpha_3,\alpha_3\! + \! \alpha_4\},\\
\mathcal{F}_r=\{\pm\alpha_2,\pm\alpha_4\},\; 
\mathcal{F}_n^{\thetaup}=\{\alpha_1,\alpha_3\},\;
\mathcal{F}_r^{\thetaup}=\emptyset,\;  \mathcal{F}_n
\cap\mathcal{R}_{\mathrm{re}}
=\{\alpha_2\! + \! \alpha_3, 
\alpha_1\! + \! \alpha_2\! + \! \alpha_3\! + \! \alpha_4\},\\
\alpha_1+(\alpha_2\! + \! \alpha_3)\in\mathcal{R},\quad 
-\overline{\alpha_1+\alpha_2+\alpha_3}=-(\alpha_2+\alpha_3+\alpha_4)\notin
\mathcal{F}\Rightarrow \alpha_1\! +\!\alpha_2\! + \! \alpha_3\notin
\mathcal{F}^{\thetaup}.
\end{gather*}
Since $\mathcal{F}_n\cap\mathcal{R}_{\mathrm{re}}$ is the
unique maximal system of strongly orthogonal real roots in $\mathcal{F}$,
this shows that $\mathfrak{v}$
does not satisfy (b). 
\begin{exam}
  Let $M$ be the minimal orbit of $\mathbf{SO}(3,5)$ in the flag $F$ of
the $3$-planes contained in the complex quadric of $\mathbb{CP}^7$,
associated to the
cross-marked Satake diagram
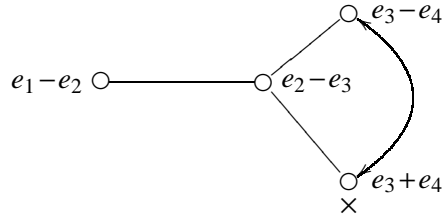
\begin{equation*}
  \xymatrix@R=1pc{&&& \!\!\medcirc \!\!\ar@{<->}@/^2pc/[dd] &
\!\!\!\!\!\!\!\!\!\!\!\!\!e_3\! - \! e_4\\
e_1\! -\! e_2\;\medcirc\!\! \ar@{-}[rr]&&\!\!\medcirc\!\! \ar@{-}[ur]\ar@{-}[dr]
&\!\!\!\!\!\!\!\!\!\!\!\!\! e_2\! -\! e_3\\
&&& \!\! {\begin{matrix}
\\[-4pt]
\medcirc\\[-4pt]
\times\end{matrix}} \!\! & \!\!\!\!\!\!\!\!\!\!\!\!\!
e_3\! + \! e_4}
\end{equation*}
We have $\mathcal{F}=\mathcal{F}_n\cup\mathcal{F}_r$ and
$\mathcal{F}^\thetaup=\mathcal{F}_n^\thetaup\cup\mathcal{F}_r^\thetaup$
with
\begin{align*}
  \mathcal{F}_n&=\{e_i+e_j\mid 1\leq{i}<j\leq{4}\}, &
\mathcal{F}_r&=\{\pm(e_i-e_j)\mid 1\leq{i}<j\leq{4}\},\\
\mathcal{F}_n^{\thetaup}&=\{e_i+e_4\mid 1\leq{i}\leq{3}\}, &
\mathcal{F}_r^{\thetaup}&=\{\pm(e_i-e_j)\mid 1\leq{i}<j\leq{3}\}.
\end{align*}
The maximal sets of strongly orthogonal real roots are the sets
$E_{i,j}=\{e_i\pm{e}_j\}$ for $1\leq{i}<j\leq{3}$. Then condition \eqref{t3}
is not satisfied and $\Li(\mathfrak{v})$ is not $\mathfrak{t}_0$-regular
for any maximal torus $\mathfrak{t}_0$ in $\kt_0\simeq\mathfrak(\mathfrak{so}_3\oplus
\mathfrak{so}_5)$.

\end{exam}
\section{Parabolic regularization} 
\label{sec:rg}
In this section, ${\kt}$ is a reductive complex Lie algebra.
Starting from any complex Lie subalgebra $\mathfrak{v}$ of ${\kt}$,
we describe a canonical \emph{regularization} procedure to associate to
$\mathfrak{v}$ a suitable \textsl{parabolic} subalgebra 
$\mathfrak{w}$ of ${\kt}$.
\begin{lem}
  Let $\mathfrak{v}$ be a Lie subalgebra of a reductive Lie algebra
${\kt}$. Define by recurrence
\begin{equation}
  \label{eq:41b}
  \begin{cases}
    \mathfrak{v}_0=\mathfrak{v},\\
\mathfrak{v}_{m+1}=\mathrm{N}_{{\kt}}(\nr(\mathfrak{v}_m)) &\text{for 
$m\geq{0}$},
  \end{cases}
\end{equation}
where
$ \mathrm{N}_{{\kt}}(\nr(\mathfrak{v}_m))=\{X\in{\kt}\mid
[X,\nr(\mathfrak{v}_m)]\subset\nr(\mathfrak{v}_m)\}$
is the normalizer of $\nr(\mathfrak{v}_m)$ in~${\kt}$.
Then we obtain two non decreasing sequences of Lie subalgebras
\begin{gather} \label{eq:42b}
  \nr(\mathfrak{v})=\nr(\mathfrak{v}_0)\subset
\nr(\mathfrak{v}_1)\subset\cdots\subset
\nr(\mathfrak{v}_m)\subset
\nr(\mathfrak{v}_{m+1})\subset\cdots,\\  \label{eq:43b}
\mathfrak{v}=\mathfrak{v}_0\subset\mathfrak{v}_1\subset\cdots
\subset \mathfrak{v}_m\subset\mathfrak{v}_{m+1}\subset\cdots
\end{gather}
\end{lem}
\begin{proof}
Since each $\nr(\mathfrak{v}_m)$ is an ideal in $\mathfrak{v}_m$, 
from $[\mathfrak{v}_m,\nr(\mathfrak{v}_m)]\subset\nr(\mathfrak{v}_m)$ we have
$\mathfrak{v}_m\subset\mathfrak{v}_{m+1}$. \par
Moreover, $\nr(\mathfrak{v}_m)$, being a nilpotent ideal of $\mathfrak{v}_{m+1}$,
is contained in $\rad(\mathfrak{v}_{m+1})$. It consists of nilpotent elements
and therefore $\nr(\mathfrak{v}_m)\subset\nr(\mathfrak{v}_{m+1})$.
\end{proof}
We have (see \cite{Weis66, BT71, Hum}):
\begin{prop}\label{prop:52}
The sequences \eqref{eq:42b} and \eqref{eq:43b} stabilize and 
\begin{equation*}
\nr(\mathfrak{v}_m)=\nr(\mathfrak{v}_{m+1})\quad\Longrightarrow\quad
\mathfrak{v}_{m+1}=\mathfrak{v}_{h+1}\;\; 
\text{and}\;\;\nr(\mathfrak{v}_m)=\nr(\mathfrak{v}_{h}),
\;\;\forall h\geq{m}.
\end{equation*}
The Lie subalgebra
$
  \mathfrak{w}={\bigcup}_{m}{\mathfrak{v}_m}
$
is parabolic in ${\kt}$.
\end{prop}
\begin{proof}
The first assertions follows because \eqref{eq:42b} and
\eqref{eq:43b} are increasing sequences of linear subspaces of a finite dimensional
vector space. \par
Assume that $\mathfrak{w}=\mathfrak{v}_h$ for all $h\geq{m}$. 
Since $\nr(\mathfrak{w})$ is a characteristic ideal of $\mathfrak{w}$, 
we have
\begin{equation*}
  \mathfrak{w}\subset\mathrm{N}_{{\kt}}(\mathfrak{w})\subset
\mathrm{N}_{{\kt}}(\nr(\mathfrak{w}))=\mathfrak{w}.
\end{equation*}
Hence $\mathfrak{w}$ is its own normalizer and equals the normalizer of 
$\nr(\mathfrak{w})$. 
 Then the statement follows
from Lemma\,\ref{lem:41} below.
\end{proof}
We have (see e.g. \cite{Mr56, Weis66, Plat69, BT71, Hum})
\begin{lem}
  \label{lem:41}
Let $\mathfrak{w}$ be a Lie subalgebra of ${\kt}$, and assume that, for
$\mathfrak{n}=\nr(\mathfrak{w})$,
\begin{equation*}
  \mathfrak{w}=\mathrm{N}_{{\kt}}(\mathfrak{w})=\{X\in{\kt}\mid
[X,\mathfrak{w}]\subset\mathfrak{w}\}=
\mathrm{N}_{{\kt}}(\mathfrak{n})=\{X\in{\kt}\mid
[X,\mathfrak{n}]\subset\mathfrak{n}\}.
\end{equation*}
Then $\mathfrak{w}$ is parabolic in ${\kt}$.
\end{lem}
The proof of Lemma\,\ref{lem:41} utilizes the following lemma.
\begin{lem}\label{lem:43} Let $\mathfrak{w}$ be a complex Lie subalgebra
of ${\kt}$ and 
$\mathfrak{r}$ any solvable subalgebra of $\mathfrak{w}$. Then there exists 
an inner automorphism $\psiup$ of $\mathfrak{w}$ such that
\begin{equation*}
  \nr(\mathfrak{r}\cap\psiup(\mathfrak{r}))\subset\nr(\mathfrak{w}).
\end{equation*}
\end{lem}
\begin{proof} First assume that $\mathfrak{w}$ is reductive. Then,
taking any Borel subalgebra $\mathfrak{b}$ of $\mathfrak{w}$ 
containing $\mathfrak{r}$, we obtain the statement as
a consequence of the existence in $\mathfrak{w}$ 
of a Borel subalgebra opposite to
$\mathfrak{b}$ (see e.g. \cite{Bou75, Mi2002}).\par
In general, we note that the quotient $\mathfrak{w}/\nr(\mathfrak{w})$ is
reductive and apply the previous argument to its solvable subalgebra
$(\mathfrak{r}+\nr(\mathfrak{w}))/\nr(\mathfrak{w})$.
\end{proof}
\begin{proof}[Proof of Lemma\,\ref{lem:41}] 
Let $\mathfrak{b},\mathfrak{b}'$ be two Borel subalgebras of $\kt$
 with the
two properties
\begin{list}{-}{}
\item $\mathfrak{b}\cap\mathfrak{w}$ and $\mathfrak{b}'\cap\mathfrak{w}$
are Borel subalgebras of $\mathfrak{w}$;
\item $ \nr(\mathfrak{b}\cap\mathfrak{b}'\cap\mathfrak{w})\subset\mathfrak{n}$.
\end{list}
The existence of such a pair $(\mathfrak{b},\mathfrak{b}')$ is granted by
Lemma\,\ref{lem:43}.
Indeed, any Borel subalgebra $\mathfrak{r}$ of $\mathfrak{w}$ is the intersection
of  $\mathfrak{w}$ with a Borel subalgebra $\mathfrak{b}$ of $\kt$,
and the inner automorphism $\psiup$ of $\mathfrak{w}$ extends to an inner
automorphism of $\kt$, so that $\psiup(\mathfrak{r})=\psiup(\mathfrak{b})\cap
\mathfrak{w}$. \par
Since all Borel subalgebras of $\mathfrak{w}$ contain $\mathfrak{n}$, we
actually have $\nr(\mathfrak{b}\cap\mathfrak{b}'\cap\mathfrak{w})=\mathfrak{n}$.
\par
We claim that, in fact, $\nr(\mathfrak{b}\cap\mathfrak{b}')=\mathfrak{n}$. 
Indeed, 
$\N_{\nr(\mathfrak{b}\cap\mathfrak{b}')}(\mathfrak{n})\subset\mathfrak{w}$ and
hence
\begin{equation*}
  \mathfrak{n}\subset\N_{\nr(\mathfrak{b}\cap\mathfrak{b}')}(\mathfrak{n})=
\N_{\nr(\mathfrak{b}\cap\mathfrak{b}')}(\mathfrak{n})\cap\mathfrak{w}\subset
\nr(\mathfrak{b}\cap\mathfrak{b}')\cap\mathfrak{w}\subset
\nr(\mathfrak{b}\cap\mathfrak{b}'\cap\mathfrak{w})=\mathfrak{n}.
\end{equation*}
We obtain that $\nr(\mathfrak{b}\cap\mathfrak{b}')=\mathfrak{n}$, because
otherwise we would have 
$\mathfrak{n}\subsetneqq\N_{\nr(\mathfrak{b}\cap\mathfrak{b}')}(\mathfrak{n})$.
It follows that
$\mathfrak{b}\cap\mathfrak{b}'\subset\mathfrak{w}$ and, in particular, that
$\mathfrak{w}$ contains a Cartan subalgebra $\mathfrak{t}$ of $\kt$. 
Let $\mathcal{K}$ be the root system of $(\kt,\mathfrak{t})$, set
$\mathcal{K}^+=\{\alpha\in\mathcal{K}\mid \kt^{\alpha}\subset\mathfrak{b}\}$,
and denote by 
$\mathcal{B}$ the set of simple roots of $\mathcal{K}^+$. To prove that
$\mathfrak{w}$ is parabolic, it suffices to show that
$\kt^{\beta}\subset\mathfrak{w}$ for all $\beta\in\mathcal{B}$, since
this implies that $\mathfrak{b}\subset\mathfrak{w}$. Assume by contradiction
that $\kt^{\beta}\not\subset\mathfrak{w}$
for some $\beta\in\mathcal{B}$. Then $\kt^{\beta}\not\subset\mathfrak{b}'$.
Hence $\kt^{-\beta}\subset\mathfrak{b}'$. Consider the Borel subalgebra
\begin{equation*}
  \mathfrak{b}''=\mathfrak{t}\oplus\kt^{-\beta}\oplus
{\sum}_{\beta\in\mathcal{K}^+\setminus\{\beta\}}\kt^{\alpha}.
\end{equation*}
We have
$\mathfrak{b}\cap\mathfrak{w}\subset\mathfrak{b}''\cap\mathfrak{w}$
and the two sets are equal because $\mathfrak{b}\cap\mathfrak{w}$ is
maximal solvable in $\mathfrak{w}$. Then $\mathfrak{b}',\mathfrak{b}''$
are Borel subalgebras of $\kt$ whose intersections with $\mathfrak{w}$
are Borel subalgebras of $\mathfrak{w}$, and 
\begin{equation*}
 \nr(\mathfrak{b}'\cap\mathfrak{b}''\cap\mathfrak{w})=
 \nr(\mathfrak{b}'\cap\mathfrak{b}\cap\mathfrak{w})
\subset\mathfrak{n}.
\end{equation*}
This implies, by the previous argument, that 
$\kt^{-\beta}\subset\mathfrak{b}'\cap\mathfrak{b}''\subset\mathfrak{w}$.
But this in turns brings $\kt^{-\beta}\subset\mathfrak{n}$, contradicting
the inclusion $\mathfrak{n}\subset\mathfrak{b}$. 
The proof is complete.
\end{proof}
\begin{defn}
We call the Lie subalgebra $\mathfrak{w}$ of Proposition\;\ref{prop:52}
the \emph{parabolic regularization} of $\mathfrak{v}$ in ${\kt}$.  
\end{defn}
\section{The deployment theorem} \label{sec:dp}
Let $\mathbf{K}_0$ be a connected compact Lie group,
$M$ a $\mathbf{K}_0$-homogeneous $CR$ manifold, with isotropy subgroup 
$\mathbf{M}_0$ and $CR$ algebra $(\kt_0,\mathfrak{v})$ 
at a point $p_0\in{M}$. We set $\kt_0=\Lie(\mathbf{K}_0)$, 
$\kt=\mathbb{C}\otimes_{\mathbb{R}}\kt_0$, $\mathfrak{m}_0=\Lie(\mathbf{M}_0)$,
$\mathfrak{m}=\mathbb{C}\otimes_{\mathbb{R}}\mathfrak{m}_0$. Note that
$\mathfrak{m}\subset\mathfrak{v}\subset\kt$. \par
Let $\Par$ be the
set of all parabolic subalgebras of $\kt$. Recall that, for $\mathfrak{q}\in\Par$,
 the intersection 
 $\Li(\mathfrak{q})=\mathfrak{q}\cap\bar{\mathfrak{q}}$
is the reductive Levi factor
of a Levi-Chevalley decomposition of $\mathfrak{q}$
 (see Lemma\,\ref{lem:11}). 
For $\mathfrak{q}\in\Par$, we set 
\begin{equation}
  \label{eq:6x4a}
  \mathbf{N}_{\mathfrak{q}}=\{g\in\mathbf{K}_0\mid \mathrm{Ad}(g)(\mathfrak{q})
=\mathfrak{q}\}, \quad N_{\mathfrak{q}}=\mathbf{K}_0/
\mathbf{N}_{\mathfrak{q}}.
\end{equation}
With the $\mathbf{K}_0$-equivariant 
$CR$ structure $(\kt_0,\mathfrak{q})$ at
the base point, $N_{\mathfrak{q}}$ is a complex flag manifold of $\kt$.
Its complex structure only depends on the pair $(\kt,\mathfrak{q})$
(see e.g. \cite{wolf69}). 
Moreover, $N_{\mathfrak{q}}$ is simply connected and the isotropy subgroup
$\mathbf{N}_{\mathfrak{q}}$ is the centralizer of a torus 
of $\kt_0$. Set
\begin{equation*}
  \ses(\mathfrak{q})=[\Li(\mathfrak{q}),\Li(\mathfrak{q})],\quad
\zt(\mathfrak{q})=\Z(\Li(\mathfrak{q}))=
\{X\in\Li(\mathfrak{q})\mid [X,\Li(\mathfrak{q})]=\{0\}\}
\end{equation*}
for the semisimple ideal and the center of its
canonical reductive factor $\Li(\mathfrak{q})$.
Then $\mathbf{N}_{\mathfrak{q}}$ is the centralizer of the torus
$\zto(\mathfrak{q})=\zt(\mathfrak{q})\cap\kt_0$, and 
$\Lie(\mathbf{N}_{\mathfrak{q}})=\Lio(\mathfrak{q})=\Li(\mathfrak{q})\cap\kt_0$.
\par\smallskip
The obvious equivariant fibrations of complex flag manifolds yield the inclusions
\begin{equation}
  \label{eq:6xx}
  \mathbf{N}_{\mathfrak{q}_1}\subset\mathbf{N}_{\mathfrak{q}_2} \quad\text{if 
$\mathfrak{q}_1,\mathfrak{q}_2\in\Par$ \; and
\; $\mathfrak{q}_1\subset\mathfrak{q}_2$.}
\end{equation}
\par 
The existence of the parabolic regularization proves the following
\begin{prop} For all complex subalgebras
 $\mathfrak{v}$ of $\kt$
  the set
  \begin{equation}\label{eq:6x0}
      {\Par}(\mathfrak{v})=\{\mathfrak{q}\in{\Par}\mid
\mathfrak{v}\subset\mathfrak{q},\;\; \nr(\mathfrak{v})
\cap\mathfrak{L}(\mathfrak{q})=
\{0\}\}
  \end{equation}
is nonempty. \qed
\end{prop}
\begin{defn}
We say that a $\mathfrak{q}\in\Par(\mathfrak{v})$ is $\mathbf{M}_0$-admissible
if $\Ad_{\kt}(\mathbf{M}_0)(\mathfrak{q})=\mathfrak{q}$. \par
This is equivalent to the condition that 
$\mathbf{M}_0\subset\mathbf{N}_{\mathfrak{q}}$. 
Set
\begin{equation}
  \label{eq:8g}
  \Par(\mathfrak{v},\mathbf{M}_0)=\{\mathfrak{q}\in\Par(\mathfrak{v})\mid
\Ad_{\kt}(\mathbf{M}_0)(\mathfrak{q})=\mathfrak{q}\}.
\end{equation}
\end{defn}
\begin{rmk}
The condition $\mathbf{M}_0\subset\mathbf{N}_{\mathfrak{q}}$ is always satisfied
if $\mathbf{M}_0$ is connected, so that $\Par(\mathfrak{v},\mathbf{M}_0)=
\Par(\mathfrak{v})$ in this special case.\par
For instance, in Example\,\ref{ex:b},
the isotropy subgroup $\mathbf{M}_0$ at the base point  
consists of the diagonal matrices of the form
\begin{equation*}
\mathrm{diag}(\exp(it_1),\exp(it_2),\exp(it_3),\exp(it_2),\exp(it_1)),
\;\;\;\text{with $2t_1+2t_2+t_3=0$.}
\end{equation*}
and hence is connected. In fact, 
this is always the case when $M$ is
the minimal orbit of a real form in a complex flag manifold (see \cite{AMN06}).
\par\smallskip
If $\mathbf{M}_0^0$ is the
connected component of the identity in $\mathbf{M}_0$, then
$\tilde{M}=\mathbf{K}_0/\mathbf{M}_0^0$ is 
a $\mathbf{K}_0$-homogeneous
$CR$ manifold with the same $CR$ algebra $(\kt_0,\mathfrak{v})$ at all
points in the fiber of $p_0$. Hence the
$\mathbf{K}_0$-equivariant covering map $\tilde{\piup}:\tilde{M}\to{M}$
is a $CR$-submersion and a local $CR$-diffeomorphism with finite discrete fiber.
\end{rmk}
\begin{lem}\label{lem64}
If $\mathfrak{q}\in\Par(\mathfrak{v},\mathbf{M}_0)$, then
\begin{equation}
  \label{eq:6x11}\left\{ \begin{aligned}
  \Ad_{\kt}(\mathbf{M}_0)(\nr(\mathfrak{q}))&=\nr(\mathfrak{q}), &
\Ad_{\kt}(\mathbf{M}_0)(\Li(\mathfrak{q}))&=\Li(\mathfrak{q}),\\
\Ad_{\kt}(\mathbf{M}_0)(\ses(\mathfrak{q}))&=\ses(\mathfrak{q}), &
\Ad_{\kt}(\mathbf{M}_0)(\zt(\mathfrak{q}))&=\zt(\mathfrak{q}).
\end{aligned}\right.
\end{equation}
\end{lem} \begin{proof}
The first equality follows because $\nr(\mathfrak{q})$
is a characteristic ideal; the second 
because $\mathbf{M}_0\subset\mathbf{K}_0$, and $\Li(\mathfrak{q})$ is the
unique $\sigmaup$-invariant reductive complement of $\nr(\mathfrak{q})$
in $\mathfrak{q}$. Here $\sigmaup$ is the conjugation of $\kt$ with respect to the
real form $\kt_0$.
Stabilizing $\Li(\mathfrak{q})$, $\Ad_{\kt}(\mathbf{M}_0)$
also stabilizes its center and its semisimple ideal.
\end{proof}
\begin{prop}\label{prop:62} 
If $\mathfrak{q}\in\Par(\mathfrak{v},\mathbf{M}_0)$, then the inclusion
$\mathbf{M}_0\subset\mathbf{N}_{\mathfrak{q}}$ defines a
$\mathbf{K}_0$-equivariant smooth submersion
\begin{equation}
  \label{eq:6x3}
  \piup_{\mathfrak{q}}:M\to{N}_{\mathfrak{q}},
\end{equation}
which is a $CR$ map. If $M$ is $\mathfrak{n}$-reductive,
then $\piup_{\mathfrak{q}}$ has totally real fibers.  
\end{prop}
\begin{proof} We only need to prove the last statement. Assume that
$M$ is $\mathfrak{n}$-reductive. 
Every $X\in\mathfrak{v}$
uniquely decomposes as a sum $X=Y+Z$, 
with $Y\in\nr(\mathfrak{v})$,
$Z\in\mathfrak{m}$. Since $\mathfrak{m}\subset\Li(\mathfrak{q})$,  $X$ belongs to
$\mathfrak{v}\cap\bar{\mathfrak{q}}$ if and only if $Y\in\mathfrak{v}\cap
\Li(\mathfrak{q})$. But then $Y=0$, because 
$\nr(\mathfrak{v})\cap\Li(\mathfrak{q})=\{0\}$.
This shows that 
$\mathfrak{v}\cap\bar{\mathfrak{q}}=\mathfrak{m}$. Therefore 
the fibers of $\piup_{\mathfrak{q}}$ are totally real by the criterion   
\eqref{eq:n1}
of Proposition\,\ref{prop:214}.
\end{proof}
\begin{prop} \label{prop:62a}
\begin{enumerate} 
\item\label{e1} Let $\mathfrak{q}_1,\mathfrak{q}_2\in
\Par(\mathfrak{v})$, with $\mathfrak{q}_1\subset\mathfrak{q}_2$.
If 
$\mathfrak{q}_1\in\Par(\mathfrak{v},\mathbf{M}_0)$, then also
\mbox{$\mathfrak{q}_2\in\Par(\mathfrak{v},\mathbf{M}_0)$}.
\item\label{e2}
The parabolic regularization $\mathfrak{w}$ of $\mathfrak{v}$
belongs to $\Par(\mathfrak{v},\mathbf{M}_0)$. 
\item\label{e3}
In particular, if $\mathfrak{q}\in\Par(\mathfrak{v})$ and 
$\mathfrak{w}\subset\mathfrak{q}$ then 
$\mathfrak{q}\in\Par(\mathfrak{v},\mathbf{M}_0)$.
\end{enumerate}
\end{prop}
\begin{proof} \eqref{e1} follows from \eqref{eq:6xx}.\par
\eqref{e2}.\;
Consider the sequence of Lie subalgebras defined by \eqref{eq:41b}. 
For each $i=1,2,\hdots$, set 
\begin{equation*}
  \mathbf{N}_{\mathfrak{v}_i}=\{g\in\mathbf{K}_0\mid\mathrm{Ad}(g)(\mathfrak{v}_i)=
\mathfrak{v}_i\},\quad N_i=\mathbf{K}_0/\mathbf{N}_{\mathfrak{v}_i}.
\end{equation*}
Since $\mathfrak{v}$ is $\mathrm{Ad}(\mathbf{M}_0)$-invariant,
also its characteristic ideal $\nr(\mathfrak{v})$ 
and its normalizer $\mathfrak{v}_1$ are $\mathrm{Ad}(\mathbf{M}_0)$-invariant.
Therefore $\mathbf{M}_0\subset\mathbf{N}_{\mathfrak{v}_1}$. 
\par 
Let $i\geq{1}$. If $g\in\mathbf{N}_{\mathfrak{v}_i}$ and $X\in\mathfrak{v}_{i+1}=
\N_{\kt}(\nr(\mathfrak{v}_i))$, we have
\begin{align*}
  [\Ad_{\kt}(g)(X),Y]=\Ad_{\kt}(g)([X,\Ad_{\kt}(g^{-1})(Y)])\in
\Ad_{\kt}(g)(\nr(\mathfrak{v}_i))=
 \nr(\mathfrak{v}_i),\;\;\\
\forall Y\in\nr(\mathfrak{v}_i),
\end{align*}
showing that $\Ad_{\kt}(g)(X)\in\mathfrak{v}_{i+1}$. Therefore we have
$\mathbf{N}_{\mathfrak{v}_i}\subset\mathbf{N}_{\mathfrak{v}_{i+1}}$. 
Hence, with $\mathfrak{w}=\mathfrak{v}_m$, we obtain
$\mathbf{M}_0\subset\mathbf{N}_{\mathfrak{v}_1}\subset\mathbf{N}_{\mathfrak{v}_2}
\subset\cdots\subset\mathbf{N}_{\mathfrak{v}_m}=\mathbf{N}_{\mathfrak{w}}$,
proving that $\mathfrak{w}\in\Par(\mathfrak{v},\mathbf{M}_0)$. 
\end{proof}
\begin{prop} \label{prop:63} 
Let $\mathfrak{v}$ be a complex Lie subalgebra of $\kt$.
If $\mathfrak{q}$ is
any maximal (with respect to inclusion) element of 
$ {\Par}(\mathfrak{v})$, then 
\begin{equation}
  \label{eq:6x2}
  \mathfrak{q}=
\lie(\nr(\mathfrak{v})+\Li(\mathfrak{q})).
\end{equation}
\end{prop}
\begin{proof}
 Let $\mathfrak{q}$ be
any maximal element of $ {\Par}(\mathfrak{v})$ and fix a maximal torus
$\mathfrak{t}_0$ in 
$\mathfrak{L}(\mathfrak{q})$. With $\zt=\zt(\mathfrak{q})$, 
following \cite{Kos2010}, we define the $\zt$-roots as the 
non zero elements
$\nuup$ of the dual $\mathfrak{z}^*_{\mathbb{R}}$ of $\zt_{\mathbb{R}}=
i\zto(\mathfrak{q})$
for which
\begin{equation*}
  \kt_{\nuup}=\{X\in\kt\mid [Z,X]=\nuup(Z)X,\;\forall 
Z\in\mathfrak{z}_{\mathbb{R}}\}\neq\{0\}.
\end{equation*}
The set $\mathcal{Z}$ of the $\mathfrak{z}$-roots is a finite subset of
$\mathfrak{z}^*_{\mathbb{R}}$, with the properties that
\begin{equation*}
  \begin{cases}
    \nuup\in\mathcal{Z}\Longrightarrow -\nuup\in\mathcal{Z};\;\;\\
\kt_{\nuup_1+\nuup_2}=[\kt_{\nuup_1},\kt_{\nuup_2}];\\
    \text{for all 
$\nuup\in\mathcal{Z}$,
$\kt_{\nuup}$ is an irreducible $\mathfrak{L}(\mathfrak{q})$-module};\\
    \text{there exists an ordering in $\mathfrak{z}^*_{\mathbb{R}}$ such that\;\;
    $\nr(\mathfrak{q})={\sum}_{\nuup>0}\kt_{\nuup}$;}\\[0.2cm]
    \begin{aligned}
      &\text{there is a system of linearly independent positive simple roots 
$\nuup_1,\hdots,\nuup_{\ell}$ such}\\[-0.2cm]
      &\text{that every positive $\nuup\in\mathcal{Z}$ uniquely decomposes
      as $\nuup={\sum}_{i=1}^{\ell}c_i\nuup_i$, with $c_i\in\mathbb{Z}_+$.}
    \end{aligned}
  \end{cases}
\end{equation*}
\par
We claim that $\kt_{\nuup_i}\cap(\nr(\mathfrak{v})+\mathfrak{L}(\mathfrak q))
\neq\{0\}$ for 
$i=1,\hdots,\ell$.
Assume by contradiction that there is $j$, with $1\leq{j}\leq{\ell}$,
such that $\kt_{\nuup_j}\cap(\nr(\mathfrak{v})+\Li(\mathfrak{q}))
=\{0\}$. Then we consider
the parabolic $\mathfrak{q}'=\mathfrak{q}\oplus\kt_{-\nu_j}$.
Since the elements $X$ of $\nr(\mathfrak{v})$ uniquely decompose
as sums $X=X_0+{\sum}_{\nuup>0}X_{\nuup}$, with 
$X_0\in\mathfrak{L}(\mathfrak{q})$ and $X_{\nuup}\in\kt_{\nuup}$,
we have 
\begin{equation*}
\mathfrak{L}(\mathfrak{q}')\cap\nr(\mathfrak{v})=(\mathfrak{L}(\mathfrak{q})\oplus
\kt_{\nu_j}\oplus\kt_{-\nu_j})\cap\nr(\mathfrak{v})
  =(\mathfrak{L}(\mathfrak{q})\oplus\kt_{\nu_j})\cap\nr(\mathfrak{v})=\{0\}.
\end{equation*}
Indeed, if $Y=X_{0}+X_{\nuup_j}=Y$, with $Y\in\nr(\mathfrak{v})$,
$X_{0}\in\mathfrak{L}(\mathfrak{q})$,
$X_{\nuup_j}\in\kt_{\nuup_j}$, 
then 
$X_{\nuup_j}=0$, because 
$X_{\nuup_j}\in (\mathfrak{L}(\mathfrak{q})+\nr(\mathfrak{v}))\cap\kt_{\nuup_j}=\{0\}$.
Hence $Y=0$, because 
$X_0=Y\in\mathfrak{L}(\mathfrak{q})\cap\nr(\mathfrak{v})=\{0\}$.\par
This would give 
$\mathfrak{q}\subsetneqq\mathfrak{q}'\in{\Par}(\mathfrak{v})$,
contradicting the maximality of $\mathfrak{q}$. \par
Therefore 
\begin{align*}
 \kt_{\nuup_i}\subset \nr(\mathfrak{v})+\Li(\mathfrak{q})+[ \Li(\mathfrak{q}),
\nr(\mathfrak{v})]+[ \Li(\mathfrak{q}),[\Li(\mathfrak{q}),
\nr(\mathfrak{v})]]+\cdots\quad \\
\subset \; \lie(\nr(\mathfrak{v})+
\Li(\mathfrak{q})),
\end{align*}
and this implies \eqref{eq:6x2}.
\end{proof}
We obtain, from Proposition\,\ref{prop:62} and Proposition\,\ref{prop:63},
\begin{thm}[the deployment theorem] \label{thm68}
If $M$ is $\mathfrak{n}$-reductive, then for all maximal 
$\mathfrak{q}$ in $\Par(\mathfrak{v},\mathbf{M}_0)$, 
\eqref{eq:6x3} 
is a $CR$-deployment. \qed
\end{thm}
The connected compact Lie group $\mathbf{K}_0$ admits a linearization and 
hence a complexification $\mathbf{K}$, which is unique modulo isomorphisms.
Let $V=\mathbf{K}/\mathbf{V}$ be the $\mathbf{K}$-realization of
$M$ (see Theorem\,\ref{thm45}) and $N_\mathfrak{q}
=\mathbf{K}/\mathbf{Q}$ the complex
flag manifold associated to the parabolic $\mathfrak{q}$ of $\kt$. 
We have: 
\begin{thm}[factoring through the realization]\label{thm69}
Let 
$\mathfrak{q}\in\Par(\mathfrak{v},\mathbf{M}_0)$ and 
$\mathbf{Q}$ be the analytic Lie subgroup of $\mathbf{K}$ with
$\Lie(\mathbf{Q})=\mathfrak{q}$. Then $\mathbf{V}\subset
\mathbf{Q}$, and we obtain a commutative diagram:
\begin{equation*}
  \begin{CD}
    M=\mathbf{K}_0/\mathbf{M}_0 @>>> V=\mathbf{K}/\mathbf{V}\,\;\\
@V{\piup_{\mathfrak{q}}}VV @VV{\tilde{\piup}_{\mathfrak{q}}}V \\
N_{\mathfrak{q}}=\mathbf{K}_0/\mathbf{N}_{\mathfrak{q}}@>>>N_{\mathfrak{q}}
=\mathbf{K}/\mathbf{Q},\;
  \end{CD}
\end{equation*}
where
the top horizontal arrow is a generic $CR$-embedding and
a $\mathbf{K}$-realization, 
the right vertical arrow
is a holomorphic submersion and the bottom horizontal arrow is an extension
of the holomorphic action of $\mathbf{K}_0$ on $N_{\mathfrak{q}}$ to a
holomorphic action of its complexification $\mathbf{K}$.
Hence $\piup_{\mathfrak{q}}$ is the restriction to $M$ of a holomorphic submersion
of the ambient space $V$.
\par
If, moreover, $\nr(\mathfrak{v})\subset\nr(\mathfrak{q})$,
then the fibers of $\tilde{\piup}_{\mathfrak{q}}$ are Stein.
\end{thm}
\begin{proof} 
We know from Theorem\,\ref{thm45} that the complex Lie subgroup
$\mathbf{V}$ has a Levi-Chevalley decomposition 
$\mathbf{V}=\mathbf{V}_n\mathbf{M}$, where 
$\mathbf{V}_n$ is the analytic Lie subgroup of $\mathbf{K}$ with
Lie algebra $\nr(\mathfrak{v})$ and 
$\mathbf{M}=\{g\exp(iX)\mid g\in\mathbf{M}_0,\; X\in\mathfrak{m}_0\}$.
If $\mathfrak{q}\in\Par(\mathfrak{v},\mathbf{M}_0)$, then
$\mathbf{M}_0\subset\mathbf{Q}$ and hence also $\mathbf{V}\subset
\mathbf{Q}$ because $\mathbf{M}_0$ is a maximal compact subgroup of
$\mathbf{V}$, and hence 
all connected components of $\mathbf{V}$ intersect
$\mathbf{M}_0$. This gives the first part of the statement.\par
Assume now that
$\mathfrak{q}\in\Par(\mathfrak{v},\mathbf{M}_0)$ and
$\nr(\mathfrak{v})\subset\nr(\mathfrak{q})$.
 Then $\mathbf{Q}$ has a 
Levi-Chevalley
decomposition $\mathbf{Q}=\mathbf{Q}_n\mathbf{Q}_r$, with
analytic complex Lie subgroups $\mathbf{Q}_n,\,\mathbf{Q}_r$
with $\Lie(\mathbf{Q}_n)=\nr(\mathfrak{q})$, 
$\Lie(\mathbf{Q}_r)=\Li(\mathfrak{q})$.
\par
Since $\mathbf{M}\subset\mathbf{Q}_r$ and $\mathbf{Q}_n$ is a normal Lie
subgroup of $\mathbf{Q}$, the product
\begin{equation*}
\mathbf{Q}_n\mathbf{M}=\{g_1g_2\mid g_1\in\mathbf{Q}_n,\; g_2\in\mathbf{M}\}  
\end{equation*}
is a closed Lie subgroup of $\mathbf{K}$. Since 
we assumed that
$\nr(\mathfrak{v})\subset
\nr(\mathfrak{q})$, we have $\mathbf{V}_n\subset\mathbf{Q}_n$ and thus
the inclusions 
$\mathbf{V}\subset \mathbf{Q}_n\mathbf{M}\subset \mathbf{Q}$
yield
a commutative diagram
\begin{equation*}
  \xymatrix{ V=\mathbf{K}/\mathbf{V} \ar[rr]^{\tilde{\piup}_1}
\ar[drr]_{\tilde{\piup}_{\mathfrak{q}}}
&&\mathbf{K}/(\mathbf{Q}_n\mathbf{M}) \ar[d]^{\tilde{\piup}_2}\\
&&N_{\mathfrak{q}}=\mathbf{K}/\mathbf{Q}}\;\;
\end{equation*}
of $\mathbf{K}$-equivariant maps. The fiber of $\tilde{\piup}_1$ 
is $(\mathbf{Q}_n\mathbf{M})/\mathbf{V}\simeq
\mathbf{Q}_n/(\mathbf{Q}_n\cap\mathbf{V})$,
and therefore is complex Euclidean. The fiber of $\piup_2$ is 
$\mathbf{Q}/(\mathbf{Q}_n\mathbf{M})\simeq\mathbf{Q}_r/\mathbf{M}$
and hence Stein because $\mathbf{Q}_r$ is Stein and 
$\mathbf{M}$  is the complexification of the compact subgroup
$\mathbf{M}_0$ (see \cite{Mts60}). Then $\mathbf{Q}/\mathbf{V}$
is the total space of a fibration $\mathbf{Q}/\mathbf{V}\to
\mathbf{Q}/(\mathbf{Q}_n\mathbf{M})$ with a Stein base
and complex Euclidean fibers and therefore is Stein by \cite{MM60}.
\end{proof}

\begin{exam} 
In the proof of Proposition\,\ref{prop:63} we show that the 
root spaces corresponding to the simple positive $\mathfrak{z}$-roots
belong to the 
$\Li(\mathfrak{q})$-module generated by 
$\nr(\mathfrak{v})\! +\!\Li(\mathfrak{q})$. 
A natural question arises whether, in the definition of $CR$-spread,
one could substitute to $\lie(\nr(\mathfrak{v})\! +\!\Li(\mathfrak{q}))$ the
$\Li(\mathfrak{q})$-module generated by 
$\nr(\mathfrak{v})\! +\!\Li(\mathfrak{q})$. We provide an example showing that 
the statement of Theorem\,\ref{thm68} 
does not hold with this stronger notion of $CR$-deployment.
\par
Let $\kt_0=\mathfrak{so}(7)$. Fix a maximal torus 
$\mathfrak{t}_0$ in  
$\kt_0$, let $\mathcal{K}=\{\pm{e}_i,\;\pm{e}_j\!\pm\!{e}_h\mid
1\leq{i}\leq{3},\;1\leq{j}<h\leq{3}\}$ the root space of 
$(\kt,\mathfrak{t})$. Take $\mathfrak{v}=\kt^{e_1-e_3}\oplus\kt^{e_2}$.
This is an Abelian Lie algebra of nilpotent elements.
The parabolic regularization $\mathfrak{w}$
of $\mathfrak{v}$ is the Borel subalgebra
associated to the standard ordering of $\mathcal{K}$, corresponding
to the cross-marked Satake diagram:
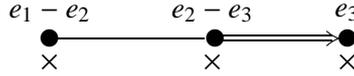
\begin{equation*}
  \xymatrix@R=-.3pc{ e_1-e_2 & e_2-e_3& e_3\\
\!\!\medbullet\!\!\ar@{-}[r]&\!\!\medbullet\!\!\!
\ar@{=>}[r]&\!\!\medbullet\!\! \\
\times& \times&\times}
\end{equation*}
There are exactly two maximal parabolic subalgebra in $\Par(\mathfrak{v})$
which contain $\mathfrak{w}$. \par
The first maximal
$\mathfrak{q}\in\Par(\mathfrak{v})$, with $\mathfrak{w}\subset\mathfrak{q}$, 
is associated to the cross-marked
Satake diagram
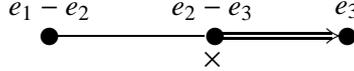
\begin{equation*}
  \xymatrix@R=-.3pc{ e_1-e_2 & e_2-e_3& e_3\\
\!\!\medbullet\!\!\ar@{-}[r]&\!\!\medbullet\!\!\!
\ar@{=>}[r]&\!\!\medbullet\!\! \\
& \times}
\end{equation*}
We have $\mathcal{Q}_n=\{e_1,e_2,e_1\pm{e}_3,e_2\pm{e_3},e_1+e_2\}$,
$\mathcal{Q}_r=\{\pm{e}_3,\pm(e_1-e_2)\}$ and, accordingly, the
centralizer $\mathfrak{z}$ of $\Li(\mathfrak{q})$ in $\mathfrak{t}$ 
consists of the elements $H\in\mathfrak{t}$ with
$e_1(H)=e_2(H)$ and $e_3(H)=0$. The $\mathfrak{z}$-roots are 
$\{\pm\nu,\pm 2\nu\}$, with 
\begin{align*}
  \kt_{\nu}\,&=\kt^{e_1}\oplus\kt^{e_2}\oplus\kt^{e_1-e_3}\oplus\kt^{e_1+e_3}\oplus
\kt^{e_2-e_3}\oplus\kt^{e_2+e_3},\\
\kt_{2\nu}&=\kt^{e_1+e_2}.
\end{align*}
We note that $\mathfrak{v}\subset\kt_{\nu}$. Hence $\mathfrak{q}$
is generated by $\mathfrak{v}+\Li(\mathfrak{q})$ as a Lie algebra,
but not as an $\Li(\mathfrak{q})$-Lie module. \par\smallskip
The other possible choice of a maximal
$\mathfrak{q}\in\Par(\mathfrak{v})$ with $\mathfrak{w}\subset\mathfrak{q}$ 
corresponds to the cross-marked
Satake diagram 
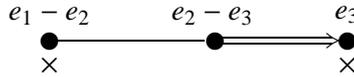
\begin{equation*}
  \xymatrix@R=-.3pc{ e_1-e_2 & e_2-e_3& e_3\\
\!\!\medbullet\!\!\ar@{-}[r]&\!\!\medbullet\!\!\!
\ar@{=>}[r]&\!\!\medbullet\!\! \\
\times && \times}
\end{equation*}
In this case, $\mathcal{Q}_n=\{e_1,e_2,e_3,e_1\pm{e}_2,e_1\pm{e}_3,e_2+e_3\}$,
$\mathcal{Q}_r=\{\pm(e_2-e_3)\}$. The corresponding set of $\mathfrak{z}$-roots
yields
\begin{align*}
  \kt_{\nuup_1}&=\{e_1-e_2,e_1-e_3\},\\
\kt_{\nuup_2}&=\{e_2,e_3\},\\
\kt_{\nuup_1+\nuup_2}&=\{e_1\},\\
\kt_{\nuup_1+2\nuup_2}&=\{e_1+e_2,e_1+e_3\},\\
\kt_{2\nuup_2}&=\{e_2+e_3\}.
\end{align*}
We note that $\mathfrak{v}\cap\kt_{\nuup_1}=\kt^{e_1-e_2}\neq\{0\}$,  
$\mathfrak{v}\cap\kt_{\nuup_2}=\kt^{e_3}\neq\{0\}$, and
$\mathfrak{v}=\mathfrak{v}\cap\kt_{\nuup_1}\oplus\mathfrak{v}\cap\kt_{\nuup_2}$.
Also in this case $\mathfrak{v}+\Li(\mathfrak{q})$
generates $\mathfrak{q}$ as a Lie subalgebra, but not as an 
$\Li(\mathfrak{q})$-Lie module.\par
This example also shows that, in general, deployments corresponding to 
different maximal elements
of $\Par(\mathfrak{v},\mathbf{M}_0)$ may  be non equivalent.
\end{exam}
\subsection{Characterization of $\Par(\mathfrak{v},\mathbf{M}_0)$ }
\label{sec:81}
For every $A\in\kt_0$, $\ad_{\kt}(A)$ is semisimple and  
has purely
imaginary eigenvalues. We associate to
$A\in\kt_0$ the parabolic subalgebra
\begin{align}
  \label{eq:8xa}
   \mathfrak{q}_A={\sum}_{\lambda\geq{0}}\{X\in\kt\mid [A,X]=i\lambda{X}\}.
\end{align}
As we already noted in \S\ref{su41},
the corresponding complex manifold $N_{\mathfrak{q}_A}$ 
is the totally complex $\mathbf{K}_0$-homogeneous
$CR$ manifold with isotropy
\begin{equation}
  \label{eq:8x1}
  \mathbf{N}_{\mathfrak{q}_A}=\{g\in\mathbf{K}_0\mid \Ad_{\kt}(g)(A)=A\}
\end{equation}
and $CR$ algebra $(\kt_0,\mathfrak{q}_A)$ at the base point $p_0$. \begin{thm} With 
\begin{equation}
  \label{eq:95}
  \mathfrak{a}_0(\mathbf{M}_0)=\{A\in\kt_0\mid \Ad_{\kt}(g)(A)=A,
\;\forall g\in\mathbf{M}_0\},
\end{equation}
we have:
\begin{equation}\label{eq:96}
 \Par(\mathfrak{v},\mathbf{M}_0)
=\{\mathfrak{q}_A\mid A\in\mathfrak{a}_0(\mathbf{M}_0),\;\; \mathfrak{v}\subset
\mathfrak{q}_A\}.
\end{equation}
 \end{thm} 
\begin{proof} Let us denote by $\mathcal{P}$ the right hand side of \eqref{eq:96}.
By \eqref{eq:8x1}, we have the inclusion
$\Par(\mathfrak{v},\mathbf{M}_0)\subset\mathcal{P}$. The opposite inclusion
is obvious, and hence we have equality.
\end{proof}
\begin{exam}\label{ex89}
Let $\mathbf{K}_0$ be a connected compact Lie group, 
$\mathfrak{t}_0$ a maximal torus
in $\mathbf{K}_0$ and $\mathbf{M}_0$ the normalizer of $\mathfrak{t}_0$  
in $\mathbf{K}_0$. If $\mathbf{T}_0$ is the analytic subgroup with
$\Lie(\mathbf{T}_0)=\mathfrak{t}_0$, then 
$\mathbf{T}_0$ is a normal subgroup of $\mathbf{M}_0$ and the quotient
$\mathbf{M}_0/\mathbf{T}_0$ is the Weyl group.
Consider on $M=\mathbf{K}_0/\mathbf{M}_0$ the totally real
$CR$ structure defined at the base point by the $CR$ algebra
$(\kt_0,\mathfrak{t})$, with $\mathfrak{t}=\mathbb{C}\otimes_{\mathbb{R}}
\mathfrak{t}_0$. Then $\Par(\mathfrak{t},\mathbf{M}_0)$ only contains
$\kt$. Thus there are no $CR$ structures of type $M_{\mathfrak{q}}$ with
positive $CR$ dimension. This shows that, in general, 
the set $\{A\in\kt_0\mid [A,\mathfrak{m}_0]=\{0\},\;
\mathfrak{v}\!\subset\!\mathfrak{q}_A\}$,
when $\mathbf{M}_0$ is not connected, 
 can be larger
than $\mathfrak{a}_0(\mathbf{M}_0)$. 
\end{exam}
\subsection{The deployment theorem for 
homogeneous foliated complex manifolds} 
With a proof similar to that of Theorem\,\ref{thm69} we obtain
\begin{thm} \label{th713}
Let $\mathbf{V}$ be a decomposable subgroup of a
complex reductive Lie group
$\mathbf{K}$, and $V=\mathbf{K}/\mathbf{V}$. 
Fix a Levi-Chevalley decomposition \eqref{eq:6x16} of $\mathbf{V}$ 
and let $(M^{\mathbb{C}},{\Dt}_{\mathfrak{v}})$, 
with ${\mathfrak{v}}=\Lie(\mathbf{V})$,
be the corresponding Stein lift of $V$ of \S\ref{sub47}, \S\ref{sub48}.
\par
Let $\mathfrak{q}$ be maximal in $\Par(\mathfrak{v})$ and
$\mathbf{Q}$ the corresponding parabolic subgroup of $\mathbf{K}$.
Then $\mathbf{V}\subset\mathbf{Q}$ and $\mathbf{Q}$
has a Levi-Chevalley decomposition 
 $\mathbf{Q}=\mathbf{Q}_n\cdot\mathbf{Q}_r$  
with $\mathbf{M}\subset\mathbf{Q}_r$,
$\mathbf{V}_n\cap\mathbf{Q}_r=\{1\}$. 
\par
Let $N_{\mathfrak{q}}=\mathbf{K}/\mathbf{Q}$ and 
$(N^{\mathbb{C}},{\Dt}_{\mathfrak{q}})$,
with $N^{\mathbb{C}}_{\mathfrak{q}}=\mathbf{K}/\mathbf{Q}_r$, its Stein lift. 
Then we have a commutative diagram
\begin{equation*}
  \begin{CD}
    M^{\mathbb{C}} @>{\tilde{\phiup}}>> N_{\mathfrak{q}}^{\mathbb{C}}\\
@V{\piup}VV       @VV{\piup'}V \\
V @>>{\phiup}> N_{\mathfrak{q}}
  \end{CD}
\end{equation*}
of $\mathbf{K}$-equivariant fibrations, and moreover
$\phiup,{\phiup}_r$ have Stein fibers, and $\tilde{\phiup}$ is a deployment
of foliated complex manifolds. \qed
\end{thm}
\section{Lifted $CR$ structures}\label{lif}
We keep the notation of \S\ref{sec:dp}. Our next goal
is to describe all \textsl{horocyclic} 
$\mathfrak{n}$-reductive $\mathbf{K}_0$-homogeneous $CR$ structures
of a given $\mathbf{K}_0$-space $M$.
In fact, we take up the slightly more general question of strengthening
a given
 $\mathfrak{n}$-reductive $\mathbf{K}_0$-homogeneous $CR$ structure on $M$. 
The totally real $M$ corresponds indeed to the choice 
of the trivial $\mathfrak{n}$-reductive $\mathbf{K}_0$-equivariant 
$CR$ structure defined at $p_0$ by~$(\kt_0,\mathfrak{m})$.
Throughout this section we will always assume 
that $M$ is $\mathfrak{n}$-reductive,
i.e. that
\begin{equation}\label{eq:7a}
  \mathfrak{v}=\mathfrak{m}\oplus\nr(\mathfrak{v}),\quad\text{with}\quad
\mathfrak{m}=\Li(\mathfrak{v})=\mathfrak{v}\cap\bar{\mathfrak{v}}. 
\end{equation}
We begin with a simple algebraic lemma.
\begin{lem}\label{lm71a}
 If $\mathfrak{q}\in\Par(\mathfrak{v})$, then 
 $\mathfrak{v}_{\mathfrak{q}}\!=\!\mathfrak{v}\! +\!\nr(\mathfrak{q})$ 
is a complex Lie algebra 
with $\nr(\mathfrak{v}_{\mathfrak{q}})=\nr(\mathfrak{v})+\nr(\mathfrak{q})$,
and admits a Levi-Chevalley decomposition 
with reductive Levi factor $\mathfrak{m}$. In particular,
$\Li(\mathfrak{v}_{\mathfrak{q}})=
\mathfrak{v}_{\mathfrak{q}}\cap\bar{\mathfrak{v}}_{\mathfrak{q}}=\mathfrak{m}$.
\end{lem} 
\begin{proof}
 If $\mathfrak{q}\in\Par(\mathfrak{v})$, then
$\mathfrak{v}\subset\mathfrak{q}$ and hence $[\mathfrak{v},\nr(\mathfrak{q})]
\subset [\mathfrak{q},\nr(\mathfrak{q})]\subset\nr(\mathfrak{q})$
shows that $\mathfrak{v}_{\mathfrak{q}}=
\mathfrak{v}+\nr(\mathfrak{q})$ is a Lie subalgebra of $\kt$,
containing $\nr(\mathfrak{q})$ as an ideal of nilpotent elements.
Moreover, $\nr(\mathfrak{v})+\nr(\mathfrak{q})$ is also an ideal in
$\mathfrak{v}_{\mathfrak{q}}$ which consists of nilpotent elements and has
the reductive complement $\mathfrak{m}$. 
\end{proof}
\begin{prop} \label{pp72} Let $\mathfrak{q}\in\Par(\mathfrak{v},\mathbf{M}_0)$,
and set
$\mathfrak{v}_{\mathfrak{q}}\! =\!\mathfrak{v}\! +\!\nr(\mathfrak{q})$.
Then the $CR$ algebra 
 $(\kt_0,\mathfrak{v}_{\mathfrak{q}})$ at the base point defines on 
$M=\mathbf{K}_0/\mathbf{M}_0$
 another $\mathfrak{n}$-reductive
$\mathbf{K}_0$-equivariant $CR$ structure, which strengthens
the one defined by $(\kt_0,\mathfrak{v})$.
\end{prop} 
\begin{proof}
 It suffices to notice that, since $\mathfrak{q}\in
\Par(\mathfrak{v},\mathbf{M}_0)$, 
the subalgebras
 $\mathfrak{v}$ and
 $\nr(\mathfrak{q})$ are both $\Ad_{\kt}(\mathbf{M}_0)$-invariant. 
Hence 
 $\mathfrak{v}_{\mathfrak{q}}$ is $\Ad_{\kt}(\mathbf{M}_0)$-invariant.
Using Lemma\,\ref{lm71a}, the statement follows from 
Proposition\,\ref{prop28}.
\end{proof}

\begin{rmk}
The set $\Par(\mathfrak{m},\mathbf{M}_0)$ 
parameterizes   the horocyclic $\mathfrak{n}$-reductive
$\mathbf{K}_0$-homogeneous $CR$ structures
on $M=\mathbf{K}_0/\mathbf{M}_0$. These are the ones which are 
defined at the base point by a
$CR$ algebra of the form $(\kt_0,\mathfrak{m}\oplus\nr(\mathfrak{q}))$,
with $\mathfrak{q}\in\Par(\mathfrak{m})$.
\end{rmk}
\par
\begin{ntz} We denote by $M_{\mathfrak{q}}$ the
$\mathbf{K}_0$-homogeneous $CR$ manifold with isotropy $\mathbf{M}_0$ and
$CR$ algebra $(\kt_0,\mathfrak{v}_{\mathfrak{q}})$ at the base point,
described in Proposition\,\ref{pp72}. 
\end{ntz}
We get from Proposition\,\ref{prop:62},
Lemma\,\ref{lm71a} and Proposition\,\ref{pp72}:
\begin{thm}[The lift theorem]
  \label{prop74}
For $\mathfrak{q}\in\Par(\mathfrak{v},\mathbf{M}_0)$, the 
$\mathbf{K}_0$-equivariant map 
\begin{equation}\label{lft}
  \phiup_{\mathfrak{q}}:M_{\mathfrak{q}}\to{N}_{\mathfrak{q}}
\end{equation}
is a $CR$ lift.
\qed
\end{thm}
\begin{prop}[realization] \label{pp79}
Let $\mathbf{K}$ be a complexification of $\mathbf{K}_0$. 
 \begin{enumerate}
\item
For every $\mathfrak{q}\in\Par(\mathfrak{v},\mathbf{M}_0)$, 
$M_{\mathfrak{q}}$ is $\mathbf{K}$-realizable.
\item 
Let $\mathbf{V}_{\mathfrak{q}}$ be a complex Lie subgroup of
$\mathbf{K}$ with  $\Lie(\mathbf{V}_{\mathfrak{q}})=\mathfrak{v}_{\mathfrak{q}}$,
$\mathbf{V}_{\mathfrak{q}}\cap\mathbf{K}_0=\mathbf{M}_0$, and
$\mathbf{Q}=\{g\in\mathbf{K}\mid \Ad_{\kt}(g)(\mathfrak{q})=\mathfrak{q}\}$
the parabolic subgroup with $\Lie(\mathbf{Q})=\mathfrak{q}$. 
Then $\mathbf{V}_{\mathfrak{q}}\subset\mathbf{Q}$ 
and we obtain a commutative diagram 
\begin{equation*}
  \begin{CD}
    M_{\mathfrak{q}}=\mathbf{K}_0/\mathbf{M}_0 
@>>> \mathbf{K}/\mathbf{V}_{\mathfrak{q}}\,\\
@V{\piup_{\mathfrak{q}}}VV @VV{\tilde{\piup}_{\mathfrak{q}}}V \\
N_{\mathfrak{q}}=\mathbf{K}_0/\mathbf{N}_{\mathfrak{q}}@>>>\mathbf{K}/\mathbf{Q},
  \end{CD}
\end{equation*}
where
the top horizontal arrow is the $\mathbf{K}$-realization of
$M_{\mathfrak{q}}$, the bottom horizontal 
arrow is an extension of holomorphic action from $\mathbf{K}_0$ 
to $\mathbf{K}$, 
and the right vertical arrow $\tilde{\piup}_{\mathfrak{q}}$
is a $\mathbf{K}$-equivariant holomorphic fibration with
Stein fibers.
\end{enumerate}
\end{prop}
\begin{proof} It suffices to apply Theorem\,\ref{thm69} 
to $M_{\mathfrak{q}}$, using \eqref{12} of Proposition\,\ref{prop:84x}.
\end{proof}
\subsection{Factorization of the lift}
\label{sec:81x} Let $\tauup_0$ be a maximal torus of $\mathfrak{m}_0$.
An element $\mathfrak{q}$ of $\Par(\mathfrak{v},\mathbf{M}_0)$ contains a
maximal torus $\mathfrak{t}_0$ 
of $\kt_0$ with $\mathfrak{t}_0\cap\mathfrak{m}_0=\tauup_0$.
Set
\begin{equation}\label{eq:73d}
\mathfrak{z}_0\! =\!
\{A\in\mathfrak{t}_0\mid\Ad_{\kt}(g)(A)=A,\,\forall{g}\in\mathbf{M}_0\}
,\;\Z_{\mathfrak{t}_0}(\mathfrak{m}_0)\! =\!  \{A\in\mathfrak{t}_0\mid
[A,\mathfrak{m}_0]=\{\! 0\!\}\!\}.
\end{equation}
We have $\mathfrak{z}_0\subset\Z_{\mathfrak{t}_0}(\mathfrak{m}_0)$,
with $\mathfrak{z}_0=\Z_{\mathfrak{t}_0}(\mathfrak{m}_0)$
when $\mathbf{M}_0$ is connected, and $\mathfrak{q}=\mathfrak{q}_A$ for some
$A\in\mathfrak{z}_0$.
The meaning of $\Z_{\mathfrak{t}_0}(\mathfrak{m}_0)$ is illustrated by the following
\begin{lem}
The necessary and sufficient condition for
$\mathfrak{m}_0$ to be $\mathfrak{t}_0$-regular is that 
\begin{equation}\label{eq:73e}
 \mathfrak{t}_0=\tauup_0
+\Z_{\mathfrak{t}_0}(\mathfrak{m}_0).
\end{equation}\end{lem} 
\begin{proof}
 The condition is clearly sufficient. 
 Vice versa, when $\mathfrak{m}_0$ is $\mathfrak{t}_0$-regular,
 $\mathfrak{z}_0$ is the center of the reductive subalgebra
 $\mathfrak{t}_0+\mathfrak{m}_0$, and $\tauup_0$ contains the maximal torus
 $\mathfrak{t}_0\cap[\mathfrak{m}_0,\mathfrak{m}_0]$ of its semisimple ideal
 $[\mathfrak{m}_0,\mathfrak{m}_0]$, so that \eqref{eq:73e} holds true.
\end{proof}
Fix $\mathfrak{q}\in\Par(\mathfrak{v},\mathbf{M}_0)$ and set 
\begin{equation}
 \mathfrak{s}_{\mathfrak{q}}=[\Lio(\mathfrak{q}),\Lio(\mathfrak{q})],\quad
 \mathbf{S}_{\mathfrak{q}}=\text{analytic Lie subgroup with $\Lie(\mathbf{S}_{\mathfrak{q}})
 =\mathfrak{s}_{\mathfrak{q}}$}.
\end{equation}
We note that the semisimple real Lie group
$\mathbf{S}_{\mathfrak{q}}$ is compact and algebraic.
\begin{lem}\label{lem712}
  We have 
$\Ad_{\kt_0}(g)(\mathfrak{s}_{\mathfrak{q}})
=\mathfrak{s}_{\mathfrak{q}}$ and $\ad(g)(\mathbf{S}_{\mathfrak{q}})=
\mathbf{S}_{\mathfrak{q}}$
for all \mbox{$g\in\mathbf{M}_0$.}
\end{lem}
\begin{proof} From $\mathbf{M}_0\subset\mathbf{N}_{\mathfrak{q}}$ we have
$\Ad_{\kt_0}(g)(\Lio(\mathfrak{q}))=\Lio(\mathfrak{q})$, because
$\Lie(\mathbf{N}_{\mathfrak{q}})=\Lio(\mathfrak{q})$.
 Then also the semisimple ideal $\mathfrak{s}_{\mathfrak{q}}$ of
 $\Lio(\mathfrak{q})$ is
 $\Ad_{\kt_0}(\mathbf{M}_0)$-invariant. Hence the analytic Lie subgroup 
 $\mathbf{S}_{\mathfrak{q}}$ generated by $\mathfrak{s}_{\mathfrak{q}}$
 is $\ad(\mathbf{M}_0)$-invariant. 
\end{proof}
\begin{cor}
  The product 
  \begin{equation}\label{77}
    \mathbf{L}_{\mathfrak{q}}=
    \mathbf{M}_0\mathbf{S}_{\mathfrak{q}}=\{g_1g_2\mid g_1\in\mathbf{M}_0,\;
g_2\in\mathbf{S}_{\mathfrak{q}}\}
  \end{equation}
is a compact subgroup of $\mathbf{K}_0$, contained in $\mathbf{N}_{\mathfrak{q}_A}$.
\end{cor}
\begin{proof}
  The product \eqref{77} is a group by Lemma\,\ref{lem712} and is contained
in $\mathbf{N}_{\mathfrak{q}_A}$ because both $\mathbf{M}_0$ and
$\mathbf{S}_{\mathfrak{q}}$ are contained in $\mathbf{N}_{\mathfrak{q}_A}$.
\end{proof}
\begin{lem}\label{lem714} The quotient
$\mathbf{L}_{\mathfrak{q}}/\mathbf{M}_0$ is a homogeneous space of 
homotopic characteristic
\begin{equation}\label{eq:810}
 c_{\mathfrak{q}}=
 \rank \mathfrak{s}_{\mathfrak{q}}-\dim (\tauup_0\cap\mathfrak{s}_{\mathfrak{q}}).
\end{equation}
In particular, $\chiup(\mathbf{L}_{\mathfrak{q}}/\mathbf{M}_0)>0$ if and only if 
$\mathfrak{t}_0\cap\mathfrak{s}^{A}_0\subset\tauup_0$. \end{lem}
\begin{proof} Indeed, $\tauup_0+(\mathfrak{t}_0\cap\mathfrak{s}_0^A)$
is a maximal torus in $\mathfrak{l}_{\mathfrak{q}}=\Lie(\mathbf{L}_{\mathfrak{q}})$.
 The homotopic characteristic of 
$\mathbf{L}_{\mathfrak{q}}/\mathbf{M}_0$ is the difference 
between the dimensions of
a maximal torus of $\mathfrak{l}_{\mathfrak{q}}$ 
and of a maximal torus of $\mathfrak{m}_0$.
This gives \eqref{eq:810}. The last claim follows because the Euler-Poincar\'e 
characteristic of $\mathbf{L}_{\mathfrak{q}}/\mathbf{M}_0$ is positive if and only
if its homotopic characteristic is $0$.
\end{proof}
Let $L_{\mathfrak{q}}$ 
be the $\mathbf{K}_0$-homogeneous $CR$ manifold with isotropy
$\mathbf{L}_{\mathfrak{q}}$ and $CR$ algebra 
$(\kt_0,\mathfrak{m}+\nr(\mathfrak{q})+
\mathbb{C}\otimes_{\mathbb{R}}\mathfrak{s}_{\mathfrak{q}})$ at the base point. 
\par
By the preceding discussion we obtain:
\begin{prop} \label{prop:812}
Let $\mathfrak{q}\in\Par(\mathfrak{v},\mathbf{M}_0)$.
 Then \eqref{lft} factors through 
$\mathbf{K}_0$-equivariant $CR$-lifts:
\begin{equation*}
  \xymatrix{M_{\mathfrak{q}}\ar[rd]_{\phiup_{\mathfrak{q}}}\ar[r]^{\etaup_{\mathfrak{q}}}
&L_{\mathfrak{q}}\ar[d]^{\psiup_{\mathfrak{q}}}\\
&N_{\mathfrak{q}}.}
\end{equation*}
in which the fibers of $\etaup_{\mathfrak{q}}$ are homogeneous spaces of 
homotopic characteristic
$c_{\mathfrak{q}}$ and the fibers of 
$\psiup_{\mathfrak{q}}$ are real tori.
\end{prop} 
\begin{proof}
In fact the typical fiber of $\etaup_{\mathfrak{q}}$ is $\mathbf{L}_{\mathfrak{q}}/\mathbf{M}_0$,
and the fibers of $\psiup_{\mathfrak{q}}$ are quotients of tori.
 \end{proof} 

\subsection{Maximal $CR$ structures}
Proposition\,\ref{prop:62a} shows that the problem of finding the
$CR$-strengthenings $M_{\mathfrak{q}}$ of 
an $\mathfrak{n}$-reductive $M$ is equivalent to find
the minimal elements of $\Par(\mathfrak{v},\mathbf{M}_0)$. We shall prove
that there is essentially a unique 
$\mathfrak{n}$-reductive maximal $CR$ structure on $M$ (see Remark\,\ref{rm77a}). 
\begin{lem}\label{lem78}
 If $\mathfrak{q}_1,\mathfrak{q}_2\in\Par(\mathfrak{v},\mathbf{M}_0)$, then 
\begin{equation}
 \label{eq:97} \mathfrak{q}=\nr(\mathfrak{q}_1)+\mathfrak{q}_1\cap\mathfrak{q}_2
\end{equation}
is a parabolic subalgebra of $\kt$ which belongs to $\Par(\mathfrak{v},\mathbf{M}_0)$. 
\end{lem} 
\begin{proof}
 Fix a Cartan subalgebra $\mathfrak{t}$ of $\kt$ 
contained in $\mathfrak{q}_1\cap\mathfrak{q}_2$ and let $\mathcal{K}$ 
be the root system of $(\kt,\mathfrak{t})$. 
  Then there are $A_1,A_2\in\mathfrak{t}_{\mathbb{R}}=\{H\in\mathfrak{t}\mid
  \ad(H)\;\text{has real eigenvalues}\}$ such that 
\begin{equation*}
 \mathfrak{q}_i=\mathfrak{t}\oplus{\sum}_{\alpha(A_i)\geq{0}}\kt^{\alpha},
 \quad\text{for \; $i=1,2$.}
\end{equation*}
 If $\epsilon>0$ is so small that $|\alpha(A_2)|<\epsilon^{-1}|\alpha(A_1)|$ for
 all roots $\alpha$ with $\alpha(A_1)\neq{0}$, then \begin{equation*}
 \mathfrak{q}=\mathfrak{t}\oplus{\sum}_{\alpha(A_1+\epsilon A_2)\geq{0}}\kt^{\alpha}.
 \end{equation*}
 Indeed $\mathfrak{q}$ contains $\mathfrak{t}$  and $\kt^{\alpha}\subset\mathfrak{q}$
 if and only if either $\alpha(A_1)>0$, or $\alpha(A_1)=0$ and $\alpha(A_2)\geq{0}$.
 This proves that $\mathfrak{q}$ is parabolic in $\kt$. 
 It belongs to $\Par(\mathfrak{v},\mathbf{M}_0)$ because,
by Lemma\,\ref{lem64},  
the right hand side of \eqref{eq:97} is a sum of $\Ad_{\kt}(\mathbf{M}_0)$-invariant
subalgebras, and $\mathfrak{v}\subset\mathfrak{q}_1\cap\mathfrak{q}_2$, \;
$\nr(\mathfrak{v})\cap\Li(\mathfrak{q})\subset
\nr(\mathfrak{v})\cap\Li(\mathfrak{q}_1)=\{0\}$.
\end{proof}
\begin{prop}\label{prop:84x} \begin{enumerate}
\item\label{11}
  If $\mathfrak{q}_1,\mathfrak{q}_2$ are minimal in
$\Par(\mathfrak{v},\mathbf{M}_0)$, then 
$\Li(\mathfrak{q}_1)$ and $\Li(\mathfrak{q}_2)$ are 
$\Ad_{\kt}(\mathbf{K}_0)$-conjugate.
\item\label{12}
For all minimal elements $\mathfrak{q}$
of $\Par(\mathfrak{v},\mathbf{M}_0)$ we have
$\nr(\mathfrak{v})\subset\nr(\mathfrak{q})$.
\end{enumerate}
\end{prop}
\begin{proof} Let $\mathfrak{q}_1,\mathfrak{q}_2$ be two elements of
$\Par(\mathfrak{v},\mathbf{M}_0)$. By Lemma\,\ref{lem78},
$  \mathfrak{q}=\nr(\mathfrak{q}_1)+
{\mathfrak{q}_1\!\cap\!\mathfrak{q}_2}$
is parabolic and belongs to $\Par(\mathfrak{v},\mathbf{M}_0)$. 
If $\mathfrak{q}_1$ is minimal, then $\mathfrak{q}=\mathfrak{q}_1$. 
Hence $\mathfrak{q}_1\cap\mathfrak{q}_2$ contains a reductive Levi-factor
$\upsilonup$
of $\mathfrak{q}_1$, which is also a maximal reductive
subalgebra of $\mathfrak{q}_1\cap\mathfrak{q}_2$.  
If also $\mathfrak{q}_2$ is minimal,
then, by the same argument, $\upsilonup$ is also a reductive Levi factor
of $\mathfrak{q}_2$. Thus we have $\Li(\mathfrak{q}_1)\sim\upsilonup\sim
\Li(\mathfrak{q}_2)$ with respect to $\Ind(\kt)$. Since both
$\Li(\mathfrak{q}_1)$ and $\Li(\mathfrak{q}_2)$ are complexifications of
subalgebras of $\kt_0$, they are also $\Ind(\kt_0)$-conjugate, and hence
$\Ad_{\kt}(\mathbf{K}_0)$-conjugate. This proves \eqref{11}.
\par
Let $\mathfrak{q}_1$ be any element of 
$\Par(\mathfrak{v},\mathbf{M}_0)$ and  $\mathfrak{q}_2$ be the
parabolic regularization of $\mathfrak{v}$. We recall that 
$\mathfrak{q}_2\in\Par(\mathfrak{v},\mathbf{M}_0)$ and
$\nr(\mathfrak{v})\subset\mathfrak{q}_2$. Then $  \mathfrak{q}=\nr(\mathfrak{q}_1)+
{\mathfrak{q}_1\!\cap\!\mathfrak{q}_2}$
 belongs to $\Par(\mathfrak{v},\mathbf{M}_0)$ and
$(\nr(\mathfrak{q}_1)\! +\!
\nr(\mathfrak{q}_2)\!\cap\!\mathfrak{q}_1)\subset\nr(\mathfrak{q})$
because it is an ideal of $\mathfrak{q}$ consisting of nilpotent elements.
If $\mathfrak{q}_1$ is minimal, we have $\mathfrak{q}=\mathfrak{q}_1$
and hence $\nr(\mathfrak{q})=\nr(\mathfrak{q}_1)$ yields 
$\nr(\mathfrak{v})\subset\nr(\mathfrak{q}_2)\cap\mathfrak{q}_1\subset
\nr(\mathfrak{q}_1)$, proving \eqref{12}.
\end{proof}
\begin{rmk}\label{rm77a}
If $\mathfrak{q}_1,\mathfrak{q}_2\in\Par$ have 
$\Ind(\kt)$-conjugate
reductive Levi factors, then there is an element
$g\in\mathbf{K}_0$ such that $\mathfrak{g}_2=\Ad_{\kt}(g)(\mathfrak{g}_1)$. \par
This follows from the observations that:
\begin{list}{-}{}
\item  $\Ad_{\kt}(\mathbf{K}_0)$ is transitive on
the $\Ind(\kt)$-orbit of $\mathfrak{q}_1$;
\item hence, if $\Li(\mathfrak{q}_1)$ and $\Li(\mathfrak{q}_2)$ are
$\Ind(\kt)$-conjugate, then there is $g\in\mathbf{K}_0$ such that
$\Ad_{\kt}(g)(\Li(\mathfrak{q}_1))=\Li(\mathfrak{q}_2)$;
\item if $\Li(\mathfrak{q}_1)=\Li(\mathfrak{q}_2)$, then there is an
element of the normalizer in $\mathbf{K}_0$ of a maximal torus
$\mathfrak{t}_0$ of $\Li(\mathfrak{q}_1)\cap\Li(\mathfrak{q}_2)\cap\kt_0$
that transforms $\mathfrak{q}_1$ into $\mathfrak{q}_2$. 
\end{list}
The $g$-translation on $N_{\mathfrak{q}_1}$ realizes a biholomorphism 
between $N_{\mathfrak{q}_1}$ and $N_{\mathfrak{q}_2}$. \par
If moreover $\mathfrak{q}_1,\mathfrak{q}_2\in\Par(\mathfrak{v},\mathbf{M}_0)$,
we have a commutative diagram
\begin{equation*}
  \begin{CD}
    M_{\mathfrak{q}_1}@>{g\cdot}>> M_{\mathfrak{q}_2}\\
@V{\piup_{\mathfrak{q}_1}}VV @VV{\piup_{\mathfrak{q}_2}}V \\
N_{\mathfrak{q}_1}@>>{g\cdot}> N_{\mathfrak{q}_2}
  \end{CD}
\end{equation*}
in which $\piup_{\mathfrak{q}_1}, g\cdot\piup_{\mathfrak{q}_1},
\piup_{\mathfrak{q}_2}, g^{-1}\cdot\piup_{\mathfrak{q}_2}$ are all
$CR$-lifts, but in which the 
left translation on $M$ in the top horizontal line is not $CR$ when 
$\mathfrak{q}_1\neq\mathfrak{q}_2$. 
\end{rmk}
\begin{cor}[maximal $CR$ structure]
 If $\mathfrak{q}$ is a minimal element of $\Par(\mathfrak{v},\mathbf{M}_0)$, 
then $M_{\mathfrak{q}}$ is maximal among the $\mathbf{K}_0$-equivariant
strengthening of the
$CR$ structure of $M$ for which there is a $CR$-lift
 $\phiup_{\mathfrak{q}}:M_{\mathfrak{q}}\to{N}_{\mathfrak{q}}$ 
from a complex flag manifold. \qed
\end{cor}
Let us use the notation of \S{\ref{sec:81x}}. 
In particular, $\mathfrak{t}_0$ is a maximal torus
of $\kt_0$ containing a maximal torus $\tauup_0$ of $\mathfrak{m}_0$ and
$\mathfrak{z}_0=\{A\in\mathfrak{t}_0\mid\Ad_{\kt}(g)(A)=A,\,\forall{g}\in\mathbf{M}_0\}$.
We have the following complement to Proposition\,\ref{prop:812}.
\begin{prop}\label{prop714} Let $A\in\mathfrak{z}_0$ and $\mathfrak{q}=\mathfrak{q}_A
\in\Par(\mathfrak{v},\mathbf{M}_0)$. \par
A necessary and sufficient condition for having $c_{\mathfrak{q}}=0$ is that
   $\mathfrak{q}$
is minimal in $\Par(\mathfrak{v},\mathbf{M}_0)$ and
\begin{equation}\label{eq:783}
\mathfrak{t}_0=\tauup_0+\mathfrak{z}_0.
\end{equation}
\par
A necessary condition in order that there exists $\mathfrak{q}
\in\Par(\mathfrak{v},\mathbf{M}_0)$
$c_{\mathfrak{q}}=0$
is that 
$(\kt_0,\mathfrak{m})$
be of type $\mathrm{I\! I}$. \par
When this condition is satisfied, and $\mathbf{M}_0$
is connected, then $c_{\mathfrak{q}}=0$ for all minimal $\mathfrak{q}_A
\in\Par(\mathfrak{m},\mathbf{M}_0)$. \qed
\end{prop}
\section{A generalization}\label{gen}
If we drop the assumption that $(\kt_0,\mathfrak{v})$ is $\mathfrak{n}$-reductive,
only part of Lemma\,\ref{lm71a} remains valid. Namely we obtain:
\begin{lem}\label{lm91}
 Let $\mathfrak{v}$ be any complex Lie algebra
and  $\mathfrak{q}\in\Par(\mathfrak{v},\mathbf{M}_0)$. 
Then
 $\mathfrak{v}_{\mathfrak{q}}=\mathfrak{v}+\nr(\mathfrak{q})$
is an $\Ad_{\kt}(\mathbf{M}_0)$-invariant Lie algebra,
with $\nr(\mathfrak{v}_{\mathfrak{q}})\! =\!\nr(\mathfrak{v})
\! +\!\nr(\mathfrak{q})$.
\end{lem} 
\begin{proof}
 If $\mathfrak{q}\in\Par(\mathfrak{v})$, then
$\mathfrak{v}\subset\mathfrak{q}$ and hence $[\mathfrak{v},\nr(\mathfrak{q})]
\subset [\mathfrak{q},\nr(\mathfrak{q})]\subset\nr(\mathfrak{q})$
shows that $\mathfrak{v}_{\mathfrak{q}}=
\mathfrak{v}+\nr(\mathfrak{q})$ is a Lie subalgebra of $\kt$,
containing $\nr(\mathfrak{q})$ as an ideal of nilpotent elements.
It follows that $\rad(\mathfrak{v}_{\mathfrak{q}})=\rad{\mathfrak{v}}+\nr(\mathfrak{q})$
and hence $\nr(\mathfrak{v}_{\mathfrak{q}})=\nr(\mathfrak{v})+\nr(\mathfrak{q})$.
Moreover, both $\mathfrak{v}$ and $\nr(\mathfrak{q})$ are
$\Ad_{\kt}(\mathbf{M}_0)$-invariant, and hence also $\mathfrak{v}_{\mathfrak{q}}$
is $\Ad_{\kt}(\mathbf{M}_0)$-invariant.
\end{proof}
\begin{lem}
  Let $\mathbf{L}_{\mathfrak{q}}$ be the analytic Lie subgroup 
of $\mathbf{K}_0$ with Lie algebra $\mathfrak{v}_{\mathfrak{q}}\cap\kt_0$.
Then
\begin{equation}
  \label{eq:8c}
 \mathbf{M}_{\mathfrak{q}}= 
\mathbf{M}_0\mathbf{L}_{\mathfrak{q}}=\{g_1g_2\mid g_1\in\mathbf{M}_0,\;
g_2\in\mathbf{L}_{\mathfrak{q}}\}
\end{equation}
is a closed subgroup of $\mathbf{N}_{\mathfrak{q}}$ 
with $\Lie(\mathbf{M}_{\mathfrak{q}})=\mathfrak{v}_{\mathfrak{q}}\cap\kt_0$
and
$\mathfrak{v}_{\mathfrak{q}}$ is an $\Ad_{\kt}(\mathbf{M}_{\mathfrak{q}})$-invariant
subalgebra of $\kt$.
\end{lem}
\begin{proof} Indeed, 
$\ad(\mathbf{M}_0)(\mathbf{L}_{\mathfrak{q}})=\mathbf{L}_{\mathfrak{q}}$
and hence $\mathbf{M}_q$ is a compact subgroup of $\mathbf{N}_{\mathfrak{q}}$.
Since $\mathfrak{m}_0\subset\mathfrak{v}_{\mathfrak{q}}\cap\kt_0$,
$\mathbf{M}_{\mathfrak{q}}$ and $\mathbf{L}_{\mathfrak{q}}$ have the same
Lie algebra $\mathfrak{v}_{\mathfrak{q}}\cap\kt_0$. Finally, both
$\mathbf{M}_0$ and $\mathbf{L}_{\mathfrak{q}}$ leave $\mathfrak{v}_{\mathfrak{q}}$
invariant, and thus the same is true for their product group
$\mathbf{M}_{\mathfrak{q}}$.
\end{proof}
Thus we obtain
\begin{prop}
For every $\mathfrak{q}\in\Par(\mathfrak{v},\mathbf{M}_0)$ there is
a $\mathbf{K}_0$-homogeneous $CR$ manifold $M_{\mathfrak{q}}$, with isotropy 
$\mathbf{M}_{\mathfrak{q}}$ and $CR$ algebra $(\kt_0,\mathfrak{v}_{\mathfrak{q}})$,
and we have a commutative diagram of $\mathbf{K}_0$-equivariant maps 
\begin{equation*}
  \xymatrix{M\ar[rr]\ar[drr]&& M_{\mathfrak{q}}\ar[d]\\
&& N_{\mathfrak{q}}.}
\end{equation*}
in which the right vertical
arrow is a $CR$-submersion onto a complex flag manifold.
\end{prop}
\section{Some examples}\label{exa}
We shall discuss in this section some examples of $CR$ manifolds which
are minimal orbits $M$ of real forms in complex flag manifolds
$F=\mathbf{G}/\mathbf{F}$, with $CR$ algebra $(\mathfrak{g}_0,\mathfrak{f})$
at the base point. 
\begin{exam} Our first example is one in which there is 
an adapted pair $(\thetaup,\mathfrak{h}_0)$ with 
a maximally
compact Cartan subalgebra $\mathfrak{h}_0$
of $\mathfrak{g}_0$, and therefore $M$ is a 
compact $\mathfrak{n}$-reductive compact manifold of type $\mathrm{I}$.
\par
Consider the real form $\mathfrak{g}_0=\mathfrak{sl}_2(\mathbb{H})$ of 
$\mathfrak{g}=\mathfrak{sl}_4(\mathbb{C})$ 
and take the minimal orbit $M$ corresponding to
the cross marked Satake diagram
\begin{equation*}
  \xymatrix@R=-.3pc{
\!\!\medbullet\!\!\ar@{-}[r]&\!\!\medcirc\!\!\ar@{-}[r]&\!\!\medbullet\!\!\ \\
\times && \times}
\end{equation*}
The maximal compact subalgebra of $\mathfrak{sl}_2(\mathbb{H})$ is
$\mathfrak{sp}_2$, with complexification  
${\kt}=\mathfrak{sp}_2(\mathbb{C})$.
We give below a matrix description of ${\kt}$ and
$\mathfrak{v}={\kt}\cap\mathfrak{f}$. Here and in the following 
we shall consistently use the convention
of labelling by $z_j$ the entries of the matrices which correspond 
to the directions in $T^{0,1}M$,
$\zeta_j\sim\bar{z}_j$ for their conjugate,
and by $w_j,\bar{w}_j\sim\eta_j,it_j$ the complex, or 
imaginary entries corresponding to directions transversal
to $HM$. We have
\begin{align*}
  {\kt}&:
  \begin{pmatrix}
\lambda_1&\zeta_1&\eta&\zeta_2\\
z_1&-\lambda_1&z_2&-w\\
w&\zeta_2&\lambda_2&\zeta_3\\
z_2&-\eta&z_3&-\lambda_2    
  \end{pmatrix}, & \mathfrak{v}={\kt}\cap\mathfrak{f}&:
  \begin{pmatrix}
    \lambda_1&\zeta_1&0&\zeta_2\\
0&-\lambda_1&0&0\\
0&\zeta_2&\lambda_2&\zeta_3\\
0&0&0&-\lambda_2
  \end{pmatrix}.
\end{align*}
If we denote by $\mathcal{K}=\{\pm{e}_1,\pm{e}_2,\pm{e}_1\pm{e}_2\}$
the root system of ${\kt}$ with respect to its Cartan subalgebra
\begin{equation*}
 \mathfrak{t}: \left(\begin{smallmatrix}
    \lambda_1\\
&-\lambda_1\\
&&\lambda_2\\
&&&-\lambda_2
  \end{smallmatrix}\right).
  \end{equation*}
Then $\mathfrak{v}$ is regular, described by
\begin{equation*}
  \mathfrak{v}={\kt}\cap\mathfrak{f}=\mathfrak{t}
\oplus{\sum}_{\alpha\in\mathcal{V}_n}{\kt}^{\alpha},\;\;
\text{with $\mathcal{V}_n=\{2e_1,2e_2,e_1+e_2\}$.}
\end{equation*}
The set $\mathcal{V}_n$ is the horocyclic set corresponding to the parabolic
set
\begin{equation*}
  \mathcal{Q}=\mathcal{Q}_n\cup\mathcal{Q}_r,\;\;\text{with
$\mathcal{Q}_n=\mathcal{V}_n$ and 
$\mathcal{Q}_r=\{\pm(e_1-e_2)\}$}
\end{equation*}
which defines the parabolic subalgebra 
\begin{equation*}
 \mathfrak{q}={\mathfrak{t}}\oplus{\sum}_{\alpha\in\mathcal{Q}}{\kt}^{\alpha}.
\end{equation*}
Correspondingly, the basis $N_{\mathfrak{q}}$ of the fibration
$M\xrightarrow{\;\piup\;}{N}_{\mathfrak{q}}$  
corresponding to the $CR$-submersion
$({\kt}_0,\mathfrak{v})\to({\kt}_0,\mathfrak{q})$
is the complex flag manifold defined by
the cross marked Satake diagram
\begin{equation*}
   \xymatrix@R=-.3pc{
\!\!\medbullet\!\!\ar@{<=}[r]&\!\!\medbullet\!\!\ \\
& \times}
\end{equation*}
This parabolic $\mathfrak{q}$ is 
the parabolic regularization of $\mathfrak{v}$. Indeed we have
\begin{align*}
 \nr( \mathfrak{v})&:\begin{pmatrix}
    0&\zeta_1&0&\zeta_2\\
0&0&0&0\\
0&\zeta_2&0&\zeta_3\\
0&0&0&0
  \end{pmatrix} &
\mathfrak{v}_1=N_{{\kt}}(\nr(\mathfrak{v}))&:\begin{pmatrix}
    \lambda_1&\zeta_1&-\eta&\zeta_2\\
0&-\lambda_1&0&-w\\
w&\zeta_2&\lambda_2&\zeta_3\\
0&\eta&0&-\lambda_2
  \end{pmatrix}=\mathfrak{v}_2=\mathfrak{q},
\end{align*}
the latest being the parabolic subalgebra corresponding to the Satake
diagram above. In this case $\nr(\mathfrak{v})=\nr(\mathfrak{v}_1)$, so 
that all terms in the sequence \eqref{eq:42b} are equal. 
\par
Note that ${\kt}_0\cap\mathfrak{q}\simeq\mathfrak{u}(2)$ and 
thus $N_{\mathfrak{q}}$ is the   Hermitian symmetric space 
$\mathbf{Sp}(2)/\mathbf{U}(2)$. The fiber is 
$\mathbf{SU}(2)/\mathbf{U(1)}\simeq{S}^2$,
realizing $M$ as an $S^2$-bundle over $N_{\mathfrak{q}}$.
\end{exam}
\begin{exam}\label{ex:92} Our next example is type $\mathrm{I\!{I}}$.\par
Consider the real form $\mathfrak{g}_0=\mathfrak{su}_{1,3}$ of 
$\mathfrak{g}=\mathfrak{sl}_4(\mathbb{C})$, take the minimal orbit 
$M$ corresponding to
the cross marked Satake diagram
\medskip
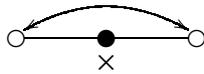
\begin{equation*}
  \xymatrix@R=-.3pc{\medcirc\!\!\ar@{-}[r]
\ar@{<->}@/^1pc/[rr]&\!\!\medbullet\!\!\ar@{-}[r]&\!\!\medcirc\\
&\times}
\end{equation*}
and denote by $(\mathfrak{g}_0,\mathfrak{f})$ the corresponding $CR$ algebra.
We write below
the traceless
matrix representation of the different Lie algebras 
involved in the construction of the
parabolic regularization of $\mathfrak{v}={\kt}\cap\mathfrak{f}$.
\begin{align*}
\mathfrak{g}_0\simeq  \mathfrak{su}_{1,3}&:
  \begin{pmatrix}
    \lambda&\zeta_1&\zeta_2&i\sigma\\
-\bar{z}_1&it_1&-\bar{z}_3&-\bar{\zeta}_1\\
z_2&z_3&it_2&-\bar{\zeta}_2\\
is&{z}_1&-\bar{z}_2&-\bar{\lambda}
  \end{pmatrix}, &
{\kt}_0\simeq\mathfrak{u}(3)&:
\begin{pmatrix}
  it_0&z_1&-\bar{z}_2&is\\
-\bar{z}_1&it_1&-\bar{z}_3&-\bar{z}_1\\
z_2&z_3&it_2&z_2\\
is&z_1&-\bar{z}_2&it_0
\end{pmatrix},\\
{\kt}\simeq\mathbb{C}\ltimes\mathfrak{sl}_3(\mathbb{C})&:
\begin{pmatrix}
 \lambda_0&z_1&\zeta_2&w\\
\zeta_1&\lambda_1&\zeta_3&\zeta_1\\
z_2&z_3&\lambda_2&z_2\\
w&z_1&\zeta_2&\lambda_0
\end{pmatrix}, &
\mathfrak{v}={\kt}\cap\mathfrak{f}&:
\begin{pmatrix}
  \lambda_0&0&\zeta_2&0\\
\zeta_1&\lambda_1&\zeta_3&\zeta_1\\
0&0&\lambda_2&0\\
0&0&\zeta_2&\lambda_0
\end{pmatrix},\\
\nr(\mathfrak{v})&:
\begin{pmatrix}
0&0&\zeta_2&0\\
\zeta_1&0&\zeta_3&\zeta_1\\
0&0&0&0\\
0&0&\zeta_2&0 
\end{pmatrix}, &
\mathfrak{v}_1&:
\begin{pmatrix}
  \lambda_0&0&\zeta_2&w\\
\zeta_1&\lambda_1&\zeta_3&\zeta_1\\
0&0&\lambda_2&0\\
w&0&\zeta_2&\lambda_0
\end{pmatrix}=\mathfrak{v}_2=\mathfrak{q}
\end{align*}
Then the parabolic $\mathfrak{q}$ corresponds to the complete complex flag 
manifold $N_{\mathfrak{q}}$ of
$\mathbf{SL}_3(\mathbb{C})$, associated to the diagram
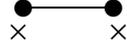
\begin{equation*}
  \xymatrix@R=-.3pc{\medbullet\!\!\ar@{-}[r]&\!\!\medbullet\\
\times&\times}
\end{equation*}
In this case this $\mathfrak{q}\in\Par(\mathfrak{v})$ is also minimal and
the fibration $M\xrightarrow{\;\piup\;}N_{\mathfrak{q}}$ associated 
to the submersion of $CR$ algebras
$({\kt}_0,\mathfrak{v})\to({\kt}_0,\mathfrak{v})$ is a circle bundle.
\end{exam}
\begin{exam}\label{ex:93} 
Consider again the real form $\mathfrak{g}_0=\mathfrak{su}_{1,3}$ of 
$\mathfrak{g}=\mathfrak{sl}_4(\mathbb{C})$
and take this time the minimal orbit $M$ corresponding to the Satake diagram   
\medskip
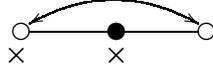
\begin{equation*}
  \xymatrix@R=-.3pc{\medcirc\!\!\ar@{-}[r]
\ar@{<->}@/^1pc/[rr]&\!\!\medbullet\!\!\ar@{-}[r]&\!\!\medcirc\\
\times&\times}
\end{equation*}
The $\mathbf{G}_0$-equivariant $CR$ structure on $M$ defined by the $CR$ algebra
$(\mathfrak{g}_0,\mathfrak{f})$, in this case, is not maximal 
(see e.g. \cite[Prop.5.9]{AMN06b}).
We give below the traceless matrix representation of the different 
Lie algebras involved
in the parabolic regularization of $\mathfrak{v}={\kt}\cap\mathfrak{f}$:
\begin{align*}
\mathfrak{g}_0\simeq  \mathfrak{su}_{1,3}&:
  \begin{pmatrix}
    \lambda&\zeta_1&\zeta_2&i\sigma\\
w&it_1&-\bar{z}_1&-\bar{\zeta}_1\\
z_2&z_1&it_2&-\bar{\zeta}_2\\
is&-\bar{w}&-\bar{z}_2&-\bar{\lambda}
  \end{pmatrix}, &
{\kt}_0\simeq\mathfrak{u}(3)&:
\begin{pmatrix}
  it_0&-\bar{w}&-\bar{z}_2&is\\
w&it_1&-\bar{z}_1&-w\\
z_2&z_1&it_2&z_2\\
is&\bar{w}&-\bar{z}_2&it_0
\end{pmatrix},\\
{\kt}\simeq\mathbb{C}\ltimes\mathfrak{sl}_3(\mathbb{C})&:
\begin{pmatrix}
 \lambda_0&\eta&\zeta_2&\mu\\
w&\lambda_1&\zeta_1&-w\\
z_2&z_1&\lambda_2&z_2\\
\mu&\eta&\zeta_2&\lambda_0
\end{pmatrix}, &
\mathfrak{v}={\kt}\cap\mathfrak{f}&:
\begin{pmatrix}
  \lambda_0&0&\zeta_2&0\\
0&\lambda_1&\zeta_1&0\\
0&0&\lambda_2&0\\
0&0&\zeta_2&\lambda_0
\end{pmatrix},\\
\nr(\mathfrak{v})&:
\begin{pmatrix}
0&0&\zeta_2&0\\
0&0&\zeta_1&0\\
0&0&0&0\\
0&0&\zeta_2&0 
\end{pmatrix}, &
\mathfrak{v}_1&:
\begin{pmatrix}
  \lambda_0&\eta&\zeta_2&\mu\\
0&\lambda_1&\zeta_1&0\\
0&0&\lambda_2&0\\
\mu&\eta&\zeta_2&\lambda_0
\end{pmatrix}=\mathfrak{v}_2=\mathfrak{q}
\end{align*}
Again, the parabolic $\mathfrak{q}$ is minimal in $\Par(\mathfrak{v})$ and
corresponds to the complete 
complex flag manifold of
$\mathbf{SL}_3(\mathbb{C})$, associated to the diagram
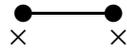
\begin{equation*}
  \xymatrix@R=-.3pc{\medbullet\!\!\ar@{-}[r]&\!\!\medbullet\\
\times&\times}\end{equation*}
The $CR$ algebra $(\kt_0,\mathfrak{v}_{\mathfrak{q}})$, with 
$\mathfrak{v}_{\mathfrak{q}}=
\nr(\mathfrak{q})+\mathfrak{v}\cap\bar{\mathfrak{v}}$,
defines a $CR$ manifold $M_{\mathfrak{q}}$ that is obtained from $M$ by a 
$\mathbf{K}_0$-equivariant strengthening of the
$CR$ structure. The corresponding 
$CR$ fibration $M_{\mathfrak{q}}
\xrightarrow{\;\piup_{\mathfrak{q}}\;}N_{\mathfrak{q}}$ is a circle bundle.
Note that the $CR$ structure of $M_{\mathfrak{q}}$ 
is not $\mathbf{G}_0$-equivariant and 
thus it is not
equivalent to the $CR$ structure of the minimal orbit of 
Example\,\ref{ex:92}.\par
\smallskip
Denote by $M_1$, $M_2$ the $CR$ manifolds $M$ and
$M_{\mathfrak{q}}$ of Examples\,\ref{ex:92},\ref{ex:93}, respectively. 
They correspond to the datum of two different $CR$ structures 
on the same underlying
smooth manifold. Let us describe geometrically their projections into the 
complete complex flag
manifold $N_{\mathfrak{q}}$ of $\mathbf{SL}_3(\mathbb{C})$.
\par Fix 
a Hermitian symmetric form $\sigmaup$ in $\mathbb{C}^4$ of signature $(1,3)$ 
and two
$\sigmaup$-orthogonal subspaces $V_1$, $V_3$ on which $\sigmaup$ is definite.
From the $CR$ manifold $M_1$ of $\sigmaup$-isotropic two-planes 
$L_2$ of $\mathbb{C}^4$ the $CR$-submersion  $M_1\to{N}_{\mathfrak{q}}$ 
is defined by 
\begin{equation*}
 L_2\longrightarrow L_2\cap{V}_3\subset (V_1+L_2)\cap{V}_3\subset{V}_3,
\end{equation*}
while the $CR$-submersion $M_2\to{N}_{\mathfrak{q}}$ from the 
$CR$ manifold $M_2$ of $L_1\subset{L}_2$ 
with a $\sigmaup$-isotropic line $L_1$ and a 
$\sigmaup$-isotropic
$2$-plane $L_2$ containing $L_1$, is given by
\begin{equation*}
L_1\subset{L}_2\longrightarrow (L_1+V_1)\cap{V}_3
\subset (L_2+V_1)\cap{V}_3\subset{V}_3.
\end{equation*}
\end{exam}
\begin{exam} We consider again Example\;\ref{ex:b}. We have
  \begin{align*}
    \mathfrak{v}={\kt}\cap\mathfrak{f}&:
    \begin{pmatrix}
      \lambda_2&\zeta_1&0&0&0\\
0&\lambda_1&0&0&0\\
0&\zeta_2&\lambda_0&\zeta_2&0\\
0&0&0&\lambda_1&0\\
0&0&0&\zeta_1&\lambda_2
    \end{pmatrix}, &
\rn(\mathfrak{v})&:\begin{pmatrix}
      0&\zeta_1&0&0&0\\
0&0&0&0&0\\
0&\zeta_2&0&\zeta_2&0\\
0&0&0&0&0\\
0&0&0&\zeta_1&0
    \end{pmatrix} \\
\mathfrak{v}_1&:\begin{pmatrix}
      \lambda_2&\zeta_1&0&\eta_2&\mu\\
0&\lambda_1&0&\mu&0\\
w_1&\zeta_2&\lambda_0&\zeta_2&w_1\\
0&\mu&0&\lambda_1&0\\
\mu&\eta_2&0&\zeta_1&\lambda_2
    \end{pmatrix}, &
\rn(\mathfrak{v}_1)&:\begin{pmatrix}
      0&\zeta_1&0&\eta_2&0\\
0&0&0&0&0\\
w_1&\zeta_2&0&\zeta_2&w_1\\
0&0&0&0&0\\
0&\eta_2&0&\zeta_1&0
    \end{pmatrix}\\
\mathfrak{v}_2&:
\begin{pmatrix}
  \lambda_2&\zeta_1&0&\eta_2&\mu_2\\
0&\lambda_1&0&\mu_1&0\\
w_1&\zeta_2&\lambda_0&\zeta_2&w_1\\
0&\mu_1&0&\lambda_1&0\\
\mu_2&\eta_2&0&\zeta_1&\lambda_2
\end{pmatrix}, &
\rn(\mathfrak{v}_2)&=\rn(\mathfrak{v}_1)\;\Longrightarrow\;
\mathfrak{v}_2=\mathfrak{v}_3=\mathfrak{q}.
\end{align*}
Therefore the parabolic regularization 
yields a $CR$ map
$M_3^{2,6}\to{N}$ of $M_3$ onto a compact complex manifold 
$N_{\mathfrak{q}}=N^{4,0}_{\mathfrak{q}}$ of complex dimension $4$. \par
The manifold $N_{\mathfrak{q}}$ is the complete flag of 
$\mathbf{S}(\mathbf{U}(2)\times\mathbf{U}(3))$.
The fibration can be described explicitly by fixing in $\mathbb{C}^5$ 
a Hermitian symmetric form $\sigmaup$ of signature $(2,3)$ and two
$\sigmaup$-orthogonal subspaces $V_2,V_3$ on which the form is definite.
We associate to the flag $(\ell_1,\ell_3)\in{M}_3$ the point
$(\ell_1+V_3)\cap{V}_2$ of $\mathbb{P}(V_2)\simeq\mathbb{CP}^1$ and
the flag $(\ell_3\cap{V}_3,(l_3\cap V_3 + (l_1+V_2)\cap V_3))
\in\mathfrak{Gr}_{1,2}(V_3)$.
The manifold $N_{\mathfrak{q}}$ 
is described by the cross marked Satake diagram
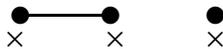
\begin{equation*}
  \xymatrix@R=-.3pc{\medbullet\!\!\ar@{-}[r]&\!\!\medbullet &\medbullet\\
\times&\times&\times}
\end{equation*} \par
We note that in this case the corresponding
$\mathbf{S}(\mathbf{U}(2)\times\mathbf{U}(3))$-equivariant 
$CR$ map $\piup:M\to{N}_{\mathfrak{q}}$ 
is a smooth submersion, and a $CR$-deployment, but not a $CR$ 
submersion. The space 
$E_q{N}_{\mathfrak{q}}={\sum}_{p\in\piup^{-1}(q)}d\piup(T^{0,1}_pM)$ 
is a hypersurface in
$T^{0,1}_qN_{\mathfrak{q}}$ for all $q\in{N}_{\mathfrak{q}}$ and a non formally
integrable complex distribution. This example explains the importance
of considering commutators in Definition\,\ref{df24}.
\end{exam}

\bibliographystyle{amsplain}

\renewcommand{\MR}[1]{}
\bibliography{homog}

\end{document}